\documentclass[a4paper]{amsart}
\usepackage[francais,english]{babel}
\usepackage{amsmath,amssymb,amsthm}
\usepackage{amsfonts}
\usepackage[latin1]{inputenc}
\usepackage[T1]{fontenc}

\usepackage{cite}
\usepackage{pifont}
\usepackage{euscript}
\usepackage{footnpag}
\usepackage{verbatim}
\RequirePackage{mathrsfs}

\def\noignorespaces{{\catcode`\ =11\ }}

\newcounter{constante}
\DeclareMathOperator{\Hom}{Hom}
\DeclareMathOperator{\GL}{GL}

\DeclareMathOperator{\spec}{Spec}
 
\DeclareMathOperator{\rg}{rg} 

\DeclareMathOperator{\Proj}{Proj}
\DeclareMathOperator{\vol}{vol}
\DeclareMathOperator{\covol}{covol}
\DeclareMathOperator{\gal}{Gal}

\DeclareMathOperator{\vr}{vr}

\theoremstyle{plain}
\newtheorem{theo}{Th\'eor\`eme}[section]

\newtheorem{lemma}[theo]{Lemme}
\newtheorem{prop}[theo]{Proposition}
\newtheorem{coro}[theo]{Corollaire}
\newtheorem{theoRumely}{Th\'eor\`eme~\cite{Rumelyetal}}

\newtheorem{propAnnexe}{Proposition}

\theoremstyle{definition}
\newtheorem{defi}[theo]{D\'efinition}
\newtheorem{defithm}[theo]{D{\'e}finition-th{\'e}or{\`e}me}
\newtheorem{rema}[theo]{Remarque}
\newtheorem{remas}[theo]{Remarques}
\newtheorem{proprietes}[theo]{Propri{\'e}t{\'e}s}

\newtheorem{exems}[theo]{Exemples}

\begin{document}
\date{Jeudi 1er f\'evrier 2007}
\title[Pentes des fibr{\'e}s vectoriels ad{\'e}liques]{Pentes des fibr{\'e}s vectoriels ad{\'e}liques sur un corps global}

\author{{\'E}ric Gaudron}

\selectlanguage{francais}

\subjclass[2000]{11H06 (11G50, 11R56, 14G99)}
\footnotetext{\textbf{Mots clefs}: Fibr\'e vectoriel ad\'elique, degr\'e ad\'elique, fibr\'e de John, fibr\'e de L{\"o}wner, quotient volumique, distance de Banach-Mazur, in\'egalit\'es de pentes ad\'eliques, minima successifs ad\'eliques, pente maximale.}

\begin{abstract}
Dans les ann\'ees $90$, J.-B.~Bost a d\'evelopp\'e tout un formalisme des pentes des fibr\'es vectoriels hermitiens sur l'anneau des entiers d'un corps de nombres. Au cours de ses recherches, une nouvelle m\'ethode d'approximation diophantienne --- dite \emph{m\'ethode des pentes} --- a \'et\'e \'elabor\'ee. Cet article propose une g\'en\'eralisation de ces travaux \`a une classe plus large de fibr\'es vectoriels, dits ad\'eliques, d\'efinis sur un corps global. Ces fibr\'es poss\`edent aux places archim\'ediennes des normes qui ne sont plus n\'ecessairement hermitiennes. Nous examinons \'egalement le lien avec la th\'eorie des minima successifs ad\'eliques. Pour parvenir \`a ces r\'esultats, nous avons recours \`a plusieurs concepts de g\'eom\'etrie des espaces de Banach de dimension finie.
 \\[.2cm] \selectlanguage{english}\textsc{Abstract}. At the end of the twentieth century, J.-B.~Bost developped a slope theory of hermitian vector bundles over number fields. A new method of diophantine approximation, the so-called \emph{slope method}, has emerged from his research. Our article proposes a generalisation to adelic vector bundles over global fields. The norms at the archimedean places are no longer supposed to be hermitian. The link with adelic successive minima is also mentioned. To get these results, we use several concepts from the geometry of finite dimensional Banach spaces.     
\end{abstract}

\selectlanguage{francais}
\maketitle
\newpage
                
\section{Introduction}
Ce texte d{\'e}crit une g{\'e}n{\'e}ralisation de la th{\'e}orie des pentes des fibr{\'e}s vectoriels hermitiens sur l'anneau des entiers d'un corps de nombres aux fibr{\'e}s vectoriels ad{\'e}liques sur un corps global.\par \'Etant donn{\'e} un corps de nombres $k$ d'anneau des entiers $\mathcal{O}_{k}$, un fibr{\'e} vectoriel hermitien $\overline{E}$ sur $\spec\mathcal{O}_{k}$ est la donn{\'e}e d'un $\mathcal{O}_{k}$-module projectif de type fini $E$ et, pour toute place archim{\'e}dienne $v$ de $k$, d'une norme $\Vert.\Vert_{v}$ sur l'espace vectoriel $E_{v}:=E\otimes_{k}k_{v}$ ($k_{v}=\mathbf{R}$ ou $\mathbf{C}$) qui est euclidienne si $v$ est une place r{\'e}elle et qui est hermitienne, invariante par conjugaison complexe, si $v$ est une place complexe. \`A une telle donn{\'e}e, l'on peut associer un nombre r{\'e}el appel{\'e} \emph{degr{\'e} d'Arakelov} (et not{\'e} $\widehat{\deg}_{\mathrm{n}}\overline{E}$), qui, au signe pr\`es, mesure une \og hauteur de $\overline{E}$\fg. Ces nombres jouent un r{\^o}le important en g\'eom\'etrie d'Arakelov. Ce sont les {\'e}l{\'e}ments primitifs {\`a} partir desquels se b{\^a}tissent d'autres invariants associ\'es \`a $\overline{E}$. Il s'agit notamment des pentes du graphe du polygone qui d{\'e}limite sup{\'e}rieurement l'enveloppe convexe des couples de nombres r{\'e}els $(\rg F,\widehat{\deg}_{\mathrm{n}}\overline{F})$ o\`u $\overline{F}$ parcourt les sous-fibr{\'e}s de $\overline{E}$. Les images de ces pentes par la fonction $x\mapsto e^{-x}$ se comparent aux minima successifs ad{\'e}liques de $\overline{E}$, d{\'e}finis par Bombieri \& Vaaler~\cite{BombieriVaaler}. \`A l'occasion de cours de $3^{\text{\`eme}}$q cycle donn{\'e}s {\`a} l'Institut Henri Poincar{\'e} (Paris) en $1997$ et $1999$, J.-B.~Bost a effectu{\'e} une {\'e}tude syst{\'e}matique des propri{\'e}t{\'e}s de ces nombres, {\'e}laborant ainsi une v{\'e}ritable \emph{th{\'e}orie des pentes}. L'objectif poursuivi initialement {\'e}tait de reformuler sous forme plus g{\'e}om{\'e}trique et intrins{\`e}que la d{\'e}monstration du th{\'e}or{\`e}me des p{\'e}riodes de D.~Masser \& G.~W{\"u}stholz (voir~\cite{bost2}). Les notes de ces cours n'ont pas {\'e}t{\'e} publi{\'e}es. N{\'e}anmoins plusieurs fragments se trouvent dans les articles~\cite{bost2,bost6,aclbourbaki,graftieaux1,graftieaux2,viada}. De l'article fondateur~\cite{bost2} est issue une m{\'e}thode --- dite \emph{m{\'e}thode des pentes} --- destin{\'e}e {\`a} prouver des {\'e}nonc{\'e}s de transcendance et d'approximation diophantienne. Elle se rapproche de la m{\'e}thode des d{\'e}terminants d'interpolation de M.~Laurent. La grande force de la m{\'e}thode des pentes est de s'adapter naturellement {\`a} un probl{\`e}me de nature g{\'e}om{\'e}trique. Cette caract\'eristique renforce la clart\'e de l'argumentation en s\'eparant distinctement les contributions, tout en faisant ressortir les invariants naturels des objets g\'eom\'etriques. Par exemple, elle a permis de mettre en lumi{\`e}re et de d{\'e}montrer un crit{\`e}re d'alg{\'e}bricit{\'e} de feuilles formelles (voir~\cite{bost6}). Nous l'avons {\'e}galement utilis{\'e}e pour fournir des minorations de formes lin{\'e}aires de logarithmes de vari{\'e}t{\'e}s ab{\'e}liennes principalement polaris{\'e}es, minorations qui sont totalement explicites en la dimension et la hauteur de Faltings de la vari{\'e}t{\'e} (voir~\cite{artepredeux}).\par \`A l'usage, il arrive parfois que le cadre des fibr{\'e}s vectoriels hermitiens sur $\spec\mathcal{O}_{k}$ dans lequel s'applique la m{\'e}thode des pentes s'av{\`e}re trop rigide. \`A la suite des travaux de S.~Zhang~\cite{Zhang2}, V.~Maillot~\cite{VMaillotsmf} ou bien encore de R.~Rumely \textit{et al.}~\cite{Rumelyetal}, il est apparu que, si l'on souhaite construire une \og hauteur canonique\fg\ sur les cycles d'une vari{\'e}t{\'e} projective $X$ munie d'un fibr{\'e} en droites ample $M$, les m{\'e}triques que l'on doit mettre sur $M$ ne sont pas en g{\'e}n{\'e}ral hermitiennes mais seulement continues. Cela emp{\^e}che alors de mettre en {\oe}uvre la m{\'e}thode des pentes telle quelle, comme on aimerait le faire par exemple avec l'espace vectoriel des sections globales $\mathrm{H}^{0}(X,M)$.\par Ces observations nous ont amen{\'e} {\`a} examiner \`a nouveau le formalisme des pentes pour des fibr\'es vectoriels munis d'une structure plus souple que celle des fibr{\'e}s vectoriels hermitiens sur $\spec\mathcal{O}_{k}$. Avant de pr{\'e}senter plus en d{\'e}tail les r{\'e}sultats de cet article, je tiens {\`a} souligner que la plupart d'entre eux --- preuves comprises --- proviennent des cours de J.-B.~Bost mentionn{\'e}s ci-dessus, au moins en ce qui concerne le cas hermitien.\par Dor{\'e}navant, nous consid{\'e}rons un corps global $k$ (corps de nombres ou corps de fonctions). En nous inspirant de~\cite{Rumelyetal}, nous d{\'e}finissons la notion de \emph{fibr{\'e} vectoriel ad{\'e}lique sur} $\spec k$ qui g{\'e}n{\'e}ralise celle de fibr{\'e} hermitien sur $\spec\mathcal{O}_{k}$. Un tel objet est la donn{\'e}e d'un $k$-espace vectoriel $E$ de dimension finie $n$, muni d'une $k$-base $\mathbf{e}$ et, pour chaque place $v$ de $k$, d'une norme $\Vert.\Vert_{v}$ sur $E\otimes_{k}\mathbf{C}_{v}$, invariante sous l'action des automorphismes continus de $\gal(\mathbf{C}_{v}/k_{v})$. On le note $\overline{E}=(E,(\Vert.\Vert_{v})_{v})$. Si $v$ est ultram{\'e}trique, la norme $\Vert.\Vert_{v}$ doit v{\'e}rifier l'in{\'e}galit{\'e} ultram{\'e}trique usuelle et, sauf pour un nombre fini de $v$, elle est {\'e}gale {\`a} la norme du sup sur $E\otimes_{k}\mathbf{C}_{v}$, cet espace \'etant identifi{\'e} {\`a} $\mathbf{C}_{v}^{n}$ au moyen de la base $\mathbf{e}$. Le trait marquant de ces fibr{\'e}s ad{\'e}liques est que les normes aux places archim{\'e}diennes ne sont plus n{\'e}cessairement hermitiennes. La collection des normes de $\overline{E}$ dote l'espace ad{\'e}lique $E\otimes_{k}k_{\mathbf{A}}$ d'une boule unit{\'e}, dont le logarithme du volume jouera le r{\^o}le de degr{\'e} d'Arakelov de $\overline{E}$ ({\`a} une constante pr{\`e}s). Il sera appel{\'e} \emph{degr{\'e} ad{\'e}lique} de $\overline{E}$ dans la suite. Une fois ces d{\'e}finitions fix{\'e}es, nous {\'e}tudions les propri{\'e}t{\'e}s de ce degr{\'e} vis-{\`a}-vis des op{\'e}rations usuelles que l'on peut effectuer sur l'ensemble des fibr{\'e}s vectoriels ad{\'e}liques (extension du corps de base, somme directe, etc.). Puis nous construisons l'analogue du \og polygone canonique\fg\ {\`a} partir duquel s'obtiennent les $n$ pentes de $E$. Le reste de l'article s'attache alors {\`a} {\'e}tudier certaines propri{\'e}t{\'e}s de ces pentes, en particulier leur comportement par transformation lin{\'e}aire, c.-{\`a}-d. lorsque deux fibr{\'e}s $\overline{E}$ et $\overline{F}$ ont leurs espaces vectoriels sous-jacent reli{\'e}s par une application lin{\'e}aire. Ceci donne naissance \`a plusieurs in{\'e}galit{\'e}s --- dites \emph{in{\'e}galit{\'e}s de pentes} --- dont l'une est au c{\oe}ur de la m{\'e}thode des pentes, \'evoqu{\'e}e au d{\'e}but de cette introduction. \par La plupart des r{\'e}sultats que nous obtenons sont bas{\'e}s sur le m{\^e}me sch{\'e}ma de preuve. On commence par s'int{\'e}resser au cas hermitien, c.-{\`a}-d. au cas d'un fibr{\'e} vectoriel ad{\'e}lique dont les normes aux places archim{\'e}diennes sont hermitiennes. Les d{\'e}monstrations sont alors tr{\`e}s proches de celles d{\'e}j{\`a} connues pour les fibr{\'e}s vectoriels hermitiens sur $\spec\mathcal{O}_{k}$. Le passage au cas g{\'e}n{\'e}ral s'effectue au moyen d'une comparaison entre une norme quelconque $\Vert.\Vert$ sur $\mathbf{R}^{n}$ et une norme euclidienne, en faisant intervenir la distance --- dite de Banach-Mazur --- entre $(\mathbf{R}^{n},\Vert.\Vert)$ et l'espace euclidien usuel $\ell^{2}_{n}$. Cette d\'emarche est fr\'equente dans l'\'etude de la g{\'e}om{\'e}trie des espaces de Banach de dimension finie (g{\'e}om{\'e}trie de Minkowski). Toutefois, un certain nombre d'adaptations et de reformulations dans le cadre ad{\'e}lique ont {\'e}t{\'e} n{\'e}cessaires. Par exemple, nous d{\'e}finissons les fibr{\'e}s vectoriels ad{\'e}liques de John et L{\"o}wner associ{\'e}s {\`a} un fibr{\'e} vectoriel ad{\'e}lique $\overline{E}$ sur $\spec k$, qui sont des fibr{\'e}s vectoriels hermitiens qui encadrent au mieux $\overline{E}$ (en termes de volumes de boules unit{\'e}s). Ces fibr\'es fournissent des formules exactes pour le degr\'e ad\'elique de $\overline{E}$. Et cela conduit de temps en temps \`a des r\'esultats plus fins, qui font intervenir le \og quotient volumique ad\'elique\fg\ de $\overline{E}$, qui est un nombre r\'eel construit \`a l'aide du fibr\'e de John de $\overline{E}$, au lieu de la distance de Banach-Mazur ad\'elique. Un aspect int{\'e}ressant de cette approche est que les termes d'erreurs induits sont tr{\`e}s bien contr{\^o}l{\'e}s. Ils sont born{\'e}s par une fonction explicite de la dimension de $\overline{E}$ et du degr{\'e} de $k$. De plus, ils disparaissent lorsque $\overline{E}$ est hermitien. Ceci assure que les {\'e}nonc{\'e}s {\'e}tablis dans cet article sont des g{\'e}n{\'e}ralisations du cas \og classique\fg, hermitien sur $\spec\mathcal{O}_{k}$.\par Pour conclure, mentionnons que la th\'eorie des pentes sous sa forme originelle rev{\^e}t un aspect assez \'el\'ementaire \`a la fois en ce qui concerne les \'enonc\'es et les preuves. Afin de pr\'eserver son caract\`ere accessible au non-sp\'ecialiste, nous avons rappel\'e quelques rudiments de th\'eorie des ad\`eles et de g\'eom\'etrie de Minkowski. Nous ne supposons de la part du lecteur aucune connaissance particuli\`ere relative \`a la th{\'e}orie des pentes \og classique\fg. Les d{\'e}monstrations sont donn{\'e}es dans leur int{\'e}gralit{\'e}, fait susceptible d'entra{\^\i}ner {\c c}{\`a} et l\`a des r\'ep\'etitions.\\[0,1cm]\indent\textbf{Remerciements}. Je remercie Ga{\"e}l R\'emond de sa lecture attentive et critique d'une premi\`ere version de ce texte et des corrections qu'il m'a sugg\'er\'ees. Je remercie \'egalement le rapporteur de l'attention qu'il a pr{\^e}t{\'e}e \`a ce texte en me signalant de nombreuses erreurs et inexactitudes.

\tableofcontents

\section{Pr{\'e}liminaires}
\label{section:preliminaires}
\subsection{Ad\`eles sur un corps global}\label{paragrapheadelescorps}
Ce paragraphe pr{\'e}sente quelques propri{\'e}t{\'e}s bien connues des corps globaux et il fixe quelques notations, utilis{\'e}es dans la suite. Une {\'e}tude syst{\'e}matique des corps globaux et de leurs propri{\'e}t{\'e}s se trouve dans les ouvrages de C.~Chevalley~\cite{Chevalley} et de A.~Weil~\cite{WeilBasic}.\par Soit $k$ un corps global. Il y a deux cas de figure selon la caract{\'e}ristique de $k$.\\ \noindent \ding{192} \emph{La caract\'eristique de $k$ est nulle}.\\ Dans ce cas, le corps $k$ est un corps de nombres. On pose $k_{0}:=\mathbf{Q}$ et $D=[k:k_{0}]$ le degr{\'e} absolu de $k$. On d{\'e}signe par $\vert.\vert_{v}$ la valeur absolue sur le compl{\'e}t{\'e} $k_{v}$ (ou $\mathbf{C}_{v}$) de $k$ en la place $v$, normalis{\'e}e de la mani{\`e}re suivante.\begin{itemize}\item[1)] Si $v$ est archim{\'e}dienne, $\vert.\vert_{v}$ est la valeur absolue usuelle sur $\mathbf{R}$ ou $\mathbf{C}$.\item[2)] Si $v$ est ultram{\'e}trique, de caract{\'e}ristique r{\'e}siduelle $p_{v}$, on a $\vert p_{v}\vert_{v}=p_{v}^{-1}$.\end{itemize}\noindent\ding{193} \emph{La caract\'eristique de $k$ est $p>0$}.\\
Le corps $k$ est une extension finie de $k_{0}:=\mathbf{F}_{p}(T)$. On note $D:=[k:k_{0}]$. Une place $v$ de $k$ est n\'ecessairement ultram\'etrique et la place $v_{0}$ de $k_{0}$ correspondante est de deux sortes~: soit elle provient d'un polyn{\^o}me irr\'eductible $\pi$ de $\mathbf{F}_{p}[T]$, soit l'id\'eal premier associ\'e \`a $v_{0}$ est engendr\'e par $T^{-1}$. Dans ce second cas la place $v_{0}$ est dite \emph{infinie}, cette d\'esignation \'etant bien s{\^u}r fonction du choix de $T$. Sur le compl\'et\'e $(k_{0})_{v_{0}}$, on consid\`ere la valeur absolue $\vert.\vert_{v_{0}}$ normalis\'ee par $\vert\pi\vert_{v_{0}}=p^{-\deg\pi}$ (premier cas) ou $\vert T\vert_{v_{0}}=p$ (second cas). On prolonge alors cette valeur absolue \`a $k_{v}$ de telle mani\`ere \`a ce que $\vert x\vert_{v}=\vert x\vert_{v_{0}}$ pour $x\in (k_{0})_{v_{0}}$. Autrement dit, si $\mathrm{N}$ d\'esigne l'application norme de l'extension $k_{v}\mid (k_{0})_{v_{0}}$, on a $\vert x\vert_{v}=\vert\mathrm{N}(x)\vert_{v_{0}}^{1/n_{v}}$ o\`u $n_{v}$ est le degr\'e local $[k_{v}:(k_{0})_{v_{0}}]$.\par
Avec ces normalisations, si l'on pose $n_{v}:=1,2,[k_{v}:\mathbf{Q}_{p_{v}}],[k_{v}:(k_{0})_{v_{0}}]$ selon que $v$ est r{\'e}elle, complexe, ultram{\'e}trique ($\mathrm{carac}\,k=0$ ou $p$ respectivement), l'application $x\in k_{v}\mapsto\vert x\vert_{v}^{n_{v}}$ est le \emph{module de Haar normalis\'e en la place $v$} et la formule du produit s'{\'e}crit alors \begin{equation*}\forall\,x\in k\setminus\{0\},\quad\prod_{\text{$v$ place de $k$}}{\vert x\vert_{v}^{n_{v}}}=1\end{equation*}(dans ce produit tous les termes sauf un nombre fini valent $1$). Si $K$ est une extension finie de $k$, il n'y a qu'un nombre fini de places $w$ de $K$ au-dessus d'une place $v$ de $k$. On dispose d'un isomorphisme de $K$-alg{\`e}bres topologiques pour les compl{\'e}t{\'e}s \begin{equation}K\otimes_{k}k_{v}\simeq \prod_{w\mid v}{K_{w}},\end{equation}qui conduit en particulier {\`a} l'{\'e}galit{\'e} des degr{\'e}s $$[K:k]=\sum_{w\mid v}{[K_{w}:k_{v}]}$$(voir le chapitre~$4$ de~\cite{Chevalley}).
\par
Soit $k_{\mathbf{A}}$ l'anneau des ad{\`e}les de $k$. En tant que groupe localement compact, $(k_{\mathbf{A}},+)$ poss{\`e}de une mesure de Haar, unique {\`a} multiplication par un nombre r{\'e}el strictement positif pr{\`e}s. Le plongement diagonal $k\hookrightarrow k_{\mathbf{A}}$ conf{\`e}re {\`a} $k$ une structure de r{\'e}seau dans $k_{\mathbf{A}}$ et l'espace compact $k_{\mathbf{A}}/k$ a une mesure de Haar finie. Plus g{\'e}n{\'e}ralement il en est de m{\^e}me pour $E\hookrightarrow E\otimes k_{\mathbf{A}}=:E_{\mathbf{A}}$ o{\`u} $E$ est un $k$-espace vectoriel de dimension finie. La mesure dite de Tamagawa est celle pour laquelle la mesure du quotient $E_{\mathbf{A}}/E$ vaut $1$.\\ \indent $\bullet$ Si $k$ est un corps de nombres, soit $\mu_{v}$ la mesure de Haar d{\'e}finie sur le compl{\'e}t{\'e} $k_{v}$ de $k$ en une place $v$ de la mani{\`e}re suivante :\begin{itemize}
\item[a)] Si $v$ est r{\'e}elle, $\mu_{v}$ est la mesure de Lebesgue usuelle sur $\mathbf{R}$.\item[b)] Si $v$ est complexe, $\mu_{v}$ s'identifie au double $2\mathrm{d} x\mathrm{d} y$ de la mesure de Lebesgue sur $\mathbf{R}^{2}$.\item[c)] Si $v$ est ultram{\'e}trique, on pose $\mu_{v}(\mathcal{O}_{v})=1$ o{\`u} $\mathcal{O}_{v}$ est l'anneau de valuation de $k_{v}$.
\end{itemize}
La mesure produit $\mu=\prod{\mu_{v}}$ est une mesure de Haar sur $k_{\mathbf{A}}$ pour laquelle la mesure de l'espace quotient $k_{\mathbf{A}}/k$ (c.-{\`a}-d. la mesure d'un domaine fondamental de celui-ci) {\'e}gale $\vert D_{k}\vert^{1/2}$ o{\`u} $D_{k}$ est le discriminant absolu de $k$ (proposition~$7$ du chapitre~$5$ de~\cite{WeilBasic}).\\ \indent $\bullet$ Si $k$ est un corps de fonctions, on note $\mu$ la mesure de Haar sur $k_{\mathbf{A}}$ telle que $\mu(\prod_{v}{\mathcal{O}_{v}})=1$. On a alors $\mu(k_{\mathbf{A}}/k)=q^{g(k)-1}$ o{\`u} $g(k)\in\mathbf{N}$ est le genre de $k$ et $q$ d\'esigne le cardinal du plus grand corps fini inclus dans $k$ (corollaire~$1$ du chapitre~$6$, \textit{ibid.}).
\par Plus g{\'e}n{\'e}ralement, si $E$ est un $k$-espace vectoriel de dimension finie, le choix d'une $k$-base de $E$ fournit un isomorphisme $E\otimes k_{\mathbf{A}}\simeq k_{\mathbf{A}}^{\dim E}$ et une mesure de Haar $\mu_{E_{\mathbf{A}}}$ sur $E_{\mathbf{A}}$. Cette mesure ne d{\'e}pend pas du choix de la base de $E$.\par
 Avec ces normalisations, pour toute place $v$ de $k$ et tout nombre r{\'e}el $r\in\vert k_{v}\vert_{v}$, on a \begin{equation}\label{equation:rayonboule}\mu_{v}\left(\left\{x\in k_{v}\,;\ \vert x\vert_{v}\le r\right\}\right)=r^{n_{v}}\mu_{v}\left(\left\{x\in k_{v}\,;\ \vert x\vert_{v}\le 1\right\}\right)\ .\end{equation}Cette \'egalit\'e reste valide si l'on remplace $\mu_{v}$ par une mesure de Haar quelconque sur $(k_{v},+)$. Enfin, il est commode de d{\'e}finir pour un ad{\`e}le $a=(a_{v})_{v}\in k_{\mathbf{A}}$ la \emph{valeur absolue ad\'elique} de $a$ comme le nombre r\'eel $$\vert a\vert_{\mathbf{A}}:=\prod_{\text{$v$ place de $k$}}{\vert a_{v}\vert_{v}^{n_{v}}}\ .$$Si $k$ est un corps de nombres, l'image par $\vert\cdotp\vert_{\mathbf{A}}$ des id{\`e}les $k_{\mathbf{A}}^{\times}$ de $k_{\mathbf{A}}$ est $\mathbf{R}^{*}_{+}$. Si $k$ est un corps de fonctions, l'image $\vert k_{\mathbf{A}}^{\times}\vert_{\mathbf{A}}$ est $\{q^{n}\,;\ n\in\mathbf{Z}\}$ (voir \cite{WeilBasic}, chap.~$7$, \S~$5$, corollaire~$6$).
\subsection{G{\'e}om{\'e}trie de Minkowski}\label{paragraphedeuxdeux} Pour comprendre les propri\'et\'es d'une norme quelconque sur $\mathbf{R}^{n}$, une possibilit\'e est d'effectuer une comparaison avec la norme de $\ell_{n}^{p}$, $p\in[1,+\infty]$, et, si possible, avec la norme euclidienne de $\ell_{n}^{2}$. La \emph{g{\'e}om{\'e}trie de Minkowski} est l'{\'e}tude des $\mathbf{R}$-espaces vectoriels norm{\'e}s $(E,\Vert.\Vert)$ de dimension finie. L'un des axes de cette {\'e}tude consiste pr{\'e}cis{\'e}ment {\`a} s'int{\'e}resser {\`a} la structure euclidienne qui se rapproche le plus de la structure d'espace vectoriel norm{\'e} de $E$, en un sens que nous pr{\'e}ciserons un peu plus loin. Dans le cas d'un corps de nombres, cette approche s'av{\'e}rera {\^e}tre la clef qui permet d'aborder la th{\'e}orie des fibr{\'e}s vectoriels ad{\'e}liques.\par Pour {\'e}crire cette synth{\`e}se, nous avons consult{\'e} les ouvrages de G.~Pisier~\cite{Pisier1989}, R.~Ryan~\cite{Ryan} et A.~Thompson~\cite{Thompson} ainsi que les articles~\cite{BourgainMilman,Rogalski1981I,SaintRaymond1981}.
\par Dans tout ce paragraphe, $E$ est un $\mathbf{R}$-espace vectoriel de dimension $n\ge 1$, muni d'une mesure de Haar $\vol$. L'espace vectoriel dual $E^{\mathsf{v}}=\Hom_{\mathbf{R}}(E,\mathbf{R})$ poss{\`e}de alors une mesure de Haar particuli{\`e}re $\vol^{*}$, associ{\'e}e au choix de $\vol$, caract{\'e}ris{\'e}e de la mani{\`e}re suivante : soit $(e_{1},\ldots,e_{n})$ une base de $E$ et le parall{\'e}lotope $$\mathscr{P}:=\left\{x_{1}e_{1}+\cdots+x_{n}e_{n}\,;\ 0\le x_{i}\le 1\right\}\cdotp$$De m{\^e}me, on dispose de la base duale $(e_{1}^{*},\ldots,e_{n}^{*})$ et du parall{\'e}lotope dual $\mathscr{P}^{*}$ associ{\'e}. Alors $\vol^{*}$ est l'unique mesure de Haar telle que \begin{equation*}\vol(\mathscr{P})\vol^{*}(\mathscr{P}^{*})=1\ .\end{equation*}Soit $C\subseteq E$ une partie convexe, d'int{\'e}rieur non vide, compact et sym{\'e}trique par rapport {\`a} l'origine. Un tel ensemble est appel{\'e} \emph{corps convexe (sym{\'e}trique\footnote{Tous les corps convexes consid{\'e}r{\'e}s dans ce texte sont sym{\'e}triques et nous omettrons de le pr{\'e}ciser {\`a} chaque fois.})}. La fonction \emph{jauge} $$\forall x\in E,\quad j(x):=\inf{\left\{\lambda>0\,;\ \frac{x}{\lambda}\in C\right\}}$$fait le lien entre une norme sur $E$ et le corps convexe qui est la boule unit{\'e} pour cette norme. Autrement dit, le couple $(E,C)$ est la donn{\'e}e d'une structure norm{\'e}e $(E,j)$ sur l'espace vectoriel $E$. Dans la suite, nous noterons aussi $\Vert.\Vert$ la norme $j$ sur $E$.
\begin{defi}Le \emph{polaire} de $C$, not{\'e} $C^{\circ}$, est l'ensemble $C^{\circ}:=\left\{\varphi\in E^{\mathsf{v}}\,;\ \vert\varphi(C)\vert\subseteq [0,1]\right\}$.
\end{defi}\noindent On v{\'e}rifie que $C^{\circ}$ est la boule unit{\'e} ferm{\'e}e de l'espace dual $(E^{\mathsf{v}},j^{\mathsf{v}})$. Lorsque $p\in [1,+\infty]$, on note $b_{n}^{p}$ (ou $b_{n,\mathbf{R}}^{p}$) la boule unit{\'e} ferm{\'e}e de l'espace de Banach $\ell_{n}^{p}:=(\mathbf{R}^{n},\vert.\vert_{p})$, o\`u \begin{equation}\label{defidenormep}\vert(x_{1},\ldots,x_{n})\vert_{p}=\begin{cases}\left(\vert x_{1}\vert^{p}+\cdots+\vert x_{n}\vert^{p}\right)^{1/p} & \text{si $p\in [1,+\infty[$},\\ \max{\{\vert x_{1}\vert,\ldots,\vert x_{n}\vert\}} & \text{si $p=+\infty$}.\end{cases}\end{equation}Si $p\in]0,1[$ et $n\ge 2$, cette application $\vert.\vert_{p}$ n'est plus une norme et $b^{p}_{n}$ n'est plus convexe, mais cela reste un ensemble mesurable au sens de Lebesgue. Ainsi, pour tout $p>0$ et pour la mesure de Lebesgue $\vol_{n}$ sur $\mathbf{R}^{n}$, on a \begin{equation}\label{formuledevolumep}\vol_{n}(b_{n}^{p})=\frac{\left(2\Gamma\left(1+\frac{1}{p}\right)\right)^{n}}{\Gamma\left(1+\frac{n}{p}\right)}\end{equation}o{\`u} $\Gamma(x):=\int_{0}^{+\infty}{t^{x-1}e^{-t}\,\mathrm{d}t}$ (voir formule~\eqref{formuleintegrale} un peu plus loin). De m{\^e}me, $b_{n,\mathbf{C}}^{p}$ d{\'e}signe la boule unit{\'e} ferm{\'e}e pour la norme $\vert.\vert_{p}$ sur $\mathbf{C}^{n}$. En identifiant $\mathbf{C}^{n}$ {\`a} $\mathbf{R}^{2n}$, la mesure de Lebesgue de $b_{n,\mathbf{C}}^{p}$ vaut \begin{equation}\label{formulevolumecomplexep}\vol_{2n}\left(b_{n,\mathbf{C}}^{p}\right)=\left(\frac{\pi}{2}\right)^{n}\vol_{n}\left(b_{n,\mathbf{R}}^{p/2}\right)\ .\end{equation}
\par {\`A} un corps convexe $C$, l'on peut associer l'invariant suivant :\begin{defi}Le \emph{produit de Mahler} du corps convexe $C$, not{\'e} $P(C)$, est le produit $\vol(C)\vol^{*}(C^{\circ})$.\end{defi}Ce nombre r{\'e}el ne d{\'e}pend pas du choix de la mesure de Haar $\vol$ sur $E$. De plus, il est invariant par isomorphisme : si $u\in\GL(E)$ alors $P(u(C))=P(C)$. C'est pourquoi l'on peut omettre la r{\'e}f{\'e}rence {\`a} $E$ dans la notation du produit de Mahler. 
\begin{theo}\label{theoreme23bourgain} Pour tout corps convexe $C$, on a \begin{equation*}P(C)\le P(b_{n}^{2})\qquad\text{(Blaschke-Santal{\'o})}\end{equation*}et l'existence d'une constante absolue $c\in\,]0,+\infty[$ telle que \begin{equation*}P(C)\ge e^{-cn}P(b_{n}^{2})\qquad\text{(Bourgain-Milman)}.\end{equation*}\end{theo}La premi{\`e}re in{\'e}galit{\'e} a {\'e}t{\'e} {\'e}tablie dans~\cite{Blaschke,Santalo}, articles auxquels on peut adjoindre le texte de J.~Saint-Raymond~\cite{SaintRaymond1981} qui comporte une preuve \og {\'e}l{\'e}mentaire\fg\ de ce r{\'e}sultat et qui d{\'e}montre qu'il y a {\'e}galit{\'e} seulement si, dans une certaine base de $E$, le convexe $C$ est la boule unit{\'e} euclidienne usuelle (on dit alors que $C$ est un \emph{ellipso{\"\i}de}). La seconde in{\'e}galit{\'e}, d{\'e}montr{\'e}e dans~\cite{BourgainMilman}, est plus difficile {\`a} obtenir. Une valeur explicite pour la constante $c$ n'est pas connue {\`a} l'heure actuelle. Aussi peut-il {\^e}tre utile de mentionner une minoration plus faible due {\`a} K.~Mahler~\cite{Mahler1939} : $$P(C)\ge\frac{4^{n}}{(n!)^{2}},$$minoration qui entra{\^\i}ne\footnote{\label{Stirling}Le calcul qui conduit {\`a} cette minoration est bas{\'e} sur la formule donn{\'e}e auparavant pour le volume de $b_{n}^{2}$ ainsi que sur l'existence d'une fonction $\eta$, d{\'e}croissante, positive et nulle \`a l'infini, telle que $$\forall x>0,\quad\Gamma(1+x)=\sqrt{2\pi x}\left(\frac{x}{e}\right)^{x}e^{\eta(x)}$$(voir~\cite{Remmert}, chap.~2, \S~4).} $$P(C)\ge e^{-n\log(en)}P(b_{n}^{2})\ .$$Une preuve de la conjecture de Mahler $P(C)\ge 4^{n}/n!$ permettrait d'obtenir une constante $c$ explicite dans l'in\'egalit\'e de Bourgain-Milman.
\subsubsection*{Ellipso{\"\i}des de John et L{\"o}wner} Un \emph{ellipso{\"\i}de} de $E$ est un ensemble de la forme $D=\{x\in E\,;\ q(x)\le 1\}$ o{\`u} $q:E\to\mathbf{R}$ est une forme quadratique d{\'e}finie positive. Si l'on fixe une base de $E$ qui permet d'identifier $E$ {\`a} $\mathbf{R}^{n}$, un ellipso{\"\i}de est l'image de la boule unit{\'e} euclidienne $b_{n}^{2}$ par un isomorphisme de $E$. En particulier, le produit de Mahler de $D$ est celui de $b_{n}^{2}$.
\begin{defithm}\label{defiellipsoidejohn}
{\'E}tant donn{\'e} un corps convexe (sym{\'e}trique) $C$, il existe un unique ellipso{\"\i}de $J(C)$, appel{\'e} \emph{ellipso{\"\i}de de John}, inclus dans $C$ et de volume maximal. De m{\^e}me, il existe un unique ellipso{\"\i}de $L(C)$, dit \emph{ellipso{\"\i}de de L{\"o}wner}, contenant $C$ et de volume minimal. 
\end{defithm}La seconde assertion d{\'e}coule de la premi{\`e}re par dualit{\'e} : $J(C)^{\circ}=L(C^{\circ})$, car $P(J(C))$ est constant, {\'e}gal {\`a} $P(b_{n}^{2})$. L'unicit{\'e} est le point difficile et remarquable de cet {\'e}nonc{\'e}. Elle entra{\^\i}ne les inclusions $C\subseteq\sqrt{n}J(C)$ et $L(C)\subseteq\sqrt{n}C$ (ce n'est pas imm{\'e}diat, voir~\cite{Thompson}, p.~$84$). En notant $\vert.\vert_{J(C)}$ (\emph{resp}. $\vert.\vert_{L(C)}$)\label{defidesmetriquesLowner} la norme euclidienne sur $E$ associ{\'e}e {\`a} $J(C)$ (\emph{resp}. $L(C)$), cela se traduit par les in{\'e}galit{\'e}s\begin{equation*}\forall x\in E,\quad \frac{1}{\sqrt{n}}\vert x\vert_{J(C)}\le j(x)\le\vert x\vert_{J(C)}\quad\text{et}\quad \vert x\vert_{L(C)}\le j(x)\le\sqrt{n}\vert x\vert_{L(C)}\ .\end{equation*}On dispose ainsi de deux structures euclidiennes qui encadrent la norme donn{\'e}e sur $E$.

\subsubsection*{Quotient volumique}
Ce paragraphe emprunte beaucoup au texte~\cite{Rogalski1981I} de M.~Rogalski. On note $\mathbf{B}(E,\Vert.\Vert)$ la boule unit\'e ferm\'ee de l'espace vectoriel norm\'e $(E,\Vert.\Vert)$.
\begin{defi}\label{defiquotientvolumique}
Le \emph{quotient volumique}\footnote{\emph{Volume ratio} en anglais.} de $(E,\Vert.\Vert)$, not{\'e} $\vr(E)$, est le nombre r{\'e}el $\ge 1$ d\'efini par\begin{equation*}\vr(E):=\inf{\left\{\left(\frac{\vol(\mathbf{B}(E,\Vert.\Vert))}{\vol(D)}\right)^{1/n}\,;\ D\ \text{ellipso{\"\i}de}\ \subseteq\mathbf{B}(E,\Vert.\Vert)\right\}}\ \cdotp\end{equation*}\end{defi}La d{\'e}finition m{\^e}me de l'ellipso{\"\i}de John entra{\^\i}ne $$\vr(E)=\left(\frac{\vol(\mathbf{B}(E,\Vert.\Vert))}{\vol(J(\mathbf{B}(E,\Vert.\Vert)))}\right)^{1/n},$$et l'inclusion $C\subseteq\sqrt{n}J(C)$ donne alors $\vr(E)\le\sqrt{n}$. Dans cette in{\'e}galit{\'e}, la fonction $n\mapsto\sqrt{n}$ ne peut pas {\^e}tre remplac{\'e}e par une fonction $f$ telle que $f(n)/\sqrt{n}\underset{n\to+\infty}{\longrightarrow}0$. En revanche, on peut montrer que $\vr(E)\le\delta\sqrt{n}$ o\`u $\delta\in]0,1[$ est explicite (par exemple $\delta=0,95$). M.~Rogalski a calcul{\'e} le quotient volumique de $\ell_{n}^{p}$ :\begin{equation*}\vr(\ell_{n}^{p})=\Phi_{p}(n)n^{\max{\left\{0,\left(\frac{1}{2}-\frac{1}{p}\right)\right\}}}\sqrt{\frac{2}{\pi}}\Gamma\left(1+\frac{1}{p}\right)e^{\frac{1}{p}-\frac{1}{2}}p^{\frac{1}{p}}\end{equation*}o\`u $\Phi_{p}(n)\in[1/3,2]$ et $\Phi_{p}(n)\underset{n\to+\infty}{\longrightarrow}1$.
\par Nous aurons {\'e}galement besoin de la variante avec l'ellipso{\"\i}de de L{\"o}wner.
\begin{defi}\label{quotientvolumiquebis}On note $\widetilde{\vr}(E)$ le nombre r{\'e}el $\ge 1$ d{\'e}fini par \begin{equation*}\widetilde{\vr}(E):=\sup{\left\{\left(\frac{\vol(D)}{\vol(\mathbf{B}(E,\Vert.\Vert))}\right)^{1/n}\,;\ \mathbf{B}(E,\Vert.\Vert)\subseteq D\ \text{ellipso{\"\i}de}\right\}}\ \cdotp\end{equation*}
\end{defi}On a donc $$\widetilde{\vr}(E)=\left(\frac{\vol(L(\mathbf{B}(E,\Vert.\Vert)))}{\vol(\mathbf{B}(E,\Vert.\Vert))}\right)^{1/n}\quad\text{et}\quad \widetilde{\vr}(E)\le\sqrt{n}\ .$$
\subsubsection*{Distance de Banach-Mazur}
\begin{defi}\label{distancebanachmazur}
\'Etant donn{\'e} deux corps convexes $C_{1}$ et $C_{2}$ de $E$, la distance $\mathrm{d}(C_{1},C_{2})$ entre ces deux ensembles est \begin{equation*}\mathrm{d}(C_{1},C_{2}):=\inf{\left\{ab\,;\quad a>0,\ b>0,\quad C_{1}\subseteq aC_{2}\quad\text{et}\quad C_{2}\subseteq bC_{1}\right\}}\ \cdotp\end{equation*}La \emph{distance} --- dite \emph{de Banach-Mazur} --- entre les espaces de Banach $E_{1}=(E,C_{1})$ et $E_{2}=(E,C_{2})$ est \begin{equation*}\mathrm{d}(E_{1},E_{2}):=\inf{\left\{\mathrm{d}(C_{1},u(C_{2}))\,;\quad u\in\GL(E)\right\}}\ \cdotp\end{equation*}Plus g{\'e}n{\'e}ralement, la distance de Banach-Mazur entre deux espaces de Banach $E$ et $F$ de m{\^e}me dimension (finie) est $\mathrm{d}(E,\varphi(F))$ o\`u $\varphi:F\to E$ est un isomorphisme quelconque entre $F$ et $E$.\end{defi}
Cette derni{\`e}re quantit{\'e} ne d{\'e}pend pas du choix de $\varphi$. La terminologie \og distance\fg\ se justifie par l'in{\'e}galit{\'e} $\mathrm{d}(E_{1},E_{3})\le\mathrm{d}(E_{1},E_{2})\mathrm{d}(E_{2},E_{3})$ o\`u $E_{1},E_{2},E_{3}$ sont des espaces vectoriels norm{\'e}s de m{\^e}me dimension. La distance qui nous int{\'e}ressera le plus dans la suite est celle entre $(E,\Vert.\Vert)$ et $\ell^{2}_{n}$. Par d{\'e}finition m{\^e}me de cette distance, pour tout $\varepsilon>0$, il existe une norme \emph{euclidienne} $\vert.\vert_{\varepsilon}$ sur $E$ telle que, pour tout $x\in E$, on ait \begin{equation}\label{encadrementhermitien}\forall x\in E,\quad\vert x\vert_{\varepsilon}\le\Vert x\Vert\le\mathrm{d}(E,\ell^{2}_{n})(1+\varepsilon)\vert x\vert_{\varepsilon}\ .\end{equation}On v{\'e}rifie alors l'encadrement \begin{equation}\label{estimationslocalesreelles}1\le\vr(E)\vr(E^{\mathsf{v}})\le\mathrm{d}(E,\ell^{2}_{n})\le\sqrt{n}\ ,\end{equation}sans qu'il y ait {\'e}galit{\'e} en g{\'e}n{\'e}ral. En r{\'e}alit{\'e} l'on peut m{\^e}me exhiber une suite d'espaces norm{\'e}s $E_{n}$ de dimension $n$ telle que \begin{equation*}\mathrm{d}(E_{n},\ell^{2}_{n})/\vr(E_{n})\underset{n\to+\infty}{\longrightarrow}+\infty\ .\end{equation*}On a aussi $1\le\widetilde{\vr}(E)\widetilde{\vr}(E^{\mathsf{v}})\le\mathrm{d}(E,\ell^{2}_{n})$.  
\subsubsection*{Somme directe}\label{paragraphedeuxsommedirecte}
Soit $E,F$ deux espaces vectoriels norm{\'e}s de dimensions respectives $n$ et $m$. Soit $\varsigma$ une norme sym{\'e}trique sur $\mathbf{R}^{2}$, invariante par changement de signes sur les coordonn{\'e}es, telle que $\varsigma(1,0)=\varsigma(0,1)=1$. On munit la somme directe $E\oplus F$ de la norme $\Vert(x,y)\Vert_{E\oplus_{\varsigma} F}:=\varsigma(\Vert x\Vert_{E},\Vert y\Vert_{F})$. L'espace vectoriel norm{\'e} obtenu sera not{\'e} $E\oplus_{\varsigma}F$. On v{\'e}rifie \begin{equation}\label{sommedirecteplus210}\forall x\in E,\ \forall y\in F,\quad\max{\{\Vert x\Vert_{E},\Vert y\Vert_{F}\}}\le\Vert(x,y)\Vert_{E\oplus_{\varsigma}F}\le\Vert x\Vert_{E}+\Vert y\Vert_{F}\ .\end{equation}
\begin{prop}\label{proposition29sommedirecte}
Avec les donn{\'e}es ci-dessus, soit $\vol_{E}$ (\emph{resp}. $\vol_{F}$) une mesure de Haar sur $E$ (\emph{resp}. $F$). On dispose de la mesure produit $\vol_{E\oplus F}:=\vol_{E}\otimes\vol_{F}$ sur $E\oplus F\simeq E\times F$. On a alors
\begin{equation*}\binom{n+m}{n}^{-1}\le\frac{\vol_{E\oplus F}(\mathbf{B}(E\oplus F,\Vert.\Vert_{E\oplus_{\varsigma}F}))}{\vol_{E}(\mathbf{B}(E,\Vert.\Vert_{E}))\vol_{F}(\mathbf{B}(F,\Vert.\Vert_{F}))}\le 1\ .\end{equation*}
\end{prop}
\begin{proof}
L'in{\'e}galit{\'e} de droite est une simple cons{\'e}quence de l'inclusion $$\mathbf{B}(E\oplus F,\Vert.\Vert_{E\oplus_{\varsigma}F})\subseteq\mathbf{B}(E,\Vert.\Vert_{E})\times\mathbf{B}(F,\Vert.\Vert_{F}),$$qui r{\'e}sulte de la majoration $\max\{\Vert x\Vert_{E},\Vert y\Vert_{F}\}\le\Vert(x,y)\Vert_{E\oplus_{\varsigma}F}$. L'autre in{\'e}galit{\'e} repose sur la formule int{\'e}grale \begin{equation}\label{formuleintegrale}\forall p>0,\quad\vol_{E}(\mathbf{B}(E,\Vert.\Vert_{E}))\times\Gamma\left(1+\frac{n}{p}\right)=\int_{E}{e^{-\Vert x\Vert_{E}^{p}}\,\mathrm{d}(\vol_{E})(x)},\end{equation}qui s'obtient gr{\^a}ce au th{\'e}or{\`e}me de Fubini en int{\'e}grant $(e^{-t}\mathrm{d}t)\otimes\mathrm{d}(\vol_{E})$ sur l'ensemble $\{(t,x)\,;\ t\ge\Vert x\Vert_{E}^{p}\}\subseteq\mathbf{R}\times E$. On applique cette formule {\`a} $E\oplus_{\varsigma}F$ et $p=1$. De la majoration $\Vert(x,y)\Vert_{E\oplus_{\varsigma}F}\le\Vert x\Vert_{E}+\Vert y\Vert_{F}$ d{\'e}coule l'in{\'e}galit{\'e} $$\vol_{E\oplus F}(\mathbf{B}(E\oplus F,\Vert.\Vert_{E\oplus_{\varsigma}F}))(n+m)!\ge\int_{E\oplus F}{e^{-\Vert x\Vert_{E}-\Vert y\Vert_{F}}\mathrm{d}(\vol_{E\oplus F})(x,y)},$$et la derni{\`e}re int{\'e}grale vaut exactement $(\vol_{E}(\mathbf{B}(E,\Vert.\Vert_{E}))n!)(\vol_{F}(\mathbf{B}(F,\Vert.\Vert_{F}))m!)$. Ceci conclut la d{\'e}monstration.
\end{proof}
\begin{rema}
Si $\varsigma(\alpha,\beta)=(\vert\alpha\vert^{p}+\vert\beta\vert^{p})^{1/p}$, $p\ge 1$, la formule int{\'e}grale~\eqref{formuleintegrale} fournit l'{\'e}galit{\'e} \begin{equation}\label{sommedirectemoins211}
\frac{\vol_{E\oplus F}(\mathbf{B}(E\oplus F,\Vert.\Vert_{E\oplus_{\varsigma}F}))}{\vol_{E}(\mathbf{B}(E,\Vert.\Vert_{E}))\vol_{F}(\mathbf{B}(F,\Vert.\Vert_{F}))}=\frac{\Gamma\left(1+\frac{n}{p}\right)\Gamma\left(1+\frac{m}{p}\right)}{\Gamma\left(1+\frac{n+m}{p}\right)}\ \cdotp\end{equation}Autrement dit, le quotient des volumes ne d{\'e}pend que des dimensions de $E$ et de $F$.
\end{rema}

\subsubsection*{Normes tensorielles}\label{paragraphenormestensorielles}
Sur le produit tensoriel de deux espaces vectoriels norm\'es de dimension finie coexistent de nombreuse normes, que l'on peut obtenir de mani\`ere \og naturelle\fg\ \`a partir des normes sur les espaces de d\'epart. Les conditions minimales que nous exigerons pour une telle norme sur le produit tensoriel sont les suivantes.

\begin{defi}\label{defidenormestensoriellesreelles} Soit $\ell$ un entier $\ge 1$. Une \emph{norme tensorielle\footnote{\emph{(Finitely generated) uniform crossnorm}, selon la terminologie de Schatten (voir~\cite{Ryan}, chap.~$6$).} d'ordre $\ell$} est la donn\'ee pour tous espaces de Banach $(E_{i},\Vert.\Vert_{E_{i}})_{1\le i\le\ell}$ (sur $\mathbf{R}$ ou $\mathbf{C}$) de dimension finie d'une norme $\alpha(\cdot;E_{1},\ldots,E_{\ell})$ sur le produit tensoriel $E_{1}\otimes\cdots\otimes E_{\ell}$ telle que~:
\begin{enumerate}\item[(i)] pour tout $i\in\{1,\ldots,\ell\}$, pour tout $e_{i}\in E_{i}$, on a \begin{equation*}\alpha(e_{1}\otimes\cdots\otimes e_{\ell};E_{1},\ldots,E_{\ell})\le\prod_{i=1}^{\ell}{\Vert e_{i}\Vert_{E_{i}}},\end{equation*}
\item[(ii)] la famille $\alpha=(\alpha(\cdot;E_{1},\ldots,E_{\ell}))_{(E_{1},\ldots, E_{\ell})}$ doit v\'erifier~: Pour tout $i\in\{1,\ldots,\ell\}$, pour tous espaces norm{\'e}s $E_{i}$ et $F_{i}$ de dimension finie et toutes applications lin\'eaires $u_{i}:E_{i}\to F_{i}$, la norme d'op\'erateur de $u_{1}\otimes\cdots\otimes u_{\ell}:E_{1}\otimes\cdots\otimes E_{\ell}\to F_{1}\otimes\cdots\otimes F_{\ell}$ est plus petite que le produit des normes d'op{\'e}rateur $\prod_{i=1}^{\ell}{\Vert u_{i}\Vert}$.
\end{enumerate}Une \emph{norme tensorielle hermitienne} est une norme tensorielle qui restreinte aux espaces hermitiens est la norme hermitienne usuelle\footnote{Autrement dit, si $E$ et $F$ sont hermitiens, munis de bases orthonorm\'ees $(e_{1},\ldots,e_{n})$ et $(f_{1},\ldots,f_{m})$ respectivement, et si $x=\sum_{i,j}{x_{i,j}e_{i}\otimes f_{j}}\in E\otimes F$ alors $\alpha(x;E,F)^{2}=\sum_{i,j}{\vert x_{i,j}\vert^{2}}$.} sur le produit tensoriel. Nous noterons $E_{1}\otimes_{\alpha}\cdots\otimes_{\alpha}E_{\ell}$ ou $\otimes_{\alpha,i=1}^{\ell}E_{i}$ l'espace vectoriel $E_{1}\otimes\cdots\otimes E_{\ell}=:\otimes_{i=1}^{\ell}E_{i}$ muni de la norme $\alpha(\cdot;E_{1},\ldots,E_{\ell})$.
\end{defi}
On montre ais\'ement que si les conditions (i) et (ii) ci-dessus sont r\'ealis\'ees alors il y a \'egalit\'e~: $\alpha(e_{1}\otimes\cdots\otimes e_{\ell})=\prod_{i=1}^{\ell}{\Vert e_{i}\Vert_{E_{i}}}$ et $\Vert u_{1}\otimes\cdots\otimes u_{\ell}\Vert=\prod_{i=1}^{\ell}{\Vert u_{i}\Vert}$.
\begin{exems}\label{normesChevetSaphard}\emph{Normes de Chevet-Saphard (d'ordre $2$)}. Soit $(E,\Vert.\Vert_{E})$ et $(F,\Vert.\Vert_{F})$ deux espaces vectoriels norm\'es de dimension finie. Soit $\ell\in\mathbf{N}\setminus\{0\}$. Soit $p\in[1,+\infty]$ et $p'$ son conjugu\'e, c.-{\`a}-d. $1/p+1/p'=1$. Si $(e_{1},\ldots,e_{\ell})\in E^{\ell}$, on note \begin{equation*}\Vert(e_{1},\ldots,e_{\ell})\Vert_{p}:=\begin{cases}\left(\sum_{i=1}^{\ell}{\Vert e_{i}\Vert_{E}^{p}}\right)^{1/p} & \text{si $1\le p<\infty$},\\ \max_{1\le i\le \ell}{\{\Vert e_{i}\Vert_{E}\}} & \text{si $p=+\infty$}.\end{cases}\end{equation*}De m{\^e}me, si $(f_{1},\ldots,f_{\ell})\in F^{\ell}$, on pose \begin{equation*}\Vert(f_{1},\ldots,f_{\ell})\Vert_{p}^{w}:=\sup{\left\{\left\Vert\sum_{i=1}^{\ell}{\lambda_{i}f_{i}}\right\Vert_{F};\ (\lambda_{1},\ldots,\lambda_{\ell})\in b_{\ell}^{p}\right\}}\end{equation*}(le $w$ est l'initiale de \emph{weak}). Alors la \emph{norme \`a gauche de Chevet-Saphard}, not\'ee $g_{p}$ ou $g_{p}(\cdot;E,F)$, est la norme sur $E\otimes F$ d\'efinie par~: pour tout $x\in E\otimes F$, \begin{equation*}g_{p}(x;E,F):=\inf{\left\{\Vert(e_{1},\ldots,e_{\ell})\Vert_{p}\Vert(f_{1},\ldots,f_{\ell})\Vert_{p'}^{w};\ x=\sum_{i=1}^{\ell}{e_{i}\otimes f_{i}}\right\}}\ \cdot\end{equation*}Cela d\'efinit une norme tensorielle d'ordre $2$ (voir~\cite{Ryan}) et, pour $p=2$, la norme $g_{2}$ est une norme tensorielle hermitienne au sens ci-dessus. 
\end{exems}
Cet exemple et une simple r\'ecurrence permettent de construire des normes tensorielles (\'eventuellement hermitiennes) d'ordre quelconque.\par L'observation suivante, bien que tr\`es simple, sera d'un usage constant dans la suite (le corps de base est $\mathbf{R}$ ou $\mathbf{C}$).
\begin{prop}\label{propositioncompnormes}Soit $\ell\in\mathbf{N}\setminus\{0\}$ et $\alpha$ une norme tensorielle d'ordre $\ell$ . Pour tout $i\in\{1,\ldots,\ell\}$, soit $E_{i}'\subseteq E_{i}$ des espaces vectoriels de dimension finie. Soit $\vert.\vert_{E_{i}'}$ et $\vert.\vert_{E_{i}}$ deux normes sur $E_{i}'$ et $E_{i}$ respectivement. Supposons que, pour tout $e_{i}\in E_{i}'$, on a $\vert e_{i}\vert_{E_{i}}\le\vert e_{i}\vert_{E_{i}'}$. Alors, pour tout $e\in\otimes_{i=1}^{\ell}E_{i}'$, on a \begin{equation*}\alpha(e;(E_{1},\vert.\vert_{E_{1}}),\ldots,(E_{\ell},\vert.\vert_{E_{\ell}}))\le\alpha(e;(E_{1}',\vert.\vert_{E_{1}'}),\ldots,(E_{\ell}',\vert.\vert_{E_{\ell}'}))\ .\end{equation*}
\end{prop}
\begin{proof}
Il s'agit de la condition (iii) caract\'erisant la norme tensorielle $\alpha$, appliqu\'ee aux morphismes d'inclusion $x\mapsto x$ de $(E_{i}',\vert.\vert_{E_{i}'})$ dans $(E_{i},\vert.\vert_{E_{i}})$ pour tout $i\in\{1,\ldots,\ell\}$.
\end{proof}
\begin{remas}
\begin{enumerate}\item[a)] Le lecteur peu accoutum\'e prendra garde aux nombreuses chausse-trapes des normes tensorielles. Pour $\alpha$ d'ordre $2$, la loi interne $(E,F)\mapsto E\otimes_{\alpha}F$ sur l'ensemble des espaces norm\'es de dimension finie n'est en g\'en\'eral ni associative ni commutative. De plus, si $F$ est un sous-espace vectoriel de $E$, la norme $\alpha(\cdot;F,F)$ sur $F\otimes F$ n'est en g\'en\'eral pas la restriction de la norme $\alpha(\cdot;E,E)$ \`a $F\otimes F$ (c.-{\`a}-d. $F\otimes_{\alpha}F$ n'est pas un sous-Banach de $E\otimes_{\alpha}E$). Il en est de m{\^e}me pour un quotient $E\twoheadrightarrow G$ ou pour la dualit\'e~: la norme sur $(E^{\mathsf{v}}\otimes_{\alpha}F^{\mathsf{v}})^{\mathsf{v}}$ n'est en g\'en\'eral pas \'egale \`a $\alpha(\cdot;E,F)$. 
\item[b)] Dans la suite, nous demanderons souvent aux normes tensorielles d'{\^e}tre \emph{hermitiennes}. Au del{\`a} de la norme $g_{2}$ de l'exemple~\ref{normesChevetSaphard}, il existe en r\'ealit\'e une infinit\'e de telles normes tensorielles, non \'equivalentes entre elles, comme l'affirme un r\'esultat de Schatten \& Puhl (voir~\cite{DefantFloret}, p.~357).   
\end{enumerate}
\end{remas}
\subsubsection*{Cas d'un espace vectoriel complexe}Dans ce paragraphe, $E$ est un $\mathbf{C}$-espace vectoriel de dimension $n\ge 1$. Soit $C\subseteq E$ un sous-ensemble convexe et compact, d'int{\'e}rieur non vide. Afin que la jauge d{\'e}finisse une norme sur $E$ (en particulier pour que la relation $j(e^{i\theta}x)=j(x)$ soit v{\'e}rifi{\'e}e pour tous $\theta\in\mathbf{R}$ et $x\in E$), on suppose \begin{equation}\label{hypothesejaugecomplexe}\forall\theta\in\mathbf{R},\quad e^{i\theta}\cdotp C=C\ .\end{equation}Notons $E_{\mathbf{R}}$ le $\mathbf{R}$-espace vectoriel sous-jacent {\`a} $E$ et $C_{\mathbf{R}}$ le sous-ensemble de $E_{\mathbf{R}}$ induit par $C$. Aux objets $E_{\mathbf{R}}$ et $C_{\mathbf{R}}$ l'on peut appliquer les r{\'e}sultats pr{\'e}c{\'e}dents et, en particulier, l'on dispose des ellipso{\"\i}des de John et L{\"o}wner associ{\'e}s {\`a} $C_{\mathbf{R}}$. Les normes euclidiennes sur $E_{\mathbf{R}}$ donn{\'e}es par ces ellipso{\"\i}des d{\'e}finissent des \emph{normes hermitiennes} sur $E$. En effet l'hypoth{\`e}se~\eqref{hypothesejaugecomplexe}, l'unicit{\'e} de ces ellipso{\"\i}des et la conservation des volumes par les applications \begin{equation*}x=x_{1}+ix_{2}\in E_{\mathbf{R}}\mapsto e^{i\theta}\cdot x=x_{1}\cos\theta-x_{2}\sin\theta+i(x_{1}\sin\theta+x_{2}\cos\theta)\in E_{\mathbf{R}}
\end{equation*}o\`u $\theta\in\mathbf{R}$, entra{\^\i}nent $$e^{i\theta}\cdot J(C_{\mathbf{R}})=J(C_{\mathbf{R}})\quad\text{et}\quad e^{i\theta}\cdot L(C_{\mathbf{R}})=L(C_{\mathbf{R}})\ .$$ Ainsi $J(C_{\mathbf{R}})$ et $L(C_{\mathbf{R}})$ d{\'e}finissent des ellipso{\"\i}des complexes de $E$, c.-{\`a}-d. {\'e}gaux {\`a} la boule unit{\'e} ferm{\'e}e $b_{n,\mathbf{C}}^{2}$ de l'espace hermitien usuel $(\mathbf{C}^{n},\vert.\vert_{2})$, apr{\`e}s le choix d'une base convenable de $E$. Comme dans le cas r{\'e}el, on les note plus simplement $J(C)$ et $L(C)$. L'hypoth{\`e}se~\eqref{hypothesejaugecomplexe} assure en r{\'e}alit{\'e} que tout se passe comme si l'on raisonnait dans $E_{\mathbf{R}}\simeq\mathbf{R}^{2n}$. Le quotient volumique d'un $\mathbf{C}$-espace vectoriel norm{\'e} $(E,\Vert.\Vert)$ est d{\'e}fini par la formule \begin{equation*}\vr(E)=\inf{\left\{\left(\frac{\vol(\mathbf{B}(E,\Vert.\Vert))}{\vol(D)}\right)^{\frac{1}{2n}}\,;\quad\text{$D$ ellipso{\"\i}de complexe $\subseteq E$}\right\}}\end{equation*}o\`u $\vol$ est une mesure de Haar sur $E$. Comme dans le cas r{\'e}el, cette quantit{\'e} est atteinte pour l'ellipso{\"\i}de de John $J(C)$. De m{\^e}me la d{\'e}finition de la distance de Banach-Mazur s'{\'e}tend au cas complexe. Et l'on dispose de l'encadrement \begin{equation}\label{estimationslocalescomplexes}1\le\vr(E)\vr(E^{\mathsf{v}})\le\mathrm{d}(E,\ell^{2}_{n,\mathbf{C}})\le\sqrt{2n}=\sqrt{\dim_{\mathbf{R}}E}\ .\end{equation} 
\section{Fibr{\'e} vectoriel ad{\'e}lique}
Suivant en partie R.~Rumely \textit{et al.}~\cite{Rumelyetal}, nous d{\'e}finissons ici une notion de fibr{\'e} vectoriel ad{\'e}lique sur $\spec k$, qui g{\'e}n{\'e}ralise celle de fibr{\'e} vectoriel hermitien sur $\spec\mathcal{O}_{k}$, qui est \`a la base m{\^e}me de la g\'eom\'etrie d'Arakelov.\par Soit $k$ un corps global et $v$ une place de $k$. On note $\mathbf{C}_{v}$ la compl\'etion d'une cl{\^o}ture alg\'ebrique de $k_{v}$. Si $K$ est une extension finie de $k$ et $w$ une place de $K$ au-dessus de $v$, on a un isomorphisme topologique de corps valu\'es $\mathbf{C}_{w}\simeq\mathbf{C}_{v}$. Soit $E$ un $k$-espace vectoriel. Une \emph{norme} sur $E\otimes_{k}\mathbf{C}_{v}$ est une application $\Vert.\Vert_{v}:E\otimes\mathbf{C}_{v}\to\mathbf{R}^{+}$ qui satisfait aux trois conditions~:\begin{enumerate}\item[(i)] $\forall x\in E\otimes\mathbf{C}_{v}$, $\Vert x\Vert_{v}=0\iff x=0$.\item[(ii)] $\forall x\in E\otimes\mathbf{C}_{v}$, $\forall\lambda\in\mathbf{C}_{v}$, $\Vert\lambda x\Vert_{v}=\vert\lambda\vert_{v}\cdotp\Vert x\Vert_{v}$.\item[(iii)] $\forall x,y\in E\otimes\mathbf{C}_{v}$, $\Vert x+y\Vert_{v}\le\Vert x\Vert_{v}+\Vert y\Vert_{v}$.\end{enumerate} 
\begin{defi}\label{definition31} 
Un \emph{fibr{\'e} vectoriel ad{\'e}lique} $\overline{E}=(E,(\Vert.\Vert_{v})_{v})$ sur $\spec k$ est la donn{\'e}e d'un $k$-espace vectoriel $E$ de dimension finie $n$ et d'une famille de normes $\Vert.\Vert_{v}$ sur $E\otimes_{k}\mathbf{C}_{v}$, aux places $v$ de $k$, soumise aux contraintes suivantes :\noignorespaces
\begin{enumerate}
\item[1)]
Il existe une $k$-base $(e_{1},\ldots,e_{n})$ de $E$ telle que, pour toute place $v$ ultram{\'e}trique en dehors d'un nombre fini, la norme sur $E\otimes_{k}\mathbf{C}_{v}$ est donn{\'e}e par \begin{equation}\label{definormun}\left\Vert\sum_{i=1}^{n}{x_{i}e_{i}}\right\Vert_{v}=\max_{1\le i\le n}{\{\vert x_{i}\vert_{v}\}}\ \cdotp\end{equation}
\item[2)] Soit $\gal(\mathbf{C}_{v}\vert k_{v})$ l'ensemble des automorphismes continus qui laissent invariants les {\'e}l{\'e}ments de $k_{v}$. Alors $\Vert.\Vert_{v}$ est invariante sous l'action de $\gal(\mathbf{C}_{v}\vert k_{v})$ : {\'e}tant donn{\'e} une $k_{v}$-base $(\alpha_{1},\ldots,\alpha_{n})$ de $E\otimes_{k}k_{v}$ et $(x_{1},\ldots,x_{n})\in\mathbf{C}_{v}^{n}$, $\sigma\in\gal(\mathbf{C}_{v}\vert k_{v})$, on a $$\Vert\sigma(x_{1})\alpha_{1}+\ldots+\sigma(x_{n})\alpha_{n}\Vert_{v}=\Vert x_{1}\alpha_{1}+\cdots+x_{n}\alpha_{n}\Vert_{v}\ .$$
\item[3)] Si $v$ est ultram{\'e}trique alors $$\forall\,x,y\in E\otimes_{k}\mathbf{C}_{v},\quad \Vert x+y\Vert_{v}\le\max{\{\Vert x\Vert_{v},\Vert y\Vert_{v}\}}$$(ultra-norme selon la terminologie de Bourbaki).
\end{enumerate}
Un \emph{fibr{\'e} en droites ad{\'e}lique} est un fibr{\'e} vectoriel ad{\'e}lique de dimension $1$. Un \emph{fibr{\'e} ad{\'e}lique hermitien} est un fibr{\'e} vectoriel ad{\'e}lique dont toutes les normes aux places archim{\'e}diennes de $k$ sont hermitiennes. Par extension, nous parlerons encore de fibr{\'e} ad{\'e}lique hermitien lorsque $k$ est un corps de fonctions. 
\end{defi}
Cette d{\'e}finition appelle quelques commentaires. Tout d'abord, rappelons qu'un \emph{fibr{\'e} vectoriel hermitien} sur $\spec\mathcal{O}_{k}$ ($k$ corps de nombres n{\'e}cessairement) est la donn{\'e}e d'un $\mathcal{O}_{k}$-module projectif de type fini $\mathcal{E}$ et de normes $\Vert.\Vert_{v}$ euclidiennes (\emph{resp}. hermitiennes) aux places archim{\'e}diennes r{\'e}elles (\emph{resp}. complexes) de $k$ sur les espaces $\mathcal{E}\otimes_{\mathcal{O}_{k}}k_{v}$. De plus, si $v$ est complexe, la norme $\Vert.\Vert_{v}$ est suppos{\'e}e invariante par conjugaison complexe. Cette hypoth{\`e}se suppl{\'e}mentaire correspond exactement {\`a} la condition 2) ci-dessus. La coh{\'e}rence avec la notion de fibr{\'e} vectoriel ad{\'e}lique sur $\spec k$ est assur{\'e}e gr{\^a}ce aux deux observations suivantes. D'une part, le module $\mathcal{E}$ fournit naturellement une \emph{structure enti{\`e}re} de $E=\mathcal{E}\otimes k$ au sens o\`u la norme en une place ultram\'etrique $v$ de $k$ que l'on choisit sur $E\otimes_{k}\mathbf{C}_{v}$ est donn\'ee par \begin{equation}\label{formulestructureentiere}\forall x\in E\otimes_{k}\mathbf{C}_{v},\quad\Vert x\Vert_{E,v}:=\inf{\left\{\vert a\vert_{v};\ a\in\mathbf{C}_{v},\ x\in a.\left(\mathcal{E}\otimes_{\mathcal{O}_{k}}\widehat{\mathcal{O}}_{v}\right)\right\}}\end{equation}o\`u l'anneau $\widehat{\mathcal{O}}_{v}$ est l'anneau de valuation de $\mathbf{C}_{v}$ (c.-{\`a}-d. sa boule unit\'e ferm\'ee). En consid\'erant une famille g\'en\'eratrice minimale de $\mathcal{E}$ sur $\mathcal{O}_{k}$, on v\'erifie que cette norme satisfait \`a la formule~\eqref{definormun}, et en particulier elle est invariante sous l'action de $\gal(\mathbf{C}_{v}/k_{v})$. D'autre part, une m{\'e}trique euclidienne sur $\mathcal{E}\otimes$ $k_{v}$ ($v$ r{\'e}elle) se prolonge naturellement et de mani{\`e}re unique en une m{\'e}trique \emph{hermitienne} sur $\mathcal{E}\otimes\mathbf{C}_{v}$, car une telle m{\'e}trique est d{\'e}termin{\'e}e par une matrice r{\'e}elle sym{\'e}trique d{\'e}finie positive, unique {\`a} conjugaison pr{\`e}s par les {\'e}l{\'e}ments du groupe orthogonal $\mathrm{O}_{n}(\mathbf{R})$. Ce dernier point est remarquable. L'unicit{\'e} du prolongement n'est plus vraie en g{\'e}n{\'e}ral lorsque $\mathcal{E}\otimes k_{v}$ est muni d'une norme quelconque. Il n'existe pas de proc{\'e}d{\'e} \og canonique\fg\ qui permette d'{\'e}tendre une norme $\Vert.\Vert$ sur une espace vectoriel r{\'e}el $E\otimes\mathbf{R}$ au complexifi{\'e} $E\otimes\mathbf{C}\simeq E\otimes\mathbf{R}\oplus i E\otimes\mathbf{R}$. Consid{\'e}rons par exemple un {\'e}l{\'e}ment $p\in[1,+\infty]$ et posons, pour $a=a_{1}+ia_{2}\in E\otimes\mathbf{C}$, \begin{equation*}\Vert a\Vert_{p}^{\widetilde{}}:=\begin{cases}2^{\min{\left\{\frac{1}{2}-\frac{1}{p},0\right\}}}\left(\Vert a_{1}\Vert^{p}+\Vert a_{2}\Vert^{p}\right)^{1/p} & \text{si $p<+\infty$},\\ \max{\left\{\Vert a_{1}\Vert,\Vert a_{2}\Vert\right\}} & \text{si $p=+\infty$}.
\end{cases}
\end{equation*}La fonction $a\mapsto\Vert a\Vert_{p}^{\widetilde{}}$ ne d{\'e}finit pas en g{\'e}n{\'e}ral une norme car $\Vert\lambda a\Vert_{p}^{\widetilde{}}$ peut {\^e}tre diff{\'e}rent de $\vert\lambda\vert\cdotp\Vert a\Vert_{p}^{\widetilde{}}$ pour $\lambda\in\mathbf{C}$. En revanche, si l'on d{\'e}finit \begin{equation*}\Vert a\Vert_{p}^{\#}:=\sup{\left\{\Vert e^{i\theta}a\Vert_{p}^{\widetilde{}};\ \theta\in[0,2\pi]\right\}},\end{equation*}on v{\'e}rifie alors que $\Vert.\Vert_{p}^{\#}$ est une norme sur $E\otimes\mathbf{C}$, invariante par conjugaison complexe et {\'e}gale {\`a} $\Vert.\Vert$ sur $E\otimes\mathbf{R}$ (\footnote{Cette derni{\`e}re propri{\'e}t{\'e} est la raison pour laquelle le facteur de normalisation $ 2^{\min{\left\{\frac{1}{2}-\frac{1}{p},0\right\}}}$ appara{\^\i}t dans la d{\'e}finition de $\Vert a\Vert_{p}^{\#}$, ce facteur {\'e}tant {\'e}gal {\`a} $\underset{0\le\theta\le 2\pi}{\mathrm{inf}}\{(\cos\theta)^{p}+(\sin\theta)^{p}\}^{-1/p}$.}). Si $p=2$ et si $\Vert.\Vert$ est euclidienne alors $\Vert.\Vert_{2}^{\#}$ est la norme hermitienne qui prolonge $\Vert.\Vert$ {\`a} $E\otimes\mathbf{C}$ (voir l'exemple $1.3$, p.~$16$, de~\cite{Rumelyetal}). Cette observation --- absence d'\emph{un} prolongement naturel --- est la raison essentielle qui justifie que les normes associ{\'e}es {\`a} un fibr{\'e} vectoriel ad{\'e}lique soient d{\'e}finies sur $E\otimes\mathbf{C}_{v}$ (ou $E\otimes\overline{k}_{v}$, ce qui revient au m{\^e}me quitte {\`a} prolonger ensuite par continuit{\'e}, de mani{\`e}re unique, {\`a} $E\otimes\mathbf{C}_{v}$). Il est alors possible d'op{\'e}rer une extension des scalaires en consid{\'e}rant une extension finie $K\vert k$ et le fibr{\'e} vectoriel ad{\'e}lique $\overline{E}_{K}$ dont l'espace sous-jacent est $E\otimes K$, avec les m\^emes normes que $\overline{E}$.\par Par ailleurs, la premi{\`e}re condition dans la d{\'e}finition~\ref{definition31}, si elle est remplie, est alors vraie pour \emph{toute} $k$-base de $E$, quitte {\'e}ventuellement {\`a} accro{\^\i}tre le nombre de places $v$ qui ne conviennent pas. En r{\'e}alit{\'e}, si l'on fixe une $k$-base $(e_{1},\ldots,e_{n})$ de $E$, il existe une matrice d'ad{\`e}les finis\label{adelesfinis} $(c_{v})_{v}\in\GL_{n}(k_{\mathbf{A},f})$ telle que, pour \emph{toute} place ultram{\'e}trique $v$ de $k$, on a \begin{equation}\label{normesegales}\forall x\in\mathbf{C}_{v}^{n},\quad\Vert x_{1}e_{1}+\cdots+x_{n}e_{n}\Vert_{v}=\underset{1\le i\le n}{\mathrm{max}}\left\{\left\vert(c_{v}.x)_{i}\right\vert_{v}\right\}\end{equation}o{\`u} $(c_{v}.x)_{i}$ d{\'e}signe la $i^{\text{{\`e}me}}$ coordonn{\'e}e du vecteur $c_{v}.x$. En effet, d'apr{\`e}s la proposition~3 du chapitre II, p.~26, de~\cite{WeilBasic}, et {\'e}tant donn{\'e} une place ultram{\'e}trique $v$, il existe une $k_{v}$-base $(\alpha_{1},\ldots,\alpha_{n})$ de $E\otimes k_{v}$ telle que \begin{equation}\label{formulenormesdeux}\forall x\in k_{v}^{n},\quad \Vert x_{1}\alpha_{1}+\cdots+x_{n}\alpha_{n}\Vert_{v}=\underset{1\le i\le n}{\mathrm{max}}\left\{\vert x_{i}\vert_{v}\right\}\ \cdotp\end{equation}Autrement dit, il existe $c_{v}\in\GL_{n}(k_{v})$ tel que la relation~\eqref{normesegales} soit satisfaite pour tout $x\in k_{v}^{n}$. La matrice $c_{v}$ est la matrice identit{\'e} pour presque tout $v$. Pour obtenir l'{\'e}galit{\'e}~\eqref{normesegales} avec $x\in\mathbf{C}_{v}^{n}$, on observe que si $K_{w}$ est une extension finie de $k_{v}$, il existe une matrice $c_{w}\in\GL_{n}(K_{w})$ telle que \begin{equation}\label{formulenormestrois}\forall x\in K_{w}^{n},\quad\Vert x_{1}e_{1}+\cdots+x_{n}e_{n}\Vert_{v}=\underset{1\le i\le n}{\mathrm{max}}\left\{\left\vert(c_{w}.x)_{i}\right\vert_{v}\right\}\ \cdotp\end{equation}En se restreignant {\`a} $x\in k_{v}^{n}$ et en choisissant les vecteurs de la base canonique de $k_{v}^{n}$, on constate que $c_{w}c_{v}^{-1}\in\GL_{n}(\mathcal{O}_{w})$ ($\mathcal{O}_{w}$ est l'anneau de valuation de $K_{w}$) et l'{\'e}galit{\'e}~\eqref{formulenormestrois} reste vraie en rempla{\c c}ant $c_{w}$ par $c_{v}$. Elle est donc valide pour $x\in\overline{k}_{v}^{n}$ et, par continuit{\'e}, pour $x\in\mathbf{C}_{v}^{n}$.\par Si, comme nous l'avons vu plus haut, un fibr\'e vectoriel hermitien $\overline{\mathcal{E}}$ sur $\spec\mathcal{O}_{k}$ fournit une structure enti\`ere pour $E=\mathcal{E}\otimes k$, la r\'eciproque est \'egalement vraie. Si l'on se donne un fibr\'e vectoriel ad\'elique $\overline{E}$ sur $\spec k$, l'ensemble \begin{equation*}\mathcal{E}:=\left\{x\in E\,;\ \forall v\nmid\infty,\ \Vert x\Vert_{E,v}\le 1\right\}\end{equation*}est une structure enti\`ere pour $E$. Dans cette \'ecriture, si $k$ est un corps de fonctions, les places $v$ exclues sont celles qui sont au-dessus de la place $\infty$ de $\mathbf{F}_{p}[T]$ d\'efinie dans les pr\'eliminaires. Notons $\mathcal{O}_{k}:=\{x\in k\,;\ \forall v\nmid\infty,\ \vert x\vert_{v}\le 1\}$ l'anneau des entiers de $k$. Cet ensemble est aussi la fermeture int\'egrale de l'anneau $$\mathcal{O}_{k_{0}}:=\begin{cases}\mathbf{Z} & \text{si $k$ est un corps de nombres}\\ \mathbf{F}_{p}[T] & \text{si $k$ est un corps de fonctions}\end{cases}$$(voir~\cite{bourbakialgcomsix}, chap.~VI, \S~$3$, corollaire~$3$), et {\`a} ce titre il est de Dedekind. L'ensemble $\mathcal{E}$ est un module sans torsion sur $\mathcal{O}_{k}$. Vu la caract\'erisation~\eqref{normesegales} des normes de $\overline{E}$ aux places ultram\'etriques, il existe un \'el\'ement $N\in\mathcal{O}_{k_{0}}\setminus\{0\}$ tel que $NE\subseteq\mathcal{E}$, ce qui entra{\^\i}ne que $\mathcal{E}$ est de type fini. Cela montre que $\mathcal{E}$ est projectif de type fini (voir~\cite{langalgnombres}, chap.~$1$, \S~$9$) et on v\'erifie que, pour toute place $v$ de $k$, qui ne domine pas la place $\infty$, la norme sur $E\otimes_{k}\mathbf{C}_{v}$ induite par $\mathcal{E}$ au moyen de la formule~\eqref{formulestructureentiere} est \'egale \`a la norme $\Vert.\Vert_{E,v}$ de d\'epart.  \par Enfin, signalons que certains auteurs accordent (ou accorderaient) aux normes $\Vert.\Vert_{v}$ de n'{\^e}tre que des semi-normes. Ici cela semble \textit{a priori} exclu en partie {\`a} cause des r{\'e}sultats du paragraphe pr{\'e}c{\'e}dent dans les espaces norm{\'e}s dont nous aurons besoin.
\subsection{Exemples de fibr{\'e}s vectoriels ad{\'e}liques}\label{fvap}L'exemple le plus simple est celui de l'espace $k$ lui-m\^eme, qui, muni des diff{\'e}rentes valeurs absolues $\vert.\vert_{v}$, poss{\`e}de une structure de fibr{\'e} vectoriel ad{\'e}lique (dite \og triviale\fg). D'autres exemples sont les fibr{\'e}s vectoriels ad{\'e}liques que l'on va noter $(k^{n},\vert.\vert_{p})$ avec $p\in[1,+\infty]$, dont la structure enti{\`e}re est donn{\'e}e par la base canonique de $k^{n}$ et o{\`u}, en une place $v$ archim{\'e}dienne, la norme $\vert.\vert_{p}$ est celle d{\'e}finie par la formule~\eqref{defidenormep}. L'on pourrait bien s{\^u}r panacher plusieurs normes de ce type aux diff{\'e}rentes places archim{\'e}diennes ou changer la structure enti{\`e}re en choisissant une autre base de $k^{n}$. Mais au-del{\`a} de ces exemples et {\`a} l'instar des fibr{\'e}s vectoriels hermitiens sur $\spec\mathcal{O}_{k}$, il existe d'autres exemples provenant de la g{\'e}om{\'e}trie alg{\'e}brique, obtenus en consid{\'e}rant l'espace des sections globales d'un fibr{\'e} en droites m\'etris\'e ad{\'e}liquement sur une vari{\'e}t{\'e} projective au-dessus de $k$ (voir \S~\ref{par:conclusion}). 
\subsection{Lien avec la notion de convexe ad{\'e}lique} {\'E}tant donn{\'e} une place ultram{\'e}trique $v$ de $k$, une partie $C_{v}$ de $E\otimes k_{v}$ est appel{\'e} $k_{v}$-r{\'e}seau si $C_{v}$ est un $\mathcal{O}_{v}$-sous-module de $E\otimes k_{v}$, {\`a} la fois ouvert et compact. Fixons une $k$-base de $E$ qui permet d'identifier $E$ {\`a} $k^{n}$. Dans toute la suite, nous \emph{supposons} que le $k_{v}$-r{\'e}seau $C_{v}$ s'identifie {\`a} $\mathcal{O}_{v}^{n}$, pour toute place $v$ en dehors d'un nombre fini. Si $v$ est une place archim{\'e}dienne de $k$, soit $C_{v}$ un sous-ensemble convexe de $E\otimes k_{v}$, d'int{\'e}rieur non vide, compact et sym{\'e}trique par rapport {\`a} l'origine. L'origine est un point int{\'e}rieur de $C_{v}$. Si $v$ est une place complexe, on suppose de plus que \begin{equation}\label{hypothesejaugecomplexedeux}\forall\theta\in\mathbf{R},\quad e^{i\theta}.C_{v}=C_{v}\ .\end{equation}
\begin{defi}
Un \emph{convexe ad{\'e}lique} est un sous-ensemble $C$ de $E_{\mathbf{A}}=E\otimes k_{\mathbf{A}}$ de la forme $\prod_{v}{C_{v}}$, o\`u les ensembles $C_{v}$ v{\'e}rifient les hypoth{\`e}ses ci-dessus.
\end{defi} 
Si $\overline{E}$ est un fibr{\'e} vectoriel ad{\'e}lique sur $\spec k$, la boule unit{\'e} ferm{\'e}e de $\overline{E}$ $$\mathbb{B}(\overline{E})=\left\{(x_{v})_{v}\in E_{\mathbf{A}}\,;\ \forall\,v,\ \Vert x_{v}\Vert_{v}\le 1\right\}$$est un convexe ad{\'e}lique. En effet, l'ensemble $C_{v}=\{x_{v}\in E\otimes k_{v}\,;\ \Vert x_{v}\Vert_{v}\le 1\}$ convient dans tous les cas et la condition 1) de la d{\'e}finition~\ref{definition31} assure que $C_{v}\simeq\mathcal{O}_{v}^{n}$ pour presque tout $v$. On a {\'e}galement une r{\'e}ciproque :
\begin{prop}\label{proposition33}Tout convexe ad{\'e}lique sur $k$ est la boule unit{\'e} ferm{\'e}e d'un fibr{\'e} vectoriel ad{\'e}lique sur $\spec k$, non n{\'e}cessairement unique.
\end{prop}   
\begin{proof}
  Si $v$ est une place archim{\'e}dienne, l'on sait que la donn{\'e}e d'un ensemble convexe $C_{v}$ comme ci-dessus {\'e}quivaut {\`a} la donn{\'e}e d'une norme sur $E\otimes k_{v}$ dont $C_{v}$ est la boule unit{\'e} correspondante, norme donn{\'e}e explicitement par l'application jauge $j_{v}:E\otimes k_{v}\to\mathbf{R}^{+}$ :\begin{equation*}j_{v}(x)=\inf{\{\lambda>0\,;\ \frac{x}{\lambda}\in C_{v}\}}\cdotp\end{equation*}La norme ainsi construite se prolonge {\`a} $E\otimes\mathbf{C}_{v}$ (sans que le prolongement ne soit unique). Comme nous l'avons d{\'e}j{\`a} mentionn{\'e}, l'hypoth{\`e}se~\eqref{hypothesejaugecomplexedeux} intervient dans le cas complexe pour assurer l'{\'e}galit{\'e} $j_{v}(ax)=\vert a\vert j_{v}(x)$ pour tous $a\in\mathbf{C}$ et $x\in E\otimes_{v}\mathbf{C}$.\par Si $v$ est une place ultram{\'e}trique, soit $\pi_{v}$ une uniformisante de $\mathcal{O}_{v}$ (c.-{\`a}-d. un g{\'e}n{\'e}rateur de l'id{\'e}al maximal principal de $\mathcal{O}_{v}$). Et, comme dans le cas pr{\'e}c{\'e}dent, consid{\'e}rons pour $x\in E\otimes k_{v}$ le nombre r{\'e}el \begin{equation*}j_{v}(x)=\inf\left\{\lambda=\vert\pi_{v}\vert_{v}^{h}\,;\ h\in\mathbf{Z}\ \text{et}\ \frac{x}{\pi_{v}^{h}}\in C_{v}\right\}\cdotp\end{equation*}Ce nombre est bien d{\'e}fini car $C_{v}$ est ouvert et la borne inf{\'e}rieure est un minimum car $C_{v}$ est compact. On v{\'e}rifie alors que $j_{v}$ d{\'e}finit une norme ultram{\'e}trique sur $E\otimes k_{v}$ dont $C_{v}$ est la boule unit{\'e} ferm{\'e}e. Comme nous l'avons vu plus haut, cette norme s'{\'e}tend {\`a} $\mathbf{C}_{v}^{n}$ \textit{via} le choix d'une matrice $c_{v}\in\GL_{n}(k_{v})$ convenable. L'hypoth{\`e}se $C_{v}\simeq\mathcal{O}_{v}^{n}$ permet de prendre $c_{v}$ {\'e}gal {\`a} la matrice identit{\'e} pour presque tout $v$. La collection des normes $(j_{v})_{v}$, prolong{\'e}es aux espaces $E\otimes\mathbf{C}_{v}$, fournit la structure ad{\'e}lique recherch{\'e}e. 
\end{proof}

\begin{rema}
Si $E$ est une droite, il n'y a plus d'ambigu{\"\i}t{\'e} selon que l'on regarde $E\otimes k_{v}$ ou $E\otimes\mathbf{C}_{v}$ ($v$ archim{\'e}dienne)\,; les m{\'e}triques aux {\'e}ventuelles places archim{\'e}diennes de $k$ sont alors automatiquement hermitiennes.
\end{rema}
 Il arrive souvent que les normes aux places ultram{\'e}triques proviennent \emph{toutes} d'une m{\^e}me base $(e_{1},\ldots,e_{n})$ de $E$, comme dans la formule~\eqref{definormun} de la d{\'e}finition~\ref{definition31}. Cette base fournit la \emph{structure enti{\`e}re} du fibr{\'e} vectoriel ad{\'e}lique $\overline{E}$. Lorsque $E=k^{n}$, la base canonique donne une structure enti{\`e}re naturelle, utilis{\'e}e par d{\'e}faut dans la suite. 

\subsection{Op{\'e}rations alg{\'e}briques sur l'ensemble des fibr{\'e}s vectoriels ad{\'e}liques}\label{paragraphetroistrois}
Soit $\overline{E},\overline{F}$ des fibr{\'e}s vectoriels ad{\'e}liques sur $\spec k$ de dimensions respectives $n$ et $m$.
\subsubsection*{Sous-fibr{\'e} et quotient}Nous dirons que $\overline{F}$ est un sous-fibr{\'e} de $\overline{E}$, et nous {\'e}crirons $\overline{F}\subseteq\overline{E}$, si $F\subseteq E$, et si, en chaque place $v$, la norme $\Vert.\Vert_{F,v}$ est la restriction de $\Vert.\Vert_{E,v}$ {\`a} $F\otimes_{k}\mathbf{C}_{v}$. De m{\^e}me, si $\overline{F}\subseteq\overline{E}$ alors le quotient $E/F$ est muni d'une structure de fibr{\'e} vectoriel ad{\'e}lique sur $\spec k$ en consid{\'e}rant les normes quotient.
\subsubsection*{Somme directe}\label{sommedirectefibreadelique}
Il n'existe pas de norme plus naturelle qu'une autre sur la somme directe (ou le produit) d'espaces vectoriels norm{\'e}s. Afin d'assurer la compatibilit{\'e} avec la somme directe hermitienne de fibr{\'e}s vectoriels hermitiens, la norme, par d{\'e}faut, que nous choisirons sur $E\oplus F$ est, pour tous $x\in E\otimes\mathbf{C}_{v}$, $y\in F\otimes\mathbf{C}_{v}$, \begin{equation*} \Vert (x,y)\Vert_{E\oplus F,v}:=\begin{cases}\max{\{\Vert x\Vert_{E,v},\Vert y\Vert_{F,v}\}} & \text{si $v$ est ultram{\'e}trique},\\ \left(\Vert x\Vert_{E,v}^{2}+\Vert y\Vert_{F,v}^{2}\right)^{1/2} &\text{si $v$ est archim{\'e}dienne}.\end{cases}\end{equation*}Nous noterons $\overline{E}\oplus_{2}\overline{F}$ ou, plus simplement, $\overline{E}\oplus\overline{F}$, le fibr{\'e} vectoriel ad{\'e}lique obtenu de la sorte. Si l'on remplace la norme $\vert.\vert_{2}$ par la norme $\vert.\vert_{p}$, $p\in[1,+\infty]$, nous noterons $\overline{E}\oplus_{p}\overline{F}$ le fibr{\'e} vectoriel ad{\'e}lique que l'on obtient. Il peut {\^e}tre utile parfois d'avoir (encore) un peu plus de souplesse dans le choix des normes. \'Etant donn{\'e} une place archim{\'e}dienne $v$ de $k$, consid{\'e}rons $\varsigma_{v}$ une norme sym{\'e}trique sur $\mathbf{R}^{2}$, invariante par changement de signes sur les coordonn{\'e}es, telle que $\varsigma_{v}(1,0)=\varsigma_{v}(0,1)=1$. Posons $\varsigma:=(\varsigma_{v})_{v\mid\infty}$. Le fibr{\'e} vectoriel ad{\'e}lique $\overline{E}\oplus_{\varsigma}\overline{F}$ est par d{\'e}finition l'espace $E\oplus F$ muni, aux places archim{\'e}diennes, des normes $\Vert (x,y)\Vert_{E\oplus_{\varsigma}F,v}:=\varsigma_{v}(\Vert x\Vert_{E,v},\Vert y\Vert_{F,v})$. Les normes aux places ultram{\'e}triques sont les m{\^e}mes que pr{\'e}c{\'e}demment.
\begin{rema}
Ces op{\'e}rations pr{\'e}servent la structure hermitienne quand les objets de d{\'e}part en sont pourvus.
\end{rema} 
\subsubsection*{Dual et norme d'op\'erateur}
 Le fibr{\'e} vectoriel dual $\overline{E}^{\mathsf{v}}$ est l'espace dual $E^{\mathsf{v}}=\Hom_{k}(E,k)$, muni des normes duales usuelles :\begin{equation*}\forall\varphi\in(E\otimes\mathbf{C}_{v})^{\mathsf{v}},\quad\Vert\varphi\Vert_{E^{\mathsf{v}},v}:=\sup{\left\{\frac{\vert\varphi(x)\vert_{v}}{\Vert x\Vert_{E,v}};\ x\in E\otimes\mathbf{C}_{v},\ x\ne 0\right\}}\ \cdotp\end{equation*}Plus g{\'e}n{\'e}ralement, l'espace $\Hom_{k}(E,F)$ des applications lin{\'e}aires entre $E$ et $F$ est muni aux diff{\'e}rentes places $v$ des normes d'op{\'e}rateurs usuelles, obtenues en rempla{\c c}ant $\vert\varphi(x)\vert_{v}$ par $\Vert\varphi(x)\Vert_{F,v}$ dans l'expression ci-dessus. D'une mani{\`e}re alternative, on a, pour tout $\varphi\in\Hom_{k}(E,F)\otimes\mathbf{C}_{v}$,\begin{equation*}\Vert\varphi\Vert_{\Hom(E,F),v}=\sup{\left\{y^{\mathsf{v}}(\varphi(x))\,;\ \Vert y^{\mathsf{v}}\Vert_{F^{\mathsf{v}},v}=\Vert x\Vert_{E,v}=1\right\}}\ \cdot\end{equation*}
\begin{rema}
Le bidual $(\overline{E^{\mathsf{v}}})^{\mathsf{v}}$ est isomorphe isom\'etriquement \`a $\overline{E}$ et l'on dispose de la formule (vraie pour toute place $v$ de $k$) \begin{equation}\label{normebiduale}\forall\,x\in E\otimes_{k}\mathbf{C}_{v},\quad \Vert x\Vert_{E,v}=\sup_{\Vert\varphi\Vert_{E^{\mathsf{v}},v}\le 1}{\vert\varphi(x)\vert_{v}}\ .\end{equation}Dans le cas archim\'edien, il s'agit d'un résultat classique de la th\'eorie des espaces vectoriels norm\'es de dimension finie (voir~\cite{BourbakiHilbert}, chap.~IV, \S~2.4). Dans le cas ultram\'etrique, l'\'egalit\'e~\eqref{normebiduale} peut se d\'emontrer au moyen de la base $(\alpha_{1},\ldots,\alpha_{n})$ de $E\otimes\mathbf{C}_{v}$ donnant la formule~\eqref{formulenormesdeux}.
\end{rema}
\subsubsection*{Produit tensoriel}La norme d'op\'erateur conf{\`e}re au produit tensoriel $E\otimes_{k}F$ une structure de fibr{\'e} vectoriel ad{\'e}lique \textit{via} l'isomorphisme $E\otimes_{k}F\simeq\Hom_{k}(E^{\mathsf{v}},F)$. Malheureusement, la norme tensorielle ainsi b{\^a}tie sur $E\otimes_{k}F$ ne donne pas en g\'en\'eral une structure hermitienne \`a $E\otimes_{k}F$ lorsque $\overline{E}$ et $\overline{F}$ sont des fibr\'es vectoriels \emph{hermitiens} (sauf si $E$ ou $F$ est une droite, ou si $k$ est un corps de fonctions).
\begin{defi}\label{defidenormestensoriellesadeliques} Soit $\ell\in\mathbf{N}\setminus\{0\}$. Une \emph{norme tensorielle ad\'elique d'ordre $\ell$} sur $k$ (corps global) est la donn\'ee d'une famille $\alpha=(\alpha_{v})_{v}$ index\'ee par les places $v$ de $k$ telle que~:\begin{enumerate}\item[1)] Si $v$ est archim\'edienne, $\alpha_{v}$ est une norme tensorielle d'ordre $\ell$ (au sens de la d\'efinition~\ref{defidenormestensoriellesreelles}).\item[2)] Si $v$ est ultram\'etrique alors $\alpha_{v}$ est la norme tensorielle ultram\'etrique d\'efinie de la mani\`ere suivante~: pour tout $i\in\{1,\ldots,\ell\}$, pour tout $\mathbf{C}_{v}$-espace vectoriel ultra-norm\'e $(E_{i},\Vert.\Vert_{E_{i}})$ de dimension finie $n_{i}$, la norme $\alpha_{v}(\cdot;E_{1},\ldots,E_{\ell})$ sur $\otimes_{i=1}^{\ell}{E_{i}}$ est l'ultra-norme d\'efinie par r\'ecurrence sur $\ell$~: $\alpha_{v}(\cdot;E_{1})=\Vert.\Vert_{E_{1}}$ et $\alpha_{v}(\cdot;E_{1},\ldots,E_{\ell})$ s'identifie \`a la norme d'op\'erateur des applications lin\'eaires entre l'espace norm\'e dual $(E_{1}^{\mathsf{v}},\Vert.\Vert_{E_{1}^{\mathsf{v}}})$ et l'espace norm\'e $(E_{2}\otimes\cdots\otimes E_{\ell},\alpha_{v}(\cdot;E_{2},\ldots,E_{\ell}))$.
\end{enumerate}
Une \emph{norme tensorielle (ad\'elique) hermitienne} est la donn\'ee d'une norme tensorielle ad\'elique $\alpha=(\alpha_{v})_{v}$ telle que, pour toute place archim\'edienne $v$ de $k$, la norme tensorielle $\alpha_{v}$ est hermitienne au sens de la d\'efinition~\ref{defidenormestensoriellesreelles}.
\end{defi}
\'Etant donn\'e des fibr\'es vectoriels ad\'eliques $\overline{E}_{1},\ldots,\overline{E}_{\ell}$ sur $\spec k$, on notera $\alpha(\cdot;\overline{E}_{1},\ldots,\overline{E}_{\ell})$ la collection des normes $\alpha_{v}(\cdot;(E_{1}\otimes_{k}\mathbf{C}_{v},\Vert.\Vert_{\overline{E}_{1},v}),\ldots,(E_{\ell}\otimes_{k}\mathbf{C}_{v},\Vert.\Vert_{\overline{E}_{\ell},v}))$ sur $E_{1}\otimes_{k}\cdots\otimes_{k}E_{\ell}\otimes_{k}\mathbf{C}_{v}$, lorsque $v$ varie parmi les places de $k$. La d\'efinition~\ref{defidenormestensoriellesadeliques} assure que si $\alpha$ est une norme tensorielle ad\'elique, alors le couple $$\otimes_{\alpha,i=1}^{\ell}\overline{E}_{i}:=(E_{1}\otimes_{k}\cdots\otimes_{k}E_{\ell},\alpha(\cdot;\overline{E}_{1},\ldots,\overline{E}_{\ell}))$$est un fibr{\'e} vectoriel ad\'elique sur $\spec k$ (lorsque tous les $\overline{E}_{i}$ sont \'egaux \`a $\overline{E}$, on notera aussi $\overline{E}^{\otimes_{\alpha}\ell}:=\otimes_{\alpha,i=1}^{\ell}\overline{E}$). Et, de plus, si tous les $\overline{E}_{i}$ et $\alpha$ sont hermitiens alors $\otimes_{\alpha,i=1}^{\ell}\overline{E}_{i}$ est un fibr\'e ad\'elique hermitien. En choisissant aux places archim\'ediennes de $k$ les normes de Chevet-Saphard d\'efinies dans l'exemple~\ref{normesChevetSaphard}, on obtient des exemples de normes tensorielles ad\'eliques d'ordre $2$, que l'on peut ensuite \'etendre \`a tout ordre. En particulier si $\alpha_{v}=g_{2}$ pour toute place $v$ archim\'edienne alors la norme tensorielle ad\'elique $(\alpha_{v})_{v}$, que nous noterons encore $g_{2}$, est hermitienne. 
\begin{rema}
Par r\'ecurrence sur $\ell$, la formule de bidualit\'e~\eqref{normebiduale} fournit une expression concr\`ete de la norme d'op\'erateur g\'en\'eralis\'ee du point $2)$ de la d\'efinition~\ref{defidenormestensoriellesadeliques}, sous la forme suivante~: pour tout entier $h\ge 1$, pour tous $x_{j}^{(i)}\in E_{j}$, $1\le j\le\ell$ et $1\le i\le h$, on a \begin{equation*}\begin{split}&\alpha_{v}\left(\sum_{i=1}^{h}{\otimes_{j=1}^{\ell}{x_{j}^{(i)}}};E_{1},\ldots,E_{\ell}\right) \\ &\qquad=\sup{\left\{\left\vert\sum_{i=1}^{h}{\prod_{j=1}^{\ell}{\varphi_{j}(x_{j}^{(i)})}}\right\vert_{v};\forall\,j\in\{1,\ldots,\ell\},\ \Vert\varphi_{j}\Vert_{E_{j}^{\mathsf{v}}}\le 1\right\}}\cdotp\end{split}\end{equation*}   
\end{rema}

\subsubsection*{Produit ext{\'e}rieur et d{\'e}terminant}
L'alg{\`e}bre ext{\'e}rieure $\bigwedge(E)$ est le quotient de l'alg{\`e}bre tensorielle $\mathbf{T}(E)$ par l'id{\'e}al bilat{\`e}re engendr{\'e} par les {\'e}l{\'e}ments $x\otimes x$ avec $x\in E$. Soit $\ell\in\{1,\ldots,n\}$. Le $\ell^{\text{\`eme}}$ produit ext{\'e}rieur $\bigwedge^{\ell}(E)$ est donc un quotient du produit tensoriel $E^{\otimes \ell}$. Il h{\'e}rite alors de la structure ad\'elique quotient induite par celle d{\'e}finie pr{\'e}c{\'e}demment sur $\overline{E}^{\otimes_{\alpha}\ell}$ au moyen d'une norme tensorielle ad\'elique $\alpha$ d'ordre $\ell$. En particulier, si $\ell=n$, le d{\'e}terminant $\det(E)=\bigwedge^{n}(E)$ de $E$ poss{\`e}dent une structure de fibr{\'e} vectoriel ad{\'e}lique sur $\spec k$. Nous noterons $\bigwedge_{\alpha}^{\ell}(\overline{E})$ et $\det_{\alpha}\overline{E}$ les fibr\'es ad\'eliques que l'on obtient ainsi. Si $\overline{E}$ et $\alpha$ sont hermitiens alors le fibr\'e vectoriel ad\'elique $\bigwedge_{\alpha}^{\ell}(\overline{E})$ est {\'e}galement hermitien et l'on a dans ce cas, pour toute place archim\'edienne $v$ de $k$, \begin{equation*}
\Vert e_{1}\wedge\cdots\wedge e_{\ell}\Vert_{\bigwedge_{\alpha}^{\ell}(E),v}=\det((e_{i},e_{j})_{v})_{i,j}^{1/2}\quad\text{pour tous $e_{1},\ldots,e_{\ell}\in E\otimes\mathbf{C}_{v}$}\end{equation*}($(\,,)_{v}$ d\'esigne le produit hermitien sur $E\otimes\mathbf{C}_{v}$).  
\subsubsection*{Produit sym{\'e}trique}
L'alg{\`e}bre sym{\'e}trique $\mathbf{S}(E)$ est le quotient de l'alg{\`e}bre tensorielle $\mathbf{T}(E)$ par l'id{\'e}al bilat{\`e}re engendr{\'e} par les {\'e}l{\'e}ments de la forme $x\otimes y-y\otimes x$. Pour $\ell\in\mathbf{N}$, la structure de fibr{\'e} vectoriel ad{\'e}lique sur la puissance sym\'etrique $S^{\ell}(E)$ est celle obtenue par quotient de celle de $\overline{E}^{\otimes_{\alpha}\ell}$, une norme tensorielle ad\'elique $\alpha$ d'ordre $\ell$ ayant \'et\'e fix\'ee au pr\'ealable. Nous noterons $S_{\alpha}^{\ell}(\overline{E})$ le fibr\'e vectoriel ad\'elique obtenu de la sorte. Si $\overline{E}$ et $\alpha$ sont hermitiens alors $S_{\alpha}^{\ell}(\overline{E})$ est \'egalement hermitien. Dans ce cas l'on peut pr\'eciser la valeur de la norme en une place archim\'edienne $v$ de $k$ d'un \'el\'ement de $S^{\ell}(E)\otimes_{v}\mathbf{C}$ de la mani\`ere suivante. Soit $\sigma_{E,v}:E^{\otimes\ell}\otimes_{v}\mathbf{C}\to E^{\otimes\ell}\otimes_{v}\mathbf{C}$ l'op\'erateur de sym\'etrisation~: pour tous $x_{1},\ldots,x_{\ell}\in E\otimes_{v}\mathbf{C}$, on a $\sigma_{E,v}(x_{1}\otimes\cdots\otimes x_{\ell})=\frac{1}{\ell!}\sum_{\eta\in\mathfrak{S}_{\ell}}{x_{\eta(1)}\otimes\cdots\otimes x_{\eta(\ell)}}$ ($\mathfrak{S}_{\ell}$ est le groupe sym\'etrique d'ordre $\ell$). L'image de $\sigma_{E,v}$, munie de la norme restreinte de $E^{\otimes\ell}\otimes_{v}\mathbf{C}$, est un espace hermitien, isomorphe isom\'etriquement \`a $S^{\ell}(E)\otimes_{v}\mathbf{C}$, muni de la norme quotient comme ci-dessus. De la sorte, si $e_{1},\ldots,e_{n}$ d\'esigne une base orthonorm\'ee de $E\otimes_{v}\mathbf{C}$, la famille $\{e_{1}^{i_{1}}\cdots e_{n}^{i_{n}};\ (i_{1},\ldots,i_{n})\in\mathbf{N}^{n}\ \text{et}\ i_{1}+\cdots+i_{n}=\ell\}$ forme une base orthogonale de $S^{\ell}(E)\otimes_{v}\mathbf{C}$. Les normes des vecteurs de cette base se calculent au moyen de la d\'efinition de $\sigma_{E,v}$ et l'on obtient \begin{equation*}\Vert e_{1}^{i_{1}}\cdots e_{n}^{i_{n}}\Vert_{S^{\ell}(E)\otimes_{v}\mathbf{C}}=\Vert\sigma_{E,v}(e_{1}^{\otimes i_{1}}\otimes\cdots\otimes e_{n}^{\otimes i_{n}})\Vert_{E^{\otimes\ell}\otimes_{v}\mathbf{C}}=\left(\frac{i_{1}!\cdots i_{n}!}{\ell !}\right)^{1/2}\end{equation*}(voir~\cite{BourbakiHilbert}, \S~$3.3$ du cinqui\`eme chapitre).
\begin{remas}
\begin{enumerate}\item[-] Il existe d'autres normes ad\'eliques \og naturelles\fg\ sur le produit sym\'etrique qui ne sont pas obtenues par quotient. Toutefois, une telle norme \og naturelle\fg\ se compare \`a une norme obtenue par quotient d'une norme tensorielle, avec des constantes de comparaison ne faisant intervenir que $\ell$ (voir~\cite{Floret}). \item[-] Lorsque $\overline{E}$ et $\alpha$ sont hermitiens, les fibr\'es ad\'eliques $\overline{E}^{\otimes_{\alpha}\ell}$, $\bigwedge_{\alpha}^{\ell}(\overline{E})$ et $S_{\alpha}^{\ell}(\overline{E})$ ne d\'ependent pas de $\alpha$ et nous simplifierons la notation en ne mettant pas l'indice $\alpha$. 
\end{enumerate}
\end{remas}

\section{Notion de degr{\'e} ad{\'e}lique et propri{\'e}t{\'e}s}
Soit $\overline{E}=(E,(\Vert.\Vert_{v})_{v})$ un fibr{\'e} vectoriel ad{\'e}lique sur $\spec k$, de dimension $n\ge 1$. La boule unit{\'e} de $\overline{E}$ est l'ensemble \begin{equation}\mathbb{B}(\overline{E}):=\left\{(x_{v})_{v}\in E_{\mathbf{A}}\,;\quad\text{$\forall\,v$ place de $k$},\quad\Vert x_{v}\Vert_{v}\le 1\right\}\cdotp \end{equation}
\begin{defi}Soit $\phi:E\to k^{n}$ un isomorphisme d'espaces vectoriels. Soit $\vol$ une mesure de Haar sur $k_{\mathbf{A}}^{n}$. Le \emph{degr{\'e} ad{\'e}lique} de $\overline{E}$ est le nombre r{\'e}el :
\begin{equation}
\widehat{\deg}\overline{E}:=\log\frac{\vol(\phi(\mathbb{B}(\overline{E})))}{\vol(\mathbb{B}(k^{n},\vert.\vert_{2}))}\cdotp
\end{equation} 
Le \emph{degr{\'e} ad{\'e}lique normalis{\'e}} est $$\widehat{\deg}_{\mathrm{n}}\overline{E}:=\frac{\widehat{\deg}\overline{E}}{[k:k_{0}]}\ \cdotp$$
\end{defi} 
Ce nombre est ind{\'e}pendant des choix de l'isomorphisme $\phi$, en raison de la formule du produit, et de la mesure de Haar $\vol$, car elles sont toutes proportionnelles. Par convention, si $E=\{0\}$, on pose $\widehat{\deg}\overline{E}:=0$. Pour la suite, il est utile de mentionner d{\`e}s {\`a} pr{\'e}sent le comportement de ce degr{\'e} par changement d'{\'e}chelle. Soit $A=(A_{v})_{v}\in\GL(E_{\mathbf{A}})$.\label{pagefibre} On note $A.\overline{E}=(E,A.\Vert.\Vert)$ le fibr{\'e} vectoriel ad{\'e}lique dont l'espace sous-jacent est $E$ et dont la norme en une place $v$ de $k$ est donn{\'e}e par $\Vert x\Vert_{v}:=\Vert A_{v}(x)\Vert_{E,v}$. Le d{\'e}terminant de $A$ calcul{\'e} dans une $k$-base de $E$ s'identifie {\`a} un id{\`e}le de $k_{\mathbf{A}}$ dont la valeur absolue $\vert\det A\vert_{\mathbf{A}}$ ne d{\'e}pend pas du choix de cette $k$-base. Compte tenu des normalisations du d{\'e}but de ce texte (en particulier la formule~\eqref{equation:rayonboule}), on d{\'e}duit la relation \begin{equation}\label{transfoechelle}\widehat{\deg}(E,A.\Vert.\Vert)=\widehat{\deg}\overline{E}-\log\vert\det A\vert_{\mathbf{A}}\ .\end{equation}Par ailleurs et bien que nous ne l'utiliserons pas dans la suite, signalons qu'il existe \'egalement une \emph{caract\'eristique d'Euler-Poincar\'e} d'un fibr{\'e} vectoriel ad{\'e}lique $\overline{E}$, que l'on peut d{\'e}finir en posant \begin{equation*}\chi(\overline{E}):=\log\frac{\vol(\mathbb{B}(\overline{E}))}{\covol(E)}\end{equation*}o\`u $\vol$ est une mesure de Haar quelconque sur $E_{\mathbf{A}}$ et $\covol(E)=\vol(E_{\mathbf{A}}/E)$. Le lien avec le degr{\'e} ad{\'e}lique s'exprime au travers d'une formule analogue {\`a} celle de Riemann-Roch~: \begin{equation*}\chi(\overline{E})=\widehat{\deg}\overline{E}+\chi(k^{n},\vert.\vert_{2})\end{equation*}(voir~\cite{CSouleSuccessive,SzpiroA}).
\subsubsection*{Notation} Soit $\overline{E_{1}}=(E,\Vert.\Vert_{1})$ et $\overline{E_{2}}=(E,\Vert.\Vert_{2})$ deux fibr{\'e}s vectoriels ad{\'e}liques ayant le m{\^e}me espace sous-jacent $E$. Si, en toute place $v$ de $k$ et pour tout $x\in E\otimes_{k}\mathbf{C}_{v}$, on a $\Vert x\Vert_{1,v}\le\Vert x\Vert_{2,v}$ alors nous {\'e}crirons $\overline{E_{1}}\preceq\overline{E_{2}}$.\par L'int\'er{\^e}t de cette {\'e}criture est le fait suivant.
\begin{lemma}\label{lemmesimplee}
Si $\overline{E_{1}}\preceq\overline{E_{2}}$ alors $\widehat{\deg}\overline{E_{2}}\le\widehat{\deg}\overline{E_{1}}$.
\end{lemma}Cela r{\'e}sulte de l'inclusion des boules unit{\'e}s $\mathbb{B}(\overline{E_{2}})\subseteq\mathbb{B}(\overline{E_{1}})$.
\begin{lemma}\label{lemmecomptrois}
Soit $\alpha$ une norme tensorielle ad\'elique d'ordre $\ell\in\mathbf{N}\setminus\{0\}$. Pour tout $i\in\{1,\ldots,\ell\}$, soit $\overline{E}_{i}\preceq\overline{F}_{i}$ deux fibr\'es vectoriels ad\'eliques sur $\spec k$. Alors $\otimes_{\alpha,i=1}^{\ell}\overline{E}_{i}\preceq\otimes_{\alpha,i=1}^{\ell}\overline{F}_{i}$. En particulier on a $\overline{E}_{1}^{\otimes_{\alpha}\ell}\preceq\overline{F}_{1}^{\otimes_{\alpha}\ell}$, $\bigwedge_{\alpha}^{\ell}\overline{E}_{1}\preceq\bigwedge_{\alpha}^{\ell}\overline{F}_{1}$, $\det_{\alpha}\overline{E}_{1}\preceq\det_{\alpha}\overline{F}_{1}$ (si $\ell=\dim E_{1}$) et $S_{\alpha}^{\ell}(\overline{E}_{1})\preceq S_{\alpha}^{\ell}(\overline{F}_{1})$. 
\end{lemma}Le premier r\'esultat d\'ecoule directement de la proposition~\ref{propositioncompnormes}. Et les suivants sont cons\'equences des d\'efinitions. \par Ces observations seront souvent utilis\'ees dans la suite.
\subsection{Exemples} Au paragraphe~\ref{fvap}, nous avons introduit le fibr{\'e} vectoriel ad{\'e}lique $(k^{n},\vert.\vert_{p})$, $p\in[1,+\infty]$. Si $k$ est un corps de fonctions alors le degr{\'e} ad{\'e}lique de $(k^{n},\vert.\vert_{p})$ est nul. Si $k$ est un corps de nombres, les calculs~\eqref{formuledevolumep} et~\eqref{formulevolumecomplexep} des volumes des boules unit{\'e}s r{\'e}elle et complexe donnent la formule exacte~:\begin{equation}\label{formuleexactedegrekp}\widehat{\deg}(k^{n},\vert.\vert_{p})=r_{1}\log\frac{\Gamma\left(1+\frac{1}{p}\right)^{n}\Gamma\left(1+\frac{n}{2}\right)}{\Gamma\left(1+\frac{n}{p}\right)(\sqrt{\pi})^{n}}+r_{2}\log\frac{\Gamma\left(1+\frac{2}{p}\right)^{n}n!}{\Gamma\left(1+\frac{2n}{p}\right)},\end{equation}o\`u $r_1$ et $r_2$ sont respectivement le nombre de places r{\'e}elles et complexes de $k$. Lorsque $n\longrightarrow+\infty$ et $p$ est fix{\'e}, la formule de Stirling, qui se trouve en note au bas de la page~\pageref{Stirling}, fournit le d{\'e}veloppement asymptotique \begin{equation}\label{formuleasymptotiquedegre}\widehat{\deg}(k^{n},\vert.\vert_{p})=an\log(n)+bn+c+\mathrm{o}(1)\end{equation}avec\begin{equation*}\begin{split}a:=& \left(\frac{1}{2}-\frac{1}{p}\right)[k:\mathbf{Q}],\\ b:=& r_{1}\log\frac{\Gamma\left(1+\frac{1}{p}\right)(pe)^{1/p}}{\sqrt{2e\pi}}+r_{2}\log\frac{\Gamma\left(1+\frac{2}{p}\right)(p/2)^{2/p}}{e^{1-2/p}},\\ c:=& \frac{(r_{1}+r_{2})}{2}\log\left(\frac{p}{2}\right)\ \cdotp\end{split}\end{equation*}
\subsection{Expressions alternatives du degr{\'e} ad{\'e}lique}
\begin{lemma}\label{lemme4225}
Soit $\overline{E}$ un fibr{\'e} en \emph{droites} ad{\'e}lique et $s\in E\setminus\{0\}$. On a \begin{equation}\label{formuledegreun}\widehat{\deg}\overline{E}=-\sum_{v}{n_{v}\log\Vert s\Vert_{v}}\ .\end{equation}
\end{lemma}
\begin{proof}
Il suffit d'observer que si $x=\lambda.s\in E\otimes k_{v}$, $\lambda\in k_{v}$, alors $\Vert x\Vert_{v}\le 1$ {\'e}quivaut {\`a} $\vert\lambda\vert_{v}\le\frac{1}{\Vert s\Vert_{v}}$. L'{\'e}galit{\'e}~\eqref{equation:rayonboule} fait le lien entre les volumes et permet de conclure. 
\end{proof}On en d{\'e}duit le
\begin{lemma}\label{lemme4226}
Pour tout fibr{\'e} ad{\'e}lique \emph{hermitien} $\overline{E}$, on a $\widehat{\deg}\overline{E}=\widehat{\deg}\det\overline{E}$.
\end{lemma}
\begin{proof}
Nous avons vu que les normes de $\overline{E}$ aux places ultram{\'e}triques peuvent s'exprimer au moyen d'un ad{\`e}le fini $(c_{v})_{v\nmid\infty}\in\GL_{n}(k_{\mathbf{A},f})$, apr{\`e}s le choix d'une $k$-base $\mathbf{e}:=(e_{1},\ldots,e_{n})$ de $E$. L'hypoth{\`e}se faite ici sur $\overline{E}$ signifie que le m{\^e}me ph{\'e}nom{\`e}ne se produit aux places archim{\'e}diennes et qu'il existe un prolongement $c=(c_{v})\in\GL_{n}(k_{\mathbf{A}})$ tel que, pour toute place $v$ de $k$, on a \begin{equation*}(E_{v}\simeq_{\mathbf{e}}k_{v}^{n},\Vert.\Vert_{v})=c_{v}.(k_{v}^{n},\vert.\vert_{2})\ .\end{equation*}Ainsi chacun des quotients $$\frac{\vol\left(\{x\in k_{v}^{n}\,;\ \Vert x\Vert_{v}\le 1\}\right)}{\vol\left(\{x\in k_{v}^{n}\,;\ \vert x\vert_{2,v}\le 1\}\right)}$$ {\'e}gale $\vert\det(c_{v})\vert_{v}^{-n_{v}}=\Vert e_{1}\wedge\cdots\wedge e_{n}\Vert_{v}^{-n_{v}}$. Le lemme pr{\'e}c{\'e}dent permet de conclure.
\end{proof}
Ces deux lemmes montrent que le degr{\'e} ad{\'e}lique co{\"\i}ncide avec le degr{\'e} d'Arakelov sur le domaine de d{\'e}finition de ce dernier, {\`a} savoir les fibr{\'e}s vectoriels hermitiens sur $\spec\mathcal{O}_{k}$. Le degr{\'e} ad{\'e}lique en est donc une g{\'e}n{\'e}ralisation.\par  Lorsque $k$ poss{\`e}de des places archim{\'e}diennes, la norme $\Vert.\Vert_{v}$ de $\overline{E}$ en une telle place est comparable {\`a} la norme de John, introduite {\`a} la suite de la d{\'e}finition~\ref{defiellipsoidejohn}. Notons $\vert.\vert_{J(\overline{E}),v}$ cette norme hermitienne {\'e}tendue {\`a} l'espace $E\otimes\mathbf{C}_{v}$. 
\begin{defi}
Avec les donn{\'e}es ci-dessus, nous appelerons \emph{fibr{\'e} vectoriel ad{\'e}lique de John}, ou plus simplement \emph{fibr{\'e} de John}, associ{\'e} {\`a} $\overline{E}$, et nous noterons  $J(\overline{E})$, le fibr{\'e} vectoriel ad{\'e}lique obtenu {\`a} partir de $\overline{E}$ en rempla{\c c}ant, aux places archim{\'e}diennes $v$ de $k$, la norme $\Vert.\Vert_{v}$ par la norme hermitienne $\vert.\vert_{J(\overline{E}),v}$.
\end{defi}Le fibr{\'e} de John est donc {\'e}gal {\`a} $\overline{E}$ lorsque ce dernier est un fibr{\'e} ad{\'e}lique hermitien. De la m{\^e}me mani{\`e}re, on d{\'e}finit le \emph{fibr{\'e} vectoriel ad{\'e}lique de L{\"o}wner} associ{\'e} {\`a} $\overline{E}$ --- que l'on note $L(\overline{E})$ --- au moyen des m{\'e}triques $\vert.\vert_{L(\mathbf{B}(E_{v},\Vert.\Vert_{v}))}$ introduites au paragraphe~\ref{paragraphedeuxdeux}, en une place archim{\'e}dienne $v$ de $k$. On notera que si $E'$ est un sous-espace vectoriel de $E$ alors $(E',\vert.\vert_{L(\overline{E})})\preceq(E',\Vert.\Vert_{E})\preceq(E',\vert.\vert_{J(\overline{E})})$ et donc, par le lemme~\ref{lemmesimplee}, on a \begin{equation}\label{comparaisonsousfibres}\widehat{\deg}(E',\vert.\vert_{J(\overline{E})})\le\widehat{\deg}(E',\Vert.\Vert_{E})\le\widehat{\deg}(E',\vert.\vert_{L(\overline{E})})\ .\end{equation}On peut m{\^e}me intercaler $\widehat{\deg}J(\overline{E'})$ entre $\widehat{\deg}(E',\vert.\vert_{J(\overline{E})})$ et $\widehat{\deg}(E',\Vert.\Vert_{E})$ en utilisant la d\'efinition des m\'etriques de $J(\overline{E'})$ qui sont les m\'etriques euclidiennes qui optimisent par valeur inf\'erieure le volume de la boule unit\'e de $\overline{E'}$. La m{\^e}me observation avec $L(\overline{E'})$ conduit \`a \begin{equation*}\widehat{\deg}(E',\vert.\vert_{J(\overline{E})})\le\widehat{\deg}J(\overline{E'})\le\widehat{\deg}(E',\Vert.\Vert_{E})\le\widehat{\deg}L(\overline{E'})\le\widehat{\deg}(E',\vert.\vert_{L(\overline{E})})\ .\end{equation*}Le fibr{\'e} de John permet de donner une formule exacte pour le degr{\'e} ad{\'e}lique de $\overline{E}$ gr{\^a}ce {\`a} la notion de quotient volumique (voir d{\'e}finition~\ref{defiquotientvolumique}).
\begin{defi}
Soit $\overline{E}$ un fibr{\'e} vectoriel ad{\'e}lique sur $\spec k$. Le \emph{quotient volumique ad{\'e}lique} de $\overline{E}$, not{\'e} $\vr(\overline{E})$, est la norme de l'ad{\`e}le de composantes $\vr(E\otimes_{k}k_{v})$ --- le quotient volumique de l'espace de Banach r{\'e}el ou complexe $(E\otimes_{k}k_{v},\Vert.\Vert_{E,v})$ --- aux places $v$ archim\'ediennes et $1$ aux autres places~:
\begin{equation*}\vr(\overline{E}):=\prod_{\text{$v$ archim\'edienne}}{\vr(E\otimes_{k}k_{v})^{n_{v}}}\ .
\end{equation*}
\end{defi}Si $k$ est un corps de fonctions, on a donc toujours $\vr(\overline{E})=1$. Ce nombre est un terme d'erreur qui intervient souvent lors de comparaisons entre la situation ad{\'e}lique g{\'e}n{\'e}rale et le cas hermitien. On a $\vr(\overline{E})=1$ si et seulement si $\overline{E}$ est un fibr{\'e} ad{\'e}lique \emph{hermitien}. Aussi, les in{\'e}galit{\'e}s {\'e}tablies dans la suite et qui font intervenir cette quantit\'e sont des {\'e}galit{\'e}s dans le cas hermitien. De mani{\`e}re ponctuelle, nous aurons besoin {\'e}galement du nombre r{\'e}el \begin{equation*}\widetilde{\vr}(\overline{E}):=\prod_{\text{$v$ archim\'edienne}}{\widetilde{\vr}(E\otimes_{k}k_{v})^{n_{v}}}
\end{equation*}(voir d{\'e}finition~\ref{quotientvolumiquebis} pour la notation $\widetilde{\vr}$). Introduisons maintenant l'analogue ad\'elique de la distance de Banach-Mazur classique (d\'efinition~\ref{distancebanachmazur}).
\begin{defi}
Soit $\overline{E}$ et $\overline{F}$ deux fibr{\'e}s vectoriels ad{\'e}liques sur $\spec k$ de m{\^e}me dimension. La \emph{distance de Banach-Mazur ad{\'e}lique} entre $\overline{E}$ et $\overline{F}$ est le nombre r\'eel\begin{equation*}\mathrm{d}(\overline{E},\overline{F}):=\prod_{\text{$v$ archim\'edienne}}{\mathrm{d}(E\otimes_{k}k_{v},F\otimes_{k}k_{v})^{n_{v}}}\ .\end{equation*}Si $k$ est un corps de fonctions, on pose $\mathrm{d}(\overline{E},\overline{F})=1$. 
\end{defi}\noindent
\textbf{Notation}. Lorsque $\overline{F}=(k^{n},\vert.\vert_{2})$, on notera plus simplement $\Delta(\overline{E})$ la distance $\mathrm{d}(\overline{E},(k^{n},\vert.\vert_{2}))$.\par On montre que si $\overline{E'}$ est un sous-fibr{\'e} vectoriel ad{\'e}lique de $\overline{E}$ alors $\Delta(\overline{E'})\le\Delta(\overline{E})$. De plus, le dual $\overline{E^{\mathsf{v}}}$ de $\overline{E}$ v{\'e}rifie $\Delta(\overline{E^{\mathsf{v}}})=\Delta(\overline{E})$. Gr{\^a}ce aux estimations locales~\eqref{estimationslocalesreelles} et~\eqref{estimationslocalescomplexes}, on dispose de l'encadrement \begin{equation}\label{eq:encadrementquotientadelique}1\le\vr(\overline{E})\vr(\overline{E^{\mathsf{v}}})\le\Delta(\overline{E})\le(2n)^{D/2}\ .\end{equation}  
\begin{prop}\label{prop:comparaisonjohnadelique}
Soit $\overline{E}$ un fibr{\'e} vectoriel ad{\'e}lique sur $\spec k$. On a \begin{equation*}\widehat{\deg}\overline{E}=\widehat{\deg}J(\overline{E})+n\log\vr(\overline{E})\quad\text{et}\quad\widehat{\deg} L(\overline{E})=\widehat{\deg}\overline{E}+n\log\widetilde{\vr}(\overline{E})\ .\end{equation*}
\end{prop}
\begin{proof}
On identifie $E$ {\`a} $k^{n}$ en choisissant une $k$-base quelconque de $E$. La contribution au degr{\'e} ad{\'e}lique de $\overline{E}$ des places ultram{\'e}triques est exactement la m{\^e}me que celle de $J(\overline{E})$, par d{\'e}finition de ce dernier fibr{\'e}. En une place archim{\'e}dienne $v$ de $k$, le quotient volumique de $(k_{v}^{n},\Vert.\Vert_{v})$ intervient \`a travers l'{\'e}galit{\'e}~:\begin{equation*}\begin{split}&\frac{\vol\left(\{x\in k_{v}^{n}\,;\ \Vert x\Vert_{v}\le 1\}\right)}{\vol\left(\{x\in k_{v}^{n}\,;\ \vert x\vert_{2}\le 1\}\right)}\\ &\qquad=\frac{\vol\left(\{x\in k_{v}^{n}\,;\ \vert x\vert_{J(\overline{E}),v}\le 1\}\right)}{\vol\left(\{x\in k_{v}^{n}\,;\ \vert x\vert_{2}\le 1\}\right)}\cdotp\vr(k_{v}^{n},\Vert.\Vert_{v})^{\dim_{\mathbf{R}}k_{v}^{n}}\ .\end{split}\end{equation*}La proposition s'en d{\'e}duit alors car $\dim_{\mathbf{R}}k_{v}^{n}=n_{v}n$. La m{\^e}me d{\'e}marche avec l'ellipso{\"\i}de de L{\"o}wner permet de montrer l'autre {\'e}galit{\'e}.\end{proof} Par ailleurs, pour tout $\varepsilon>0$ et gr{\^a}ce {\`a} l'encadrement~\eqref{encadrementhermitien}, on a l'existence d'une norme hermitienne $\vert.\vert_{\varepsilon,v}$ sur $E\otimes_{k}k_{v}$ ($v$ archim\'edienne) telle que, pour tout sous-espace vectoriel $E'$ de $E$, \begin{equation}\label{eqref:important}(E',(\vert.\vert_{\varepsilon,v})_{v})\preceq\overline{E'}\preceq\left(E',(\mathrm{d}(E\otimes_{k}k_{v},\ell^{2}_{n,k_{v}})(1+\varepsilon)\vert.\vert_{\varepsilon,v})_{v}\right)\end{equation}(les normes aux places ultram{\'e}triques des fibr{\'e}s ad{\'e}liques de gauche et de droite sont celles de $\overline{E}$).\par\noindent\textbf{Notation}. On d{\'e}signe par $\overline{E}_{\varepsilon}$ le fibr{\'e} vectoriel ad{\'e}lique d'espace vectoriel sous-jacent $E$ et dont les normes sont\begin{equation*}\Vert.\Vert_{\overline{E}_{\varepsilon},v}:=\begin{cases}\Vert.\Vert_{E,v} & \text{si $v$ est ultram\'etrique},\\ \vert.\vert_{\varepsilon,v} & \text{si $v$ est archim\'edienne}.
\end{cases}\end{equation*}On notera aussi plus simplement $\vert.\vert_{\varepsilon}$ la collection $(\Vert.\Vert_{\overline{E}_{\varepsilon},v})_{v}$ des normes de $\overline{E}_{\varepsilon}$. Il est bien clair que cette d\'efinition repose sur le choix des m\'etriques hermitiennes $\vert.\vert_{\varepsilon,v}$ et que le fibr\'e $\overline{E}_{\varepsilon}$ n'est pas d\'etermin\'e uniquement par $\overline{E}$ et $\varepsilon$.\par Les d{\'e}finitions m{\^e}me des fibr{\'e}s de John et L{\"o}wner entra{\^\i}nent $\widehat{\deg}_{\mathrm{n}}\overline{E}_{\varepsilon}\le\widehat{\deg}_{\mathrm{n}}J(\overline{E})$ et \begin{equation*}\widehat{\deg}_{\mathrm{n}}L(\overline{E})\le\widehat{\deg}_{\mathrm{n}}\left(E,(\mathrm{d}(E\otimes_{k}k_{v},\ell^{2}_{n,k_{v}})(1+\varepsilon)\vert.\vert_{\varepsilon,v})_{v}\right),\end{equation*}puis, en faisant tendre $\varepsilon$ vers $0$, l'on obtient \begin{equation}\label{eqref:estimationdifference}0\le\widehat{\deg}_{\mathrm{n}}L(\overline{E})-\widehat{\deg}_{\mathrm{n}}J(\overline{E})\le\frac{n}{D}\log\Delta(\overline{E})\ \cdotp\end{equation}De l'encadrement~\eqref{eqref:important}, on d{\'e}duit aussi \begin{equation*}-\frac{n}{D}\log\Delta(\overline{E})-n\log(1+\varepsilon)\le\widehat{\deg}_{\mathrm{n}}\overline{E'}-\widehat{\deg}_{\mathrm{n}}\left(E',(\vert.\vert_{\varepsilon,v})_{v}\right)\le 0\ \cdotp\end{equation*}Ceci permet d'approcher le degr{\'e} ad{\'e}lique de $\overline{E'}$ par le degr{\'e} d'un fibr{\'e} ad{\'e}lique \emph{hermitien}, avec un terme d'erreur de l'ordre de $-\frac{n}{D}\log\Delta(\overline{E})$. Cette observation sera constamment utilis{\'e}e dans la suite.
\begin{coro}\label{lemmedecomparaison47}Pour tout fibr{\'e} vectoriel ad{\'e}lique $\overline{E}$ et toute norme tensorielle ad\'elique hermitienne $\alpha$ d'ordre $\dim E$, on a\begin{equation*}\left\vert\widehat{\deg}_{\mathrm{n}}\overline{E}-\widehat{\deg}_{\mathrm{n}}\mathrm{det}_{\alpha}\overline{E}\right\vert\le\widehat{\deg}_{\mathrm{n}}L(\overline{E})-\widehat{\deg}_{\mathrm{n}}J(\overline{E})\ .\end{equation*}
\end{coro}
\begin{proof}
La diff{\'e}rence ne se situe {\'e}ventuellement qu'en une place archim{\'e}dienne $v$ de $k$. On a vu que dans ce cas  $\vert.\vert_{L(\overline{E}),v}\le\Vert.\Vert_{E,v}\le\vert.\vert_{J(\overline{E}),v}$. L'hypoth{\`e}se $\alpha$ hermitienne et les lemmes~\ref{lemmesimplee} et~\ref{lemmecomptrois} entra{\^\i}nent alors l'in{\'e}galit{\'e} \begin{equation*}\widehat{\deg}_{\mathrm{n}}\det J(\overline{E})\le\widehat{\deg}_{\mathrm{n}}\mathrm{det}_{\alpha}\overline{E}\le\widehat{\deg}_{\mathrm{n}}\det L(\overline{E})\ .\end{equation*}Le lemme~\ref{lemme4226} et la proposition~\ref{prop:comparaisonjohnadelique} donnent \begin{equation*}-\frac{n}{D}\log\widetilde{\vr}(\overline{E})\le\widehat{\deg}_{\mathrm{n}}\overline{E}-\widehat{\deg}_{\mathrm{n}}\mathrm{det}_{\alpha}\overline{E}\le\frac{n}{D}\log\vr(\overline{E}),\end{equation*}qui est un encadrement un peu plus pr\'ecis que celui annonc\'e puisque $\max{\{\vr(\overline{E}),\widetilde{\vr}(\overline{E})\}}\le\vr(\overline{E})\widetilde{\vr}(\overline{E})$.
\end{proof}On rappelle que $\overline{k}$ d\'esigne une cl{\^o}ture alg\'ebrique de $k$.
\begin{defi}\label{definition48hauteur}
Soit $\overline{E}=(E,(\Vert.\Vert_{v})_{v})$ un fibr{\'e} vectoriel ad{\'e}lique. La \emph{hauteur} d'un {\'e}l{\'e}ment $x\in (E\otimes_{k}\overline{k})\setminus\{0\}$ relative {\`a} $\overline{E}$, not{\'e}e $h_{\overline{E}}(x)$, est la somme normalis{\'e}e \begin{equation}h_{\overline{E}}(x):=\frac{1}{D}\sum_{v}{n_{v}\log\Vert x\Vert_{E,v}}\end{equation}(la somme porte sur toutes les places $v$ de $k$ mais seuls un nombre fini de termes ne sont pas nuls).
\end{defi}
C'est aussi l'oppos{\'e} du degr{\'e} normalis{\'e} du fibr{\'e} en droites $(K.x,(\Vert.\Vert_{E,v})_{v})$ o\`u $K$ est une extension finie de $k$ sur laquelle est d\'efini $x$. La hauteur de $x$ ne d\'epend pas du choix de $K$ en raison de la formule $\sum_{w\mid v}{[K_{w}:k_{v}]}=[K:k]$ (voir aussi le paragraphe suivant). La terminologie est coh{\'e}rente avec la notion de hauteur au sens usuel gr{\^a}ce au r{\'e}sultat suivant.
\begin{prop}\label{propositionsurlahauteur}
Soit $a\in\mathbf{R}$. L'ensemble des droites $\{k.x\,;\ x\in E\ \text{et}\ h_{\overline{E}}(x)\le a\}$ est fini.
\end{prop}
\begin{proof}
Soit $(e_{1},\ldots,e_{n})$ une $k$-base de $E$. Il existe une constante $\alpha>0$ (qui d{\'e}pend de $\overline{E}$ et de $k$) telle que, pour toute place $v$ de $k$ et pour tout $x=\sum_{i=1}^{n}{x_{i}e_{i}}\in E$, on ait $\Vert x\Vert_{\overline{E},v}\ge\alpha\max_{1\le i\le n}{\{\vert x_{i}\vert_{v}\}}$. Aux places ultram{\'e}triques (et il suffit de consid{\'e}rer un nombre \emph{fini} d'entre elles car $\overline{E}$ est ad{\'e}lique), cela r{\'e}sulte de la formule~\eqref{normesegales}. Et aux places archim{\'e}diennes, il s'agit d'une cons{\'e}quence de l'{\'e}quivalence des normes en dimension finie. Si $x\ne 0$, la borne $h_{\overline{E}}(x)\le a$ entra{\^\i}ne une majoration de la hauteur de Weil du point projectif $(x_{1}:\cdots:x_{n})$ et donc un nombre fini de tels points (th\'eor\`eme de Northcott) et de vecteurs $x$ correspondant.  
\end{proof} 
\begin{prop}Soit $\overline{E}$ un fibr\'e vectoriel ad\'elique sur $\spec k$ de dimension $n$ et $\alpha$ une norme tensorielle ad\'elique d'ordre $n$. Soit $(e_{1},\ldots,e_{n})$ une $k$-base de $E$. Alors on a\begin{equation*}0\le \widehat{\deg}_{\mathrm{n}}\mathrm{det}_{\alpha}\overline{E}+h_{\overline{E}}(e_{1})+\cdots+h_{\overline{E}}(e_{n})\ .\end{equation*}De plus, on a\begin{equation*}-\frac{n}{D}\log\Delta(\overline{E})\le \widehat{\deg}_{\mathrm{n}}\overline{E}+h_{\overline{E}}(e_{1})+\cdots+h_{\overline{E}}(e_{n})\ .\end{equation*}
\end{prop}
\begin{proof}
La premi{\`e}re assertion est une cons\'equence de l'in\'egalit\'e d'Hadamard \begin{equation*}\Vert e_{1}\wedge\cdots\wedge e_{n}\Vert_{\mathrm{det}_{\alpha}\overline{E},v}\le\Vert e_{1}\otimes\cdots\otimes e_{n}\Vert_{\overline{E}^{\otimes_{\alpha}n},v}\le\prod_{i=1}^{n}{\Vert e_{i}\Vert_{E,v}},\end{equation*}vraie pour toute place $v$ de $k$. En effet, la premi{\`e}re majoration vient de la d\'efinition de norme quotient et la seconde est la propri\'et\'e (i) de la d\'efinition~\ref{defidenormestensoriellesreelles} des normes tensorielles (il y a m{\^e}me \'egalit\'e). Pour obtenir la seconde assertion, on choisit pour $\alpha$ une norme hermitienne (construite, par exemple, au moyen de la norme hermitienne $g_{2}$ de Chevet-Saphard) dans le corollaire~\ref{lemmedecomparaison47} et on utilise l'estimation~\eqref{eqref:estimationdifference}. 
\end{proof} Il est souvent utile de conna{\^\i}tre {\'e}galement une majoration de la somme des hauteurs d'une base de $E\otimes\overline{k}$ o\`u $\overline{k}$ d\'esigne une cl{\^o}ture alg\'ebrique de $k$ . C'est l'objet du \emph{lemme de Siegel absolu} suivant, {\'e}tabli dans le cas des corps de fonctions par D.~Roy \& J.L.~Thunder~\cite{roy-thunder} et issu d'un {\'e}nonc{\'e} d{\^u} {\`a} S.~Zhang~\cite{Zhang2} dans le cas d'un corps de nombres (les deux \'etant \'ecrits dans le cas hermitien).
\begin{theo}\label{lemmedesiegelabsolu}
Soit $\overline{E}$ un fibr{\'e} vectoriel ad{\'e}lique, de dimension $n\ge 1$. Posons $\delta:=0$ si $k$ est un corps de fonctions et $\delta:=1$ si $k$ est un corps de nombres. Pour tout $\varepsilon>0$, il existe une base $(e_{1},\ldots,e_{n})$ de $E\otimes\overline{k}$ telle que \begin{equation}\label{inegalite23}h_{\overline{E}}(e_{1})+\cdots+h_{\overline{E}}(e_{n})+\widehat{\deg}_{\mathrm{n}}\overline{E}\le\frac{\delta n}{2}\log n+\frac{n}{D}\log\vr(\overline{E})+\varepsilon\ .\end{equation} 
\end{theo}
Dans le cas d'un corps de nombres, la preuve compl{\`e}te se trouve dans~\cite{casrationnel}. Il est possible de demander que la base $(e_{1},\ldots,e_{n})$ soit d{\'e}finie sur $k$ et, en contrepartie, de rajouter un terme d{\'e}pendant du discriminant absolu $D_{k}$ de $k$ dans le membre de droite de l'in{\'e}galit{\'e}~\eqref{inegalite23} (voir~\cite{roy-thunder}). Ce terme est $(n/D)\log\vert D_{k}\vert$ si $k$ est un corps de nombres. Le cas g\'en\'eral d'un fibr\'e vectoriel ad\'elique non n\'ecessairement hermitien s'obtient en appliquant le th\'eor\`eme ci-dessus au fibr\'e (hermitien) de John $J(\overline{E})$ et en utilisant $\overline{E}\preceq J(\overline{E})$ ainsi que la proposition~\ref{prop:comparaisonjohnadelique}.

\subsection{Extension des scalaires}En toute g{\'e}n{\'e}ralit{\'e}, le degr{\'e} ad{\'e}lique normalis{\'e} n'est pas un invariant stable par extension des scalaires. Si l'on consid{\`e}re par exemple le fibr{\'e} vectoriel ad{\'e}lique $(\mathbf{Q}^{n},\vert.\vert_{\infty})$, la formule~\eqref{formuledevolumep} entra{\^\i}ne l'{\'e}galit{\'e} $$\widehat{\deg}_{\mathrm{n}}(\mathbf{Q}^{n},\vert.\vert_{\infty})=\log\frac{\Gamma\left(1+\frac{n}{2}\right)}{\pi^{n/2}}\ \cdotp$$Or si l'on choisit l'extension quadratique $\mathbf{Q}(i)$, qui poss{\`e}de une seule place complexe et aucune r{\'e}elle, un rapide calcul donne $$\widehat{\deg}_{\mathrm{n}}(\mathbf{Q}(i)^{n},\vert.\vert_{\infty})=\frac{1}{2}\log n!,$$nombre toujours diff{\'e}rent de celui de la formule pr{\'e}c{\'e}dente, d\`es que $n\ge 2$. Cette anomalie se voit {\'e}galement au moyen de la formule asymptotique~\eqref{formuleasymptotiquedegre} dont les coefficients $b$ et $c$ d{\'e}pendent en g{\'e}n{\'e}ral de $r_1$ \emph{et} $r_2$ et non seulement de la somme $r_1+2r_{2}=[k:\mathbf{Q}]$.\par  
L'explication de ce ph{\'e}nom{\`e}ne est l'observation suivante, qui, r{\'e}duite {\`a} l'essentiel, marque la diff{\'e}rence entre le volume d'un cube et celui d'une boule.
\begin{lemma}\label{lemme:normescomplexereelle}
Soit $\Vert.\Vert$ une norme sur $\mathbf{C}^{n}$, invariante par conjugaison complexe. Il existe une constante $\kappa=\kappa(n,\Vert.\Vert)$, comprise entre $4^{-n}$ et $4^{n}$, telle que \begin{equation}\label{egaliteprovisoire}\frac{\vol\left(\left\{x\in\mathbf{C}^{n}\,;\ \Vert x\Vert\le 1\right\}\right)}{\vol\left(\left\{x\in\mathbf{C}^{n}\,;\ \vert x\vert_{2}\le 1\right\}\right)}=\kappa\,\left(\frac{\vol\left(\left\{x\in\mathbf{R}^{n}\,;\ \Vert x\Vert\le 1\right\}\right)}{\vol\left(\left\{x\in\mathbf{R}^{n}\,;\ \vert x\vert_{2}\le 1\right\}\right)}\right)^{2}\ \cdotp
\end{equation}La notation $\vol$ d{\'e}signe indiff{\'e}remment une mesure de Haar (quelconque) sur $\mathbf{R}^{n}$ ou $\mathbf{C}^{n}$. Si la norme $\Vert.\Vert$ est hermitienne alors $\kappa=1$.
\end{lemma}
\begin{proof}
Il s'agit d'une simple cons{\'e}quence des in{\'e}galit{\'e}s~:
\begin{equation*}
\max{\left\{\Vert x_{1}\Vert,\Vert x_{2}\Vert\right\}}\le\frac{1}{2}\left(\Vert x_{1}+i x_{2}\Vert+\Vert x_{1}-ix_{2}\Vert\right)=\Vert x_{1}+ix_{2}\Vert
\end{equation*}et $\Vert x_{1}+ix_{2}\Vert\le\Vert x_{1}\Vert+\Vert x_{2}\Vert\le 2\max{\left\{\Vert x_{1}\Vert,\Vert x_{2}\Vert\right\}}$, pour tous $x_{1},x_{2}\in\mathbf{R}^{n}$. Le r{\'e}sultat s'ensuit en choisissant sur $\mathbf{C}^{n}$ la mesure de Haar produit $\vol_{\mathbf{R}^{n}}\otimes\vol_{\mathbf{R}^{n}}$, obtenue au moyen de l'isomorphisme naturel $\mathbf{C}^{n}\simeq\mathbf{R}^{n}\times\mathbf{R}^{n}$. Si la norme est hermitienne, elle est d{\'e}termin{\'e}e par la donn\'ee d'une matrice sym\'etrique $A\in\mathrm{M}_{n}(\mathbf{R})$ (il n'y a pas de partie imaginaire car la norme induite doit {\^e}tre invariante par conjugaison complexe). Les deux membres de l'\'equation~\eqref{egaliteprovisoire} valent alors $(\det A)^{-1}$ et $\kappa=1$.
\end{proof}
\begin{rema}\label{remarque43}
Les calculs exacts des volumes de $b_{n,\mathbf{C}}^{2}$ et $b_{n,\mathbf{R}}^{2}$ donn{\'e}es par les formules~\eqref{formuledevolumep} et~\eqref{formulevolumecomplexep} fournissent l'encadrement plus pr{\'e}cis $$\frac{1}{4^{n}}\cdotp\frac{n!}{\Gamma\left(1+\frac{n}{2}\right)^{2}}\le\kappa\le\frac{n!}{\Gamma\left(1+\frac{n}{2}\right)^{2}}\ \cdotp$$On peut noter que la formule de Stirling (voir note au bas de la page~\pageref{Stirling}) donne l'{\'e}quivalent asymptotique $$\frac{n!}{\Gamma\left(1+\frac{n}{2}\right)^{2}}\sim\frac{2^{n+1/2}}{\sqrt{\pi n}}\quad\text{lorsque $n\longrightarrow+\infty$.}$$
\end{rema}
Apr{\`e}s ces consid{\'e}rations, nous pouvons {\'e}noncer le r{\'e}sultat principal de ce paragraphe.
\begin{prop}\label{proposition43}
Soit $K$ une extension finie de $k$ et $\overline{E}$ un fibr{\'e} vectoriel ad{\'e}lique sur $k$, de dimension $n$. Alors \begin{equation}\label{majorationdelaproposition43}
\left\vert\widehat{\deg}_{\mathrm{n}}\overline{E}_{K}-\widehat{\deg}_{\mathrm{n}}\overline{E}\right\vert\le n\log(4)\ .
\end{equation}Si, de plus, $\overline{E}$ est un fibr{\'e} ad{\'e}lique hermitien alors $\widehat{\deg}_{\mathrm{n}}\overline{E}_{K}=\widehat{\deg}_{\mathrm{n}}\overline{E}$.    
\end{prop}
\begin{proof}
Quitte {\`a} fixer une $k$-base de $E$, on peut supposer $E=k^{n}$. On distingue deux cas.\\ \noindent \texttt{a}) \emph{Soit $v$ une place ultram{\'e}trique de $k$}. On a vu qu'il existe une matrice $c_{v}\in\GL_{n}(k_{v})$ telle que :\begin{equation*}\forall x=(x_{1},\ldots,x_{n})\in\mathbf{C}_{v}^{n},\quad\Vert x\Vert_{E,v}=\underset{1\le i\le n}{\mathrm{max}}\left\{\left\vert (c_{v}.x)_{i}\right\vert_{v}\right\}\ \cdotp\end{equation*}Notons $B_{v}:=\{x\in k_{v}^{n}\,;\ \Vert x\Vert_{v}\le 1\}$. Cet ensemble est l'image par $c_{v}^{-1}$ de $\mathcal{O}_{v}^{n}$ et, d'apr{\`e}s les rappels du paragraphe~\ref{section:preliminaires}, on a $$\vol_{v}(B_{v})=\vert\det c_{v}\vert_{v}^{-n_{v}}\vol_{v}(\mathcal{O}_{v}^{n})$$(ici $\vol_{v}$ est une mesure de Haar sur $k_{v}^{n}$). Soit $w$ une place de $K$ au-dessus de $v$. Le compl{\'e}t{\'e} $K_{w}$ est un $k_{v}$-espace vectoriel de dimension $[K_{w}:k_{v}]$, et si l'on note $B_{w}$ la boule unit{\'e} ferm{\'e}e dans $K_{w}^{n}$, on a $$\frac{\vol_{w}(B_{w})}{\vol_{w}(\mathcal{O}_{w}^{n})}=\vert\det c_{v}\vert_{v}^{-n_{w}}=\left(\frac{\vol_{v}(B_{v})}{\vol_{v}(\mathcal{O}_{v}^{n})}\right)^{[K_{w}:k_{v}]}\ .$$Comme cela est rappel{\'e} dans les pr\'eliminaires, l'on dispose de l'isomorphisme d'espaces vectoriels topologiques \begin{equation}\label{isomseparable}K\otimes_{k}k_{v}\simeq\prod_{w\mid v}{K_{w},}\end{equation}qui, outre l'{\'e}galit{\'e} des dimensions $[K:k]=\sum_{w\mid v}{[K_{w}:k_{v}]}$, entra{\^\i}ne \begin{equation*}\begin{split}&\frac{\vol\left(\{x\in(K\otimes k_{v})^{n}\,;\ \Vert x\Vert_{v}\le 1\}\right)}{\vol\left(\left\{x\in(K\otimes k_{v})^{n}\,;\ \underset{1\le i\le n}{\mathrm{max}}{\{\vert x_{i}\vert_{v}\}\le 1\}}\right\}\right)}\\ &\quad=\prod_{w\mid v}{\frac{\vol_{w}(B_{w})}{\vol_{w}(\mathcal{O}_{w}^{n})}}=\left(\frac{\vol_{v}(B_{v})}{\vol_{v}(\mathcal{O}_{v}^{n})}\right)^{[K:k]}\end{split}\end{equation*}(la mesure $\vol$ sur $(K\otimes k_{v})^{n}$ peut {\^e}tre choisie comme l'image r{\'e}ciproque de la mesure produit $\otimes_{w\mid v}\vol_{w}$ sur $\prod_{w\mid v}{K_{w}^{n}}$ par l'isomorphisme~\eqref{isomseparable}). Soulignons que cette {\'e}galit{\'e} ne d{\'e}pend ni du choix de l'isomorphisme~\eqref{isomseparable}, ni du choix des mesures $\vol_{w}$, $\vol_{v}$ ou bien encore $\vol$. Ainsi le comportement du quotient des volumes aux places ultram{\'e}triques par extension des scalaires est exactement celui souhait{\'e} et les parties (aux places) finies des degr{\'e}s ad{\'e}liques normalis{\'e}s de $\overline{E}_{K}$ et $\overline{E}$ sont identiques. En particulier, ces degr{\'e}s (dans leur totalit{\'e}) sont {\'e}gaux si $k$ est un corps de fonctions. \\ \noindent \texttt{b}) \emph{Soit $v$ une place archim{\'e}dienne de $k$}. Comme pr{\'e}c{\'e}demment, soit $w$ une place de $K$ au-dessus de $v$. Si $w$ et $v$ sont de m{\^e}me nature (toutes les deux r{\'e}elles ou complexes), le quotient des volumes $$\frac{\vol(\{x\in k_{v}^{n}\,;\ \Vert x\Vert_{v}\le 1\})}{\vol(\{x\in k_{v}^{n}\,;\ \vert x\vert_{2}\le 1\})}$$ reste inchang{\'e}. En revanche, si $v$ est r{\'e}elle et $w$ est complexe, nous sommes dans le cas du figure du lemme~\ref{lemme:normescomplexereelle} et la proposition~\ref{proposition43} s'en d{\'e}duit imm{\'e}diatement (le cas d'{\'e}galit{\'e} lorsque $\overline{E}$ est un fibr{\'e} vectoriel hermitien provenant du cas d'{\'e}galit{\'e} $\kappa=1$ dans ce m{\^e}me lemme).
\end{proof}
\begin{rema} L'argument qui est {\`a} la fin de cette d{\'e}monstration permet d'{\'e}tablir en r{\'e}alit{\'e} une borne plus faible pour la diff{\'e}rence des degr{\'e}s puisqu'il suffit de consid{\'e}rer les couples de places $(v,w)$ telles que $w\mid v$ et qui ne sont pas de m{\^e}me nature. Si $N(k,K)$ est le nombre de tels couples, la majoration~\eqref{majorationdelaproposition43} reste vraie en multipliant $n\log(4)$ par le quotient $N(k,K)/[K:\mathbf{Q}]\le 1$.  
\end{rema}
\subsection{Somme directe}Comme nous l'avons vu au paragraphe~\ref{paragraphetroistrois}, il est possible de d{\'e}finir la somme directe $\overline{E}\oplus_{\varsigma}\overline{F}$ de deux fibr{\'e}s vectoriels ad{\'e}liques, relative {\`a} une famille $\varsigma=(\varsigma_{v})_{v\mid\infty}$ de normes (particuli{\`e}res) sur $\mathbf{R}^{2}$. La proposition suivante {\'e}value le degr{\'e} d'une telle somme directe.
\begin{prop}\label{degresommedirecte}
Soit $\overline{E}$ et $\overline{F}$ des fibr{\'e}s vectoriels ad{\'e}liques sur $\spec k$, de dimensions respectives $n$ et $m$. On a alors \begin{equation*}\left\vert\widehat{\deg}_{\mathrm{n}}(\overline{E}\oplus_{\varsigma}\overline{F})-\widehat{\deg}_{\mathrm{n}}\overline{E}-\widehat{\deg}_{\mathrm{n}}\overline{F}\right\vert\le\log\binom{n+m}{n}\ \cdotp\end{equation*}Pour la somme hermitienne (c.-{\`a}-d. $\varsigma_{v}(x,y)=(x^{2}+y^{2})^{1/2}$ aux places archim\'ediennes $v$ de $k$), on a l'{\'e}galit{\'e} $\widehat{\deg}_{\mathrm{n}}(\overline{E}\oplus_{2}\overline{F})=\widehat{\deg}_{\mathrm{n}}\overline{E}+\widehat{\deg}_{\mathrm{n}}\overline{F}$.
\end{prop}
En particulier, si $k$ est un corps de fonctions, on a $\widehat{\deg}_{\mathrm{n}}(\overline{E}\oplus\overline{F})=\widehat{\deg}_{\mathrm{n}}\overline{E}+\widehat{\deg}_{\mathrm{n}}\overline{F}$.
\begin{proof}
En vertu de l'encadrement~\eqref{sommedirecteplus210}, on a $\overline{E}\oplus_{\infty}\overline{F}\preceq\overline{E}\oplus_{\varsigma}\overline{F}\preceq\overline{E}\oplus_{1}\overline{F}$ et donc $$\widehat{\deg}(\overline{E}\oplus_{1}\overline{F})\le\widehat{\deg}(\overline{E}\oplus_{\varsigma}\overline{F})\le\widehat{\deg}(\overline{E}\oplus_{\infty}\overline{F})$$(lemme~\ref{lemmesimplee}). Or pour $p\in[1,+\infty]$ on dispose d'une formule exacte pour $\widehat{\deg}(\overline{E}\oplus_{p}\overline{F})$. En effet, si l'on applique l'{\'e}galit{\'e}~\eqref{sommedirectemoins211}, {\`a} $E_{v}\oplus_{p}F_{v}$ o\`u $v$ est une place archim{\'e}dienne de $k$, la diff{\'e}rence $\widehat{\deg}(\overline{E}\oplus_{p}\overline{F})-\widehat{\deg}\overline{E}-\widehat{\deg}\overline{F}$ {\'e}gale \begin{equation*}\sum_{\text{$v$ archim\'edienne}}{\log\frac{\Gamma\left(1+\frac{nn_{v}}{p}\right)\Gamma\left(1+\frac{mn_{v}}{p}\right)\Gamma\left(1+\frac{(m+n)n_{v}}{2}\right)}{\Gamma\left(1+\frac{(m+n)n_{v}}{p}\right)\Gamma\left(1+\frac{nn_{v}}{2}\right)\Gamma\left(1+\frac{mn_{v}}{2}\right)}}\ \cdotp\end{equation*}Il n'y a pas de contributions aux places ultram{\'e}triques car, par choix de la norme sup, la boule unit{\'e} de $E_{v}\oplus_{p}F_{v}$ est le produit des boules unit{\'e}s de $E_{v}$ et $F_{v}$. Et les volumes se multiplient {\'e}galement pourvu que l'on choisisse la mesure produit $\vol_{E,v}\otimes\vol_{F,v}$ sur la somme directe $E_{v}\oplus F_{v}$. On obtient ainsi au passage l'{\'e}galit{\'e} $\widehat{\deg}_{\mathrm{n}}(\overline{E}\oplus_{2}\overline{F})=\widehat{\deg}_{\mathrm{n}}\overline{E}+\widehat{\deg}_{\mathrm{n}}\overline{F}$ lorsque $p=2$. La majoration du degr{\'e} de $\overline{E}\oplus_{\infty}\overline{F}$ s'obtient alors en remarquant $$\frac{\Gamma\left(1+\frac{m+n}{2}\right)}{\Gamma\left(1+\frac{n}{2}\right)\Gamma\left(1+\frac{m}{2}\right)}\le\binom{n+m}{n},$$majoration que l'on utilise aux places r{\'e}elles de $k$. Quant {\`a} la minoration du degr{\'e} de $\overline{E}\oplus_{1}\overline{F}$, on se sert de $\Gamma\left(1+\frac{m+n}{2}\right)\ge\Gamma\left(1+\frac{n}{2}\right)\Gamma\left(1+\frac{m}{2}\right)$ (places r{\'e}elles) et de $\binom{2n+2m}{2n}\le\binom{n+m}{n}^{3}$ (places complexes). Cette derni\`ere estimation peut se d\'emontrer par r\'ecurrence sur $m$ et les pr\'ec\'edentes au moyen de la formule d'Euler pour la fonction Gamma. 
\end{proof}
\begin{rema}En observant que, pour tout $x\in\mathbf{C}^{n}$, on a $\vert x\vert_{2}\le\vert x\vert_{p}$ si $p\le 2$ et $\vert x\vert_{2}\ge\vert x\vert_{p}$ si $p\ge 2$, on obtient ais{\'e}ment les majorations \begin{equation*}\widehat{\deg}_{\mathrm{n}}(\overline{E}\oplus_{p}\overline{F})\le\widehat{\deg}_{\mathrm{n}}\overline{E}+\widehat{\deg}_{\mathrm{n}}\overline{F}\quad\text{si $p\le 2$}\end{equation*}et\begin{equation*}\widehat{\deg}_{\mathrm{n}}\overline{E}+\widehat{\deg}_{\mathrm{n}}\overline{F}\le\widehat{\deg}_{\mathrm{n}}(\overline{E}\oplus_{p}\overline{F})\quad\text{si $p\ge 2$}.\end{equation*} 
\end{rema}
\subsection{Dualit{\'e}}
\begin{prop}\label{propositiondualitedegres}
Il existe une constante absolue $c>0$ telle que, pour tout fibr{\'e} vectoriel ad{\'e}lique $\overline{E}$ de dimension $n$, on a $$\widehat{\deg}_{\mathrm{n}}\overline{E}+\widehat{\deg}_{\mathrm{n}}\overline{E}^{\mathsf{v}}\in[-cn,0]\ .$$Si $\overline{E}$ est un fibr{\'e} ad{\'e}lique hermitien alors $\widehat{\deg}_{\mathrm{n}}\overline{E}^{\mathsf{v}}=-\widehat{\deg}_{\mathrm{n}}\overline{E}$.
\end{prop}
\begin{proof}
On peut supposer $E=k^{n}$. Soit $v$ une place de $k$. S'il existe une matrice $u_{v}\in\GL_{n}(k_{v})$ telle que $(k_{v}^{n},\Vert.\Vert_{v})=u_{v}(k_{v}^{n},\vert.\vert_{2,v})$ alors le quotient des volumes $$\frac{\vol(\mathbf{B}(k_{v}^{n},\Vert.\Vert_{v}))}{\vol(\mathbf{B}(k_{v}^{n},\vert.\vert_{2,v}))}$${\'e}gale $\vert\det u_{v}\vert_{v}^{-n_{v}}$. Il suffit alors d'observer que les normes duales $\Vert.\Vert_{v}^{\mathsf{v}}$ et $\vert.\vert_{2,v}^{\mathsf{v}}$ sont reli{\'e}es par la matrice inverse $u_{v}^{-1}$ et le quotient des volumes pour les espaces duaux est invers{\'e}. En particulier cela assure l'{\'e}galit{\'e} $\widehat{\deg}\overline{E}+\widehat{\deg}\overline{E}^{\mathsf{v}}=0$ si $\overline{E}$ est un fibr{\'e} ad{\'e}lique hermitien. Dans le cas g{\'e}n{\'e}ral, comme il ne reste que les places archim{\'e}diennes {\`a} traiter, on utilise le th{\'e}or{\`e}me~\ref{theoreme23bourgain} en chacune de ces places (avec $n=\dim_{\mathbf{R}}E_{v}$) et le r{\'e}sultat s'ensuit. 
\end{proof}
\subsection{Quotient}
\begin{prop}\label{propositiondegrequotient46}
Soit $\overline{E}_{2}\subseteq\overline{E}_{1}$ deux sous-fibr{\'e}s ad{\'e}liques de $\overline{E}$. Alors on a\begin{equation*}\left\vert\widehat{\deg}_{\mathrm{n}}\overline{E}_{1}-\widehat{\deg}_{\mathrm{n}}\overline{E}_{2}-\widehat{\deg}_{\mathrm{n}}\overline{E_{1}/E_{2}}\right\vert\le\frac{\dim E_{1}}{D}\log\Delta(\overline{E})\ .\end{equation*}\end{prop}
\begin{proof}
Si $\overline{E}$ est un fibr{\'e} ad{\'e}lique hermitien, la m{\'e}trique quotient sur $E_{1}/E_{2}$ co{\"\i}ncide en une place archim{\'e}dienne $v$ avec la m{\'e}trique sur l'orthogonal de $E_{2,v}$ relatif au produit hermitien $(\,,)_{E,v}$ de $E_{v}$. En une place ultram{\'e}trique, la m{\'e}trique quotient est celle d'un suppl{\'e}mentaire quelconque de $E_{2}$ dans $E_{1}$. Ainsi le fibr{\'e} $\overline{E}_{1}$ est isomorphe isom{\'e}triquement {\`a} $\overline{E}_{2}\oplus_{2}\overline{E_{1}/E_{2}}$ et la proposition~\ref{degresommedirecte} donne l'{\'e}galit{\'e} $$\widehat{\deg}_{\mathrm{n}}\overline{E}_{1}=\widehat{\deg}_{\mathrm{n}}\overline{E}_{2}+\widehat{\deg}_{\mathrm{n}}\overline{E_{1}/E_{2}}\ .$$La propri{\'e}t{\'e} {\'e}tant {\'e}tablie dans le cas hermitien, le passage au cas g{\'e}n{\'e}ral s'effectue gr{\^a}ce aux fibr{\'e}s de John et L{\"o}wner associ{\'e}s {\`a} $\overline{E}_{1}$. Comme $\Vert.\Vert_{E_{1},v}\le\vert.\vert_{J(\overline{E}_{1}),v}$ pour toute place $v$ de $k$, cela entra{\^\i}ne $$\widehat{\deg}_{\mathrm{n}}(E_{2},\vert.\vert_{J(\overline{E}_{1}),E_{2}})\le\widehat{\deg}_{\mathrm{n}}\overline{E}_{2}$$et $$\widehat{\deg}_{\mathrm{n}}(E_{1}/E_{2},\vert.\vert_{J(\overline{E}_{1}),E_{1}/E_{2}})\le\widehat{\deg}_{\mathrm{n}}\left(\overline{E_{1}/E_{2}}\right)$$o\`u $\vert.\vert_{J(\overline{E}_{1}),E_{2}}$ (\emph{resp}. $\vert.\vert_{J(\overline{E}_{1}),E_{1}/E_{2}}$) d{\'e}signe la famille des normes de John de $\overline{E}_{1}$ restreintes {\`a} $E_{2}$ (\emph{resp}. normes quotient sur $E_{1}/E_{2}$ issues des normes de John de $\overline{E}_{1}$). L'on obtient ainsi\begin{equation*}\begin{split}&\widehat{\deg}_{\mathrm{n}}\overline{E}_{1}-\widehat{\deg}_{\mathrm{n}}\overline{E}_{2}-\widehat{\deg}_{\mathrm{n}}\overline{E_{1}/E_{2}}\\ &\quad\le\widehat{\deg}_{\mathrm{n}}\overline{E}_{1}-\left(\widehat{\deg}_{\mathrm{n}}(E_{2},\vert.\vert_{J(\overline{E}_{1}),E_{2}})+\widehat{\deg}_{\mathrm{n}}(E_{1}/E_{2},\vert.\vert_{J(\overline{E}_{1}),E_{1}/E_{2}})\right)\end{split}\end{equation*}et l'expression dans cette derni{\`e}re parenth{\`e}se {\'e}gale $\widehat{\deg}_{\mathrm{n}}J(\overline{E}_{1})$ car $J(\overline{E}_{1})$ est un fibr{\'e} ad{\'e}lique hermitien. On raisonne de la m{\^e}me fa{\c c}on avec $L(\overline{E_{1}})$ pour la majoration en sens inverse. Gr{\^a}ce {\`a}~\eqref{eqref:estimationdifference}, on obtient la borne \begin{equation*}\begin{split}\left\vert\widehat{\deg}_{\mathrm{n}}\overline{E}_{1}-\widehat{\deg}_{\mathrm{n}}\overline{E}_{2}-\widehat{\deg}_{\mathrm{n}}\overline{E_{1}/E_{2}}\right\vert &\le\widehat{\deg}_{\mathrm{n}}L(\overline{E}_{1})-\widehat{\deg}_{\mathrm{n}}J(\overline{E}_{1})\\ &\le\frac{\dim E_{1}}{D}\log\Delta(\overline{E}_{1})\end{split}\end{equation*}et l'on conclut avec $\Delta(\overline{E}_{1})\le\Delta(\overline{E})$.
\end{proof}
\subsection{Somme}
Soit $\overline{E}$ un fibr{\'e} vectoriel ad{\'e}lique sur $\spec k$ et $E_{1},E_{2}$ deux sous-espaces vectoriels de $E$, de dimensions respectives $n_{1}$ et $n_{2}$. Chacun des quatre espaces vectoriels $E_{1}$, $E_{2}$, $E_{1}+E_{2}$, $E_{1}\cap E_{2}$ est muni de la structure de sous-fibr{\'e} vectoriel ad{\'e}lique induite par $\overline{E}$.
\begin{prop}\label{prop:sommededegrefibres}
Avec les donn{\'e}es ci-dessus, la diff{\'e}rence des degr{\'e}s ad\'eliques normalis{\'e}s \begin{equation*}\widehat{\deg}_{\mathrm{n}}(\overline{E_{1}+E_{2}})+\widehat{\deg}_{\mathrm{n}}\overline{E_{1}\cap E_{2}}-\widehat{\deg}_{\mathrm{n}}\overline{E}_{1}-\widehat{\deg}_{\mathrm{n}}\overline{E}_{2}\end{equation*}est sup{\'e}rieure {\`a}\begin{equation*}-\frac{(n_{1}+n_{2})}{D}\log\Delta(\overline{E_{1}+E_{2}})\ \cdotp\end{equation*}
\end{prop}
\begin{proof}L'application naturelle $\iota:E_{2}/(E_{1}\cap E_{2})\to(E_{1}+E_{2})/E_{1}$ v{\'e}rifie, pour toute place $v$ de $k$ et tout $x\in k_{v}$, $$\Vert\iota(x)\Vert_{(E_{1}+E_{2})/E_{1},v}\le\Vert x\Vert_{E_{2}/(E_{1}\cap E_{2}),v}\ .$$Par cons{\'e}quent, la boule unit{\'e} du premier quotient est incluse dans celle du second et l'on a $$\widehat{\deg}_{\mathrm{n}}\left(\overline{E_{2}/(E_{1}\cap E_{2})}\right)\le\widehat{\deg}_{\mathrm{n}}\left(\overline{(E_{1}+E_{2})/E_{1}}\right)\ \cdotp$$Si $\overline{E_{1}+E_{2}}$ est hermitien, la proposition~\ref{propositiondegrequotient46} permet de conclure puisque dans ce cas le degr{\'e} d'un quotient est la diff{\'e}rence des degr{\'e}s. Dans le cas g{\'e}n{\'e}ral, on consid{\`e}re $\varepsilon>0$ et une norme hermitienne $\vert.\vert_{\varepsilon,v}$ sur $(E_{1}+E_{2})\otimes_{k}k_{v}$, ce qui permet de d{\'e}finir le fibr{\'e} ad{\'e}lique hermitien $\overline{(E_{1}+E_{2})}_{\varepsilon}$ (voir~\eqref{eqref:important}). On a alors \begin{equation*}\begin{split}\widehat{\deg}_{\mathrm{n}}\overline{E}_{1}+\widehat{\deg}_{\mathrm{n}}\overline{E}_{2}&\le\widehat{\deg}_{\mathrm{n}}(E_{1},\vert.\vert_{\varepsilon})+\widehat{\deg}_{\mathrm{n}}(E_{2},\vert.\vert_{\varepsilon})\\ & \le\widehat{\deg}_{\mathrm{n}}(E_{1}+E_{2},\vert.\vert_{\varepsilon})+\widehat{\deg}_{\mathrm{n}}(E_{1}\cap E_{2},\vert.\vert_{\varepsilon})\ .\end{split}\end{equation*}La premi\`ere in\'egalit\'e qui vient apr\`es~\eqref{eqref:estimationdifference} montre que le degr{\'e} ad{\'e}lique normalis{\'e} de $\overline{(E_{1}+E_{2})}_{\varepsilon}$ est plus petit que $$\widehat{\deg}_{\mathrm{n}}\overline{E_{1}+E_{2}}+\dim(E_{1}+E_{2})\left(\log(1+\varepsilon)+\frac{1}{D}\log\Delta(\overline{E_{1}+E_{2}})\right)\ .$$Il en est de m{\^e}me pour le degr{\'e} de $(E_{1}\cap E_{2},\vert.\vert_{\varepsilon})$ en rempla{\c c}ant $\dim(E_{1}+E_{2})$ par $\dim E_{1}\cap E_{2}$. La proposition~\ref{prop:sommededegrefibres} s'ensuit en faisant tendre $\varepsilon$ vers $0$.
\end{proof}
En notant $H(\overline{E})=\exp{\left\{-\widehat{\deg}_{\mathrm{n}}\overline{E}\right\}}$, cet {\'e}nonc{\'e} peut se voir comme une g{\'e}n{\'e}ralisation de l'in{\'e}galit{\'e} $H(\overline{E_{1}+E_{2}})H(\overline{E_{1}\cap E_{2}})\le H(\overline{E_{1}})H(\overline{E_{2}})$ de W.~Schmidt~\cite{Schmidt2}, T.~Struppeck \& J.D.~Vaaler~\cite{Struppeck}, lorsque les m{\'e}triques sont hermitiennes.

\section{Th{\'e}orie des pentes et pentes maximales}

\begin{defi}\label{definitionpenteadelique} La \emph{pente ad{\'e}lique} de $\overline{E}$ est \begin{equation}\widehat{\mu}(\overline{E}):=\frac{\widehat{\deg}\overline{E}}{\dim E}\cdotp\end{equation}Par convention, si $E=\{0\}$, on pose $\widehat{\mu}(\overline{E}):=-\infty$. De plus, lorsque l'on divise ce nombre r{\'e}el par le degr{\'e} $D$ du corps $k$, on parle de pente ad{\'e}lique \emph{normalis{\'e}e} que l'on note $\widehat{\mu}_{\mathrm{n}}(\overline{E})$. 
\end{defi}
Cette quantit\'e est un avatar ad{\'e}lique et logarithmique du \emph{rayon volumique} $$\left(\frac{\vol(C)}{\vol(b_{n}^{2})}\right)^{1/n}$$d'un corps convexe $C$ de $\mathbf{R}^{n}$, utilis\'e dans l'\'etude des espaces de Banach. \par Voici quelques propri{\'e}t{\'e}s que l'on d{\'e}duit imm{\'e}diatement de celles du degr{\'e} ad{\'e}lique, d{\'e}montr{\'e}es dans le paragraphe pr{\'e}c{\'e}dent. Dans la liste qui suit, $\overline{E}$ d{\'e}signe un fibr{\'e} vectoriel ad{\'e}lique de dimension $n$, $\overline{E}_{\varepsilon}=(E,\vert.\vert_{\varepsilon})$ un fibr{\'e} ad{\'e}lique hermitien associ{\'e} {\`a} $\overline{E}$ (voir~\eqref{eqref:important}) et $K$ est une extension finie de $k$.\par
\emph{Propri{\'e}t{\'e}s}~:
 \begin{enumerate}
\item[1)] Pour tout sous-espace vectoriel $E'$ de $E$, on a \begin{equation*}\label{encadrementpentessimples}-\log(1+\varepsilon)-\frac{1}{D}\log\Delta(\overline{E})\le\widehat{\mu}_{\mathrm{n}}(\overline{E'})-\widehat{\mu}_{\mathrm{n}}(E',\vert.\vert_{\varepsilon})\le 0\ .\end{equation*}
\item[2)] Il existe une constante absolue $c>0$ telle que $-c\le\widehat{\mu}_{\mathrm{n}}(\overline{E})+\widehat{\mu}_{\mathrm{n}}(\overline{E}^{\mathsf{v}})\le 0$.
\item[3)] $\left\vert\widehat{\mu}_{\mathrm{n}}(\overline{E}_{K})-\widehat{\mu}_{\mathrm{n}}(\overline{E})\right\vert\le 2\log(2)$.
\end{enumerate}
Si de plus $\overline{E}$ est un fibr{\'e} ad{\'e}lique hermitien alors $\widehat{\mu}_{\mathrm{n}}(\overline{E}^{\mathsf{v}})=-\widehat{\mu}_{\mathrm{n}}(\overline{E})$ et $\widehat{\mu}_{\mathrm{n}}(\overline{E}_{K})=\widehat{\mu}_{\mathrm{n}}(\overline{E})$.\par {\`A} ces propri{\'e}t{\'e}s s'ajoute une estimation de la pente d'un produit tensoriel.
\begin{prop}\label{penteproduittensorielfibre52}Soit $\ell\in\mathbf{N}\setminus\{0\}$. Pour tout $i\in\{1,\ldots,\ell\}$, soit $\overline{E}_{i}$ un fibr{\'e} vectoriel ad{\'e}lique sur $\spec k$. Soit $\alpha$ une norme tensorielle ad\'elique d'ordre $\ell$. Si tous les $\overline{E}_{i}$ sauf au plus un sont de dimension $1$ alors $\widehat{\mu}(\otimes_{\alpha,i=1}^{\ell}\overline{E}_{i})=\sum_{i=1}^{\ell}{\widehat{\mu}(\overline{E}_{i})}$. En dimension quelconque, si $\alpha$ est hermitienne alors \begin{equation*}\left\vert\widehat{\mu}(\otimes_{\alpha,i=1}^{\ell}\overline{E}_{i})-\sum_{i=1}^{\ell}{\widehat{\mu}(\overline{E}_{i})}\right\vert\le\sum_{i=1}^{\ell}{\log\Delta(\overline{E}_{i})}\ .\end{equation*}
\end{prop}On d{\'e}duit aussi de cette in{\'e}galit{\'e} la relation $\widehat{\mu}(\otimes_{\alpha,i=1}^{\ell}\overline{E}_{i})=\sum_{i=1}^{\ell}{\widehat{\mu}(\overline{E}_{i})}$ lorsque tous les $\overline{E}_{i}$ et $\alpha$ sont hermitiens. Toutefois, la preuve que nous proposons commence par {\'e}tablir ce cas.
\begin{proof} Le premier cas d'\'egalit\'e lorsque par exemple $\dim E_{i}=1$ pour $2\le i\le\ell$ s'obtient en observant que l'application $x\mapsto x\otimes e_{2}\otimes\cdots\otimes e_{\ell}$ de $(E,(\Vert.\Vert_{E_{1},v}\times\prod_{i=2}^{\ell}{\Vert e_{i}\Vert_{E_{i},v}}))$ dans $\otimes_{\alpha,i=1}^{\ell}\overline{E}_{i}$ est un isomorphisme isom\'etrique ($e_{i}$ est une $k$-base quelconque de $E_{i}$). Pla{\c c}ons nous alors en dimension quelconque, avec $\alpha$ hermitienne. Supposons $\ell=2$ dans un premier temps. Si, en une place $v$ de $k$, les normes $\Vert.\Vert_{E_{1},v}$ et $\Vert.\Vert_{E_{2},v}$ sont simultan{\'e}ment hermitiennes ou ultram{\'e}triques, l'isomorphisme canonique $$\det(E_{1,v}\otimes E_{2,v})\simeq\det(E_{1,v})^{\otimes m}\otimes\det(E_{2,v})^{\otimes n}$$est une isom{\'e}trie. Par cons{\'e}quent, lorsque $\overline{E}_{1}$ et $\overline{E}_{2}$ sont des fibr{\'e}s ad{\'e}liques hermitiens, on a $$\widehat{\deg}\det(\overline{E}_{1}\otimes\overline{E}_{2})=m\widehat{\deg}\det(\overline{E}_{1})+n\widehat{\deg}\det(\overline{E}_{2})$$et le lemme~\ref{lemme4226} entra{\^\i}ne l'{\'e}galit{\'e} $\widehat{\mu}(\overline{E}_{1}\otimes\overline{E}_{2})=\widehat{\mu}(\overline{E}_{1})+\widehat{\mu}(\overline{E}_{2})$. Le r\'esultat pour le produit de $\ell$ fibr\'es hermitiens s'en d\'eduit par r\'ecurrence sur $\ell$. Une fois celui-ci {\'e}tabli, le passage au cas g{\'e}n{\'e}ral s'effectue au moyen des fibr{\'e}s de John et L{\"o}wner associ{\'e}s aux $\overline{E}_{i}$. Gr{\^a}ce au lemme~\ref{lemmecomptrois}, on a $$\otimes_{i=1}^{\ell}{L(\overline{E}_{i})}\preceq\otimes_{\alpha,i=1}^{\ell}{\overline{E}_{i}}\preceq\otimes_{i=1}^{\ell}{J(\overline{E}_{i})}\ .$$On en d\'eduit un encadrement de $\widehat{\mu}(\otimes_{\alpha,i=1}^{\ell}{\overline{E}_{i}})$ \textit{via} le lemme~\ref{lemmesimplee}. L'hypoth\`ese $\alpha$ hermitienne permet d'utiliser le cas hermitien \'etabli juste avant et les formules de la proposition~\ref{prop:comparaisonjohnadelique} conduisent \`a l'estimation \begin{equation*}\left\vert\widehat{\mu}(\otimes_{\alpha,i=1}^{\ell}{\overline{E}_{i}})-\sum_{i=1}^{\ell}{\widehat{\mu}(\overline{E}_{i})}\right\vert\le\sum_{i=1}^{\ell}{\left\{\widehat{\mu}(L(\overline{E}_{i}))-\widehat{\mu}(J(\overline{E}_{i}))\right\}}\ .\end{equation*}On conclut avec~\eqref{eqref:estimationdifference}.
\end{proof}
\subsection{Existence de la pente maximale} 
\begin{prop}\label{propofinitude}
Soit $\overline{E}$ un fibr{\'e} vectoriel ad{\'e}lique sur $\spec k$. Il existe une constante $c(\overline{E},k)$ telle que, pour tout sous-fibr{\'e} vectoriel ad{\'e}lique $\overline{F}$ de $\overline{E}$, on a $\widehat{\mu}_{\mathrm{n}}(\overline{F})\le c(\overline{E},k)$. Plus pr{\'e}cis{\'e}ment, pour tout nombre r{\'e}el $a$, il n'existe qu'un nombre \emph{fini} de sous-fibr{\'e}s vectoriels ad{\'e}liques $\overline{F}$ de $\overline{E}$ tels que $\widehat{\mu}_{\mathrm{n}}(\overline{F})\ge a$. 
\end{prop}
\begin{proof}
Soit $\overline{F}\subseteq\overline{E}$ et $m=\dim F$. Si $\overline{E}$ est un fibr{\'e} ad{\'e}lique hermitien alors il en est de m{\^e}me pour $\overline{F}$ et l'on a \begin{equation*}\widehat{\mu}_{\mathrm{n}}(\overline{F})=\frac{\widehat{\deg}_{\mathrm{n}}\det\overline{F}}{\dim F}=-\frac{h_{\bigwedge^{m}\overline{E}}(f_{1}\wedge\cdots\wedge f_{m})}{m}\end{equation*}o\`u $f_{1},\ldots,f_{m}$ est une $k$-base de $\overline{F}$ et $h_{\bigwedge^{m}\overline{E}}$ est la hauteur sur les {\'e}l{\'e}ments de $\bigwedge^{m}\overline{E}$ donn\'ee par la d{\'e}finition~\ref{definition48hauteur}. En vertu de la proposition~\ref{propositionsurlahauteur}, il existe une constante $c(m,\overline{E},k)$ telle que $$\forall x\in\bigwedge^{m}E,\quad h_{\bigwedge^{m}\overline{E}}(x)\ge c(m,\overline{E},k)\ ;$$de la sorte on obtient la majoration $\widehat{\mu}_{\mathrm{n}}(\overline{F})\le-c(m,\overline{E},k)/m$ et la constante $c(\overline{E},k):=\max{\left\{-c(m,\overline{E},k)D/m\,;\ 1\le m\le\dim E\right\}}$ convient dans ce cas. En r{\'e}alit{\'e}, si $a\le\widehat{\mu}_{\mathrm{n}}(\overline{F})$ alors la hauteur de $f_{1}\wedge\cdots\wedge f_{m}$ est major{\'e}e et la proposition~\ref{propositionsurlahauteur} entra{\^\i}ne la finitude du nombre de sous-espaces $F$ possibles. Si $\overline{E}$ est un fibr{\'e} vectoriel ad{\'e}lique quelconque, on utilise les m{\'e}triques de L{\"o}wner, qui satisfont aux in{\'e}galit{\'e}s $\vert.\vert_{L(\overline{E}),v}\le\Vert.\Vert_{E,v}$ en toute place $v$ de $k$, pour obtenir $\widehat{\deg}_{\mathrm{n}}(\overline{F})\le\widehat{\deg}_{\mathrm{n}}(F,\vert.\vert_{L(\overline{E})})$. La conclusion d{\'e}coule alors du cas hermitien.
\end{proof}
Cela justifie la d{\'e}finition suivante.
\begin{defi}\label{definitiondelapentemaximale}
La \emph{pente maximale} d'un fibr{\'e} vectoriel ad{\'e}lique $\overline{E}$, not{\'e}e $\widehat{\mu}_{\mathrm{max}}(\overline{E})$, est le nombre r{\'e}el \begin{equation*}\widehat{\mu}_{\mathrm{max}}(\overline{E}):=\max{\left\{\widehat{\mu}_{\mathrm{n}}(\overline{F})\ ;\ \overline{F}\subseteq\overline{E}\right\}}\ \cdotp\end{equation*}
\end{defi}On observe que si $\overline{E}_{1}\preceq\overline{E}_{2}$ alors $\widehat{\mu}_{\mathrm{max}}(\overline{E}_{2})\le\widehat{\mu}_{\mathrm{max}}(\overline{E}_{1})$. La premi{\`e}re propri{\'e}t{\'e} de la pente normalis{\'e}e donn{\'e}e {\`a} la suite de la d\'efinition~\ref{definitionpenteadelique} entra{\^\i}ne \begin{equation*}-\log(1+\varepsilon)-\frac{1}{D}\log\Delta(\overline{E})\le\widehat{\mu}_{\mathrm{max}}(\overline{E})-\widehat{\mu}_{\mathrm{max}}(\overline{E}_{\varepsilon})\le 0\ ,\end{equation*}pour tout fibr{\'e} vectoriel ad{\'e}lique $\overline{E}$ et tout $\varepsilon>0$.
\begin{rema}\label{remarquecinqcinq}
De la m{\^e}me mani\`ere, on peut consid\'erer la \emph{pente minimale} d'un fibr\'e vectoriel ad\'elique $\overline{E}$ d\'efinie par $\widehat{\mu}_{\mathrm{min}}(\overline{E}):=\min{\left\{\widehat{\mu}_{\mathrm{n}}(\overline{E/F})\right\}}$ o\`u $F$ parcourt les sous-espaces vectoriels stricts de $E$. Il s'agit bien d'un minimum car il n'y a qu'un nombre fini de quotients $E/F$ de pentes plus petites qu'un nombre r\'eel donn\'e (on se ram\`ene \`a la pente maximale \textit{via} la proposition~\ref{propositiondegrequotient46}). On dispose aussi de l'encadrement  \begin{equation*}-\log(1+\varepsilon)-\frac{1}{D}\log\Delta(\overline{E})\le\widehat{\mu}_{\mathrm{min}}(\overline{E})-\widehat{\mu}_{\mathrm{min}}(\overline{E}_{\varepsilon})\le 0\ .\end{equation*}
\end{rema}
La dualit\'e marque le lien avec la pente maximale.
\begin{lemma}\label{lemmacomparaisonmumaxmumin}
Soit $\overline{E}$ un fibr\'e vectoriel ad\'elique sur $\spec k$. On a \begin{equation*}\left\vert\widehat{\mu}_{\mathrm{max}}(\overline{E^{\mathsf{v}}})+\widehat{\mu}_{\mathrm{min}}(\overline{E})\right\vert\le\frac{1}{D}\log\Delta(\overline{E})\ .\end{equation*}  
\end{lemma}
\begin{proof}Supposons dans un premier temps que $\overline{E}$ est hermitien. Soit $\overline{F}$ un sous-fibr\'e vectoriel strict de $\overline{E}$. L'espace vectoriel dual $(E/F)^{\mathsf{v}}$ s'injecte dans $E^{\mathsf{v}}$. La d\'efinition de la pente maximale de $\overline{E^{\mathsf{v}}}$ et les propositions~\ref{propositiondualitedegres} et~\ref{propositiondegrequotient46} dans le cas hermitien donnent \begin{equation*}\dim(E/F)\,\widehat{\mu}_{\mathrm{max}}(\overline{E^{\mathsf{v}}})\ge\widehat{\deg}_{\mathrm{n}}\overline{(E/F)^{\mathsf{v}}}=-\widehat{\deg}_{\mathrm{n}}\overline{E/F}\end{equation*}puis $\widehat{\mu}_{\mathrm{max}}(\overline{E^{\mathsf{v}}})+\widehat{\mu}_{\mathrm{n}}(\overline{E/F})\ge 0$ et donc $\widehat{\mu}_{\mathrm{max}}(\overline{E^{\mathsf{v}}})+\widehat{\mu}_{\mathrm{min}}(\overline{E})\ge 0$. On montre de la m{\^e}me mani\`ere l'in\'egalit\'e en sens inverse. On a donc $\widehat{\mu}_{\mathrm{max}}(\overline{E^{\mathsf{v}}})=-\widehat{\mu}_{\mathrm{min}}(\overline{E})$ si $\overline{E}$ est hermitien. Dans le cas g\'en\'eral, on utilise le fibr\'e vectoriel $\overline{E}_{\varepsilon}$ associ\'e \`a $\overline{E}$. Compte tenu du fait que la dualit\'e renverse les in\'egalit\'es pour les normes~: \begin{equation*}\forall G\subseteq E^{\mathsf{v}},\quad (G,(\mathrm{d}(E\otimes_{k}k_{v},\ell^{2}_{n,k_{v}})^{-1}(1+\varepsilon)^{-1}\vert.\vert_{\varepsilon,v}^{\mathsf{v}})_{v})\preceq(G,\Vert.\Vert_{E}^{\mathsf{v}})\preceq(G,\vert.\vert_{\varepsilon}^{\mathsf{v}}),\end{equation*}on a \begin{equation*}0\le\widehat{\mu}_{\mathrm{max}}(\overline{E^{\mathsf{v}}})-\widehat{\mu}_{\mathrm{max}}(\overline{E^{\mathsf{v}}_{\varepsilon}})\le\frac{1}{D}\log\Delta(\overline{E})+\log(1+\varepsilon)\ .\end{equation*}On somme avec l'estimation du m{\^e}me type entre pentes minimales qui pr\'ec\`ede l'\'enonc\'e du lemme~\ref{lemmacomparaisonmumaxmumin}, puis l'on fait tendre $\varepsilon$ vers $0$ et le lemme s'en d\'eduit. 
\end{proof}
\begin{proprietes}\label{proprietes55}
Soit $\overline{E}$, $\overline{F}$ des fibr{\'e}s vectoriels ad{\'e}liques sur $\spec k$.
\begin{enumerate}
\item[1)] Si $\dim F=1$ alors, pour toute norme tensorielle ad\'elique $\alpha$ d'ordre $2$, on a \begin{equation*}\widehat{\mu}_{\mathrm{max}}(\overline{E}\otimes_{\alpha}\overline{F})=\widehat{\mu}_{\mathrm{max}}(\overline{E})+\widehat{\deg}_{\mathrm{n}}\overline{F}\ .\end{equation*}
\item[2)] Pour la somme directe hermitienne $\overline{E}\oplus_{2}\overline{F}$, on a $$0\le\widehat{\mu}_{\mathrm{max}}(\overline{E}\oplus_{2}\overline{F})-\max{\left\{\widehat{\mu}_{\mathrm{max}}(\overline{E}),\widehat{\mu}_{\mathrm{max}}(\overline{F})\right\}}\le\frac{1}{D}\log\max{\left\{\Delta(\overline{E}),\Delta(\overline{F})\right\}}\ \cdotp$$
\item[3)] Dans le cas g{\'e}n{\'e}ral, on a \begin{equation*}0\le\widehat{\mu}_{\mathrm{max}}(\overline{E}\oplus_{\varsigma}\overline{F})-\max{\left\{\widehat{\mu}_{\mathrm{max}}(\overline{E}),\widehat{\mu}_{\mathrm{max}}(\overline{F})\right\}}\le\frac{1}{D}\log\left(\sqrt{2}\max{\left\{\Delta(\overline{E}),\Delta(\overline{F})\right\}}\right)\end{equation*}(la somme directe $\overline{E}\oplus_{\varsigma}\overline{F}$ est celle d{\'e}finie au paragraphe~\ref{paragraphetroistrois}).
\end{enumerate} 
\end{proprietes}
\begin{proof}
\noindent
1) Soit $\overline{E'}$ un sous-fibr{\'e} vectoriel ad{\'e}lique de $\overline{E}$. D'apr{\`e}s la proposition~\ref{penteproduittensorielfibre52}, la pente normalis\'ee de $\overline{E'}\otimes_{\alpha}\overline{F}$ est la somme de $\widehat{\mu}_{\mathrm{n}}(\overline{E'})$ et de $\widehat{\mu}_{\mathrm{n}}(\overline{F})=\widehat{\deg}_{\mathrm{n}}\overline{F}$, car $F$ est une droite. De plus, d'apr\`es la proposition~\ref{propositioncompnormes}, on a $(E'\otimes F,\alpha(\cdot;\overline{E},\overline{F}))\preceq\overline{E'}\otimes_{\alpha}\overline{F}$ et le lemme~\ref{lemmesimplee} entra{\^\i}ne $$\widehat{\mu}_{\mathrm{n}}(\overline{E'}\otimes_{\alpha}\overline{F})\le\widehat{\mu}_{\mathrm{max}}(\overline{E}\otimes_{\alpha}\overline{F})\ .$$On a donc $\widehat{\mu}_{\mathrm{n}}(\overline{E'})\le\widehat{\mu}_{\mathrm{max}}(\overline{E}\otimes_{\alpha}\overline{F})-\widehat{\deg}_{\mathrm{n}}\overline{F}$ et, en faisant varier $\overline{E'}$, on obtient $$\widehat{\mu}_{\mathrm{max}}(\overline{E})\le\widehat{\mu}_{\mathrm{max}}(\overline{E}\otimes_{\alpha}\overline{F})-\widehat{\deg}_{\mathrm{n}}\overline{F}\ .$$Si l'on applique maintenant cette majoration {\`a} $\overline{E}\otimes_{\alpha}\overline{F}$ et $\overline{F^\mathsf{v}}$ au lieu, respectivement, de $\overline{E}$ et $\overline{F}$, on obtient l'in{\'e}galit{\'e} en sens inverse, gr{\^a}ce {\`a} la propri{\'e}t{\'e} $\widehat{\deg}_{\mathrm{n}}\overline{F^\mathsf{v}}=-\widehat{\deg}_{\mathrm{n}}\overline{F}$, valide car $\overline{F}$ est automatiquement hermitien, et car $(\overline{E}\otimes_{\alpha}\overline{F})\otimes_{\alpha}\overline{F^{\mathsf{v}}}$ est isomorphe isom\'etriquement \`a $\overline{E}$.
\par\noindent 
2) La positivit{\'e} de la diff{\'e}rence des pentes maximales est imm{\'e}diate en notant qu'un sous-espace vectoriel de $E$ ou $F$ est aussi un sous-espace de $E\oplus F$, avec compatibilit{\'e} des normes. Ceci s'applique {\'e}galement pour $\overline{E}\oplus_{\varsigma}\overline{F}$ puisque la norme $\Vert.\Vert_{\overline{E}\oplus_{\varsigma}\overline{F},v}$ restreinte {\`a} $E$ ou $F$ est celle qui est naturellement sur ces espaces.\par\indent Pour la majoration, supposons dans un premier temps que $\overline{E}$ et $\overline{F}$ sont \emph{hermitiens}. Soit $\overline{H}$ un sous-fibr{\'e} vectoriel ad{\'e}lique de $\overline{E}\oplus_{2}\overline{F}$. Posons $H_{1}:=H\cap E$. Par d\'efinition des m\'etriques quotient, on a $\overline{H/H_{1}}=\overline{H}/\overline{H_{1}}$ et donc, puisque ces fibr\'es ad\'eliques sont hermitiens, $$\widehat{\deg}_{\mathrm{n}}\overline{H}=\widehat{\deg}_{\mathrm{n}}\overline{H}_{1}+\widehat{\deg}_{\mathrm{n}}(\overline{H/H_{1}})\ .$$Le premier terme $\widehat{\deg}_{\mathrm{n}}\overline{H}_{1}$ est naturellement inf{\'e}rieur {\`a} $(\dim H_{1})\widehat{\mu}_{\mathrm{max}}(\overline{E})$. Pour le second terme, consid{\'e}rons l'injection $\iota:H/H_{1}\hookrightarrow F$, $\iota(x\oplus y\mod H_{1})=y$. L'application $\iota$ n'est pas n{\'e}cessairement une isom{\'e}trie en toute place $v$ de $k$ (et donc $\overline{H/H_{1}}$ n'est pas en g{\'e}n{\'e}ral un sous-fibr{\'e} ad{\'e}lique de $\overline{F}$). N{\'e}anmoins, comme $\Vert x\oplus y\Vert_{E\oplus_{2} F,v}\ge\Vert y\Vert_{F,v}$ pour tout $x\oplus y\in H_{v}$, on a \begin{equation*}
\widehat{\deg}_{\mathrm{n}}(\overline{H/H_{1}})\le\widehat{\deg}_{\mathrm{n}}\left(\iota(H/H_{1}),\Vert.\Vert_{F}\right)\le\left(\dim H/H_{1}\right)\widehat{\mu}_{\mathrm{max}}(\overline{F})\ .
\end{equation*}On en d{\'e}duit \begin{equation*}\widehat{\mu}_{\mathrm{n}}(\overline{H})=\frac{\widehat{\deg}_{\mathrm{n}}(\overline{H})}{\dim H}\le\frac{(\dim H_{1})\widehat{\mu}_{\mathrm{max}}(\overline{E})+\left(\dim H/H_{1}\right)\widehat{\mu}_{\mathrm{max}}(\overline{F})}{\dim H}\end{equation*}puis $\widehat{\mu}_{\mathrm{max}}(\overline{E}\oplus_{2}\overline{F})\le\max{\left\{\widehat{\mu}_{\mathrm{max}}(\overline{E}),\widehat{\mu}_{\mathrm{max}}(\overline{F})\right\}}$, et avec la remarque du d\'ebut de ce point 2), on obtient l'\'egalit\'e souhait{\'e}e dans le cas hermitien.\par Revenons maintenant au cas g{\'e}n{\'e}ral de deux fibr{\'e}s vectoriels ad{\'e}liques quelconques. Soit $\varepsilon>0$ et $\overline{E}_{\varepsilon}=(E,\vert.\vert_{\varepsilon,E,v})$ (\emph{resp}. $\overline{F}_{\varepsilon}=(F,\vert.\vert_{\varepsilon,F,v})$) un fibr{\'e} ad{\'e}lique hermitien associ{\'e} {\`a} $\overline{E}$ (\emph{resp}. $\overline{F}$), comme celui introduit {\`a} la suite de~\eqref{eqref:important}. Par d{\'e}finition, pour tous $x\in E_{v}$, $y\in F_{v}$, on a $\vert x\vert_{\varepsilon,E,v}^{2}+\vert y\vert_{\varepsilon,F,v}^{2}\le\Vert(x,y)\Vert_{\overline{E}\oplus_{2}\overline{F},v}^{2}$ et donc
\begin{equation*}\widehat{\mu}_{\mathrm{max}}(\overline{E}\oplus_{2}\overline{F})\le\widehat{\mu}_{\mathrm{max}}\left(\overline{E}_{\varepsilon}\oplus_{2}\overline{F}_{\varepsilon}\right)=\max{\left\{\widehat{\mu}_{\mathrm{max}}\left(\overline{E}_{\varepsilon}\right),\widehat{\mu}_{\mathrm{max}}\left(\overline{F}_{\varepsilon}\right)\right\}}\ \cdotp\end{equation*}La pente maximale de $\overline{E}_{\varepsilon}$ est {\`a} son tour plus petite que $\widehat{\mu}_{\mathrm{max}}(\overline{E})+\frac{1}{D}\log\Delta(\overline{E})+\log(1+\varepsilon)$. Le m{\^e}me type de majoration pour $\widehat{\mu}_{\mathrm{max}}\left(\overline{F}_{\varepsilon}\right)$ puis $\varepsilon\to 0$ donnent l'in{\'e}galit{\'e} souhait{\'e}e.
\par\noindent
3) En d{\'e}finissant $\varsigma$, on avait observ{\'e} que $\Vert.\Vert_{\overline{E}\oplus_{\infty}\overline{F}}\le\Vert.\Vert_{\overline{E}\oplus_{\varsigma}\overline{F}}$. Comme $(\sqrt{2})^{-1}\Vert.\Vert_{\overline{E}\oplus_{2}\overline{F}}\le\Vert.\Vert_{\overline{E}\oplus_{\infty}\overline{F}}$, le probl{\`e}me se ram{\`e}ne {\`a} majorer la pente maximale de $\overline{E}\oplus_{2}\overline{F}$, ce que nous venons de faire.
\end{proof}
\begin{rema}Soit $N\in\mathbf{N}\setminus\{0\}$ et $\overline{E_{1}},\ldots,\overline{E_{N}}$ des fibr\'es vectoriels ad\'eliques sur $\spec k$. Alors on a \begin{equation*}0\le\widehat{\mu}_{\mathrm{max}}\left(\overline{E_{1}}\oplus_{2}\cdots\oplus_{2}\overline{E}_{N}\right)-\max_{1\le i\le N}{\left\{\widehat{\mu}_{\mathrm{max}}(\overline{E_{i}})\right\}}\le\frac{1}{D}\log\max_{1\le i\le N}{\left\{\Delta(\overline{E_{i}})\right\}}\ \cdotp\end{equation*}Et, pour une somme directe plus g\'en\'erale comme dans le point $3)$ des propri\'et\'es~\ref{proprietes55}, on doit rajouter un $\sqrt{N}$ devant le maximum des $\Delta(\overline{E_{i}})$. Ces g\'en\'eralisations se d\'emontrent en approchant chacun des fibr\'es vectoriels ad\'eliques $\overline{E_{i}}$ par un fibr\'e hermitien $\overline{E_{i}}_{,\varepsilon}$. La pente maximale de la somme directe hermitienne des $\overline{E_{i}}_{,\varepsilon}$  est le maximum de chacune des pentes maximales $\widehat{\mu}_{\mathrm{max}}(\overline{E_{i}}_{,\varepsilon})$ et on proc\`ede alors comme dans la d\'emonstration ci-dessus. Le nombre $\sqrt{N}$ provient de la comparaison entre les normes de $\ell^{2}_{N}$ et $\ell^{\infty}_{N}$. 
\end{rema}
Au \S~\ref{subsectionpenteproduitsymetrique}, nous obtiendrons {\'e}galement une estimation de la pente maximale de la $\ell^{\text{{\`e}me}}$ puissance sym{\'e}trique $S^{\ell}(\overline{E})$. 
\subsection{Autres pentes}
Soit $\overline{E}$ un fibr{\'e} vectoriel ad{\'e}lique sur $\spec k$, de dimension $n$. La proposition~\ref{propofinitude} implique que l'intersection de l'ensemble $\{(\dim F,\widehat{\deg}_{\mathrm{n}}\overline{F})\,;\ \overline{F}\subseteq\overline{E}\}$ avec tout demi-plan sup{\'e}rieur est \emph{fini}. L'enveloppe convexe de cet ensemble est d{\'e}limit{\'e}e par une fonction affine par morceaux et \emph{concave} $P_{\overline{E}}:[0,n]\to\mathbf{R}$. Les sommets du graphe de $P_{\overline{E}}$ sont atteints par des sous-fibr{\'e}s ad{\'e}liques de $\overline{E}$. On a $P_{\overline{E}}(0)=0$ et $P_{\overline{E}}(n)=\widehat{\deg}_{\mathrm{n}}\overline{E}$. On note que si $\overline{F}\subseteq\overline{E}$ alors, pour tout $x\in[0,\dim F]$, on a $P_{\overline{F}}(x)\le P_{\overline{E}}(x)$.
\begin{defi}
Pour $i\in\{1,\ldots,n\}$, la $i^{\text{{\`e}me}}$ pente normalis{\'e}e de $\overline{E}$, not{\'e}e $\widehat{\mu}_{i}(\overline{E})$, est le nombre r{\'e}el $$\widehat{\mu}_{i}(\overline{E}):=P_{\overline{E}}(i)-P_{\overline{E}}(i-1)\ .$$\end{defi}De la sorte, on a $$\widehat{\deg}_{\mathrm{n}}\overline{E}=\sum_{i=1}^{n}{\widehat{\mu}_{i}(\overline{E})}\ .$$Soit $\varepsilon>0$. Compte tenu de l'encadrement des normes de $\overline{E}$ en fonction de celles du fibr{\'e} ad{\'e}lique hermitien $\overline{E}_{\varepsilon}$ (voir~\eqref{eqref:important}), on a $$\forall\, x\in[0,n],\quad P_{\overline{E}_{\varepsilon}}(x)-\frac{x}{D}\left(\log\Delta(\overline{E})+D\log(1+\varepsilon)\right)\le P_{\overline{E}}(x)\le P_{\overline{E}_{\varepsilon}}(x),$$ce qui implique \begin{equation}\begin{split}\label{comparaisonpentesijohn}&-\frac{i}{D}\left(\log\Delta(\overline{E})+D\log(1+\varepsilon)\right)\\ &\qquad\le\widehat{\mu}_{i}(\overline{E})-\widehat{\mu}_{i}(\overline{E}_{\varepsilon})\le\frac{i-1}{D}\left(\log\Delta(\overline{E})+D\log(1+\varepsilon)\right)\ .\end{split}\end{equation}Cet encadrement se r{\'e}v{\`e}le particuli{\`e}rement pr{\'e}cieux pour g{\'e}n{\'e}raliser les propri{\'e}t{\'e}s des pentes d'un fibr{\'e} ad{\'e}lique hermitien {\`a} un fibr{\'e} vectoriel ad{\'e}lique quelconque.\par Par ailleurs, la concavit{\'e} de $P_{\overline{E}}$ entra{\^\i}ne la d{\'e}croissance de la suite $(\widehat{\mu}_{i}(\overline{E}))_{1\le i\le n}$.
\begin{lemma}\label{lemmamuunmumax} On a 
\begin{equation*}
\widehat{\mu}_{1}(\overline{E})=\widehat{\mu}_{\mathrm{max}}(\overline{E})\ .
\end{equation*}
\end{lemma}
\begin{proof}
Pour tout sous-fibr{\'e} $\overline{F}$ de $\overline{E}$ de dimension $m$, on a $$\widehat{\deg}_{\mathrm{n}}\overline{F}\le P_{\overline{E}}(m)=\sum_{i=1}^{m}{\widehat{\mu}_{i}(\overline{E})}\le m\widehat{\mu}_{1}(\overline{E})$$donc $\widehat{\mu}_{\mathrm{n}}(\overline{F})\le\widehat{\mu}_{1}(\overline{E})$, puis $\widehat{\mu}_{\mathrm{max}}(\overline{E})\le\widehat{\mu}_{1}(\overline{E})$. Inversement, si l'on choisit pour $\overline{F}$ le premier sommet non trivial de $P_{\overline{E}}$, on a $\widehat{\mu}_{1}(\overline{E})=\widehat{\mu}_{\mathrm{n}}(\overline{F})$ et donc $\widehat{\mu}_{1}(\overline{E})\le\widehat{\mu}_{\mathrm{max}}(\overline{E})$.
\end{proof}
 Je ne sais pas si l'\'egalit\'e $\widehat{\mu}_{n}(\overline{E})=\widehat{\mu}_{\mathrm{min}}(\overline{E})$ est vraie en g\'en\'eral. Cependant on a le r\'esultat suivant qui assure l'\'egalit\'e dans le cas d'un fibr\'e vectoriel hermitien.
\begin{lemma}\label{lemmeliti}
Soit $\overline{E}$ un fibr\'e vectoriel ad\'elique sur $\spec k$, de dimension $n\ge 1$. On a \begin{equation*}\left\vert \widehat{\mu}_{n}(\overline{E})-\widehat{\mu}_{\mathrm{min}}(\overline{E})\right\vert\le\frac{n}{D}\log\Delta(\overline{E})\ \cdotp\end{equation*}
\end{lemma}
Il faut prendre garde \`a ne pas confondre dans cette in\'egalit\'e la $n^{\text{\`eme}}$ pente $\widehat{\mu}_{n}(\overline{E})$ de $\overline{E}$ avec sa pente normalis\'ee $\widehat{\mu}_{\mathrm{n}}(\overline{E})$.
\begin{proof}Lorsque $\overline{E}$ est hermitien, nous montrerons un peu plus loin que $\widehat{\mu}_{n}(\overline{E})=-\widehat{\mu}_{1}(\overline{E^{\mathsf{v}}})$ (voir propri\'et\'e~\ref{proprietes511}, (2)). Dans ce cas, l'\'egalit\'e $\widehat{\mu}_{n}(\overline{E})=\widehat{\mu}_{\mathrm{min}}(\overline{E})$ r\'esulte des lemmes~\ref{lemmacomparaisonmumaxmumin} et~\ref{lemmamuunmumax}. Le cas d'un fibr\'e vectoriel ad\'elique quelconque s'en d\'eduit au moyen du fibr\'e hermitien $\overline{E}_{\varepsilon}$, de l'encadrement~\eqref{comparaisonpentesijohn} et de la remarque~\ref{remarquecinqcinq}.
\end{proof}  
\subsubsection*{Filtration canonique dans le cas hermitien} Lorsque $\overline{E}$ est un fibr{\'e} ad{\'e}lique \emph{hermitien}, les sous-fibr{\'e}s de $\overline{E}$ qui correspondent aux sommets du graphe de $P_{\overline{E}}$ sont uniques et ils forment une filtration de $\overline{E}$.
\begin{lemma}\label{lemmefiltrationcanonique} Soit $\overline{E}$ un fibr{\'e} ad{\'e}lique hermitien et $i\in\{1,\ldots,n-1\}$ un point de non-d\'erivabilit\'e de $P_{\overline{E}}$. Alors il existe un unique sous-espace vectoriel $E_{i}$ de $E$, de dimension $i$, tel que $\widehat{\deg}_{\mathrm{n}}\overline{E_{i}}=P_{\overline{E}}(i)$. 
\end{lemma}
\begin{proof}L'existence de $E_{i}$ est une cons\'equence du fait qu'il n'y a qu'un nombre fini de sous-espaces de $E$ dont le degr\'e ad\'elique est plus grand qu'une constante. Pour montrer l'unicit\'e, on raisonne par l'absurde et on suppose qu'il existe $E_{i}$ et $E_{i}'$ deux sous-espaces vectoriels distincts de $E$, de m{\^e}me dimension $i$, et de degr\'es ad\'eliques \'egaux \`a $P_{\overline{E}}(i)$. Soit $\delta_{1}:=\dim(E_{i}+E_{i}')$ et $\delta_{2}:=\dim(E_{i}\cap E_{i}')$. Par concavit{\'e} de $P_{\overline{E}}$ et comme $\delta_{1}\ne i$, on a $$\frac{P_{\overline{E}}(\delta_{1})-P_{\overline{E}}(i)}{\delta_{1}-i}<\frac{P_{\overline{E}}(i)-P_{\overline{E}}(\delta_{2})}{i-\delta_{2}}$$et donc, compte tenu de la relation $\delta_{1}-i=i-\delta_{2}$, on a \begin{equation}\label{inegaliteinterinter}\widehat{\deg}_{\mathrm{n}}\overline{E_{i}+E_{i}'}+\widehat{\deg}_{\mathrm{n}}\overline{E_{i}\cap E_{i}'}\le P_{\overline{E}}(\delta_{1})+P_{\overline{E}}(\delta_{2})<2P_{\overline{E}}(i)\ .\end{equation}D'apr{\`e}s la proposition~\ref{prop:sommededegrefibres} et comme $\overline{E}$ est hermitien, le membre de gauche de~\eqref{inegaliteinterinter} est minor{\'e} par la somme $\widehat{\deg}_{\mathrm{n}}\overline{E_{i}}+\widehat{\deg}_{\mathrm{n}}\overline{E_{i}'}=2P_{\overline{E}}(i)$ et ceci contredit~\eqref{inegaliteinterinter}. 
\end{proof}
\begin{prop}\label{proposition:filtrationcanonique}Soit $\overline{E}$ un fibr{\'e} ad{\'e}lique hermitien sur $\spec k$. Il existe une \emph{unique} filtration $\{0\}=E_{0}\subsetneq E_{1}\subsetneq\cdots\subsetneq E_{g}:=E$ telle que, pour tout $j\in\{0,\ldots,g\}$, on ait $P_{\overline{E}}(\dim E_{j})=\widehat{\deg}_{\mathrm{n}}\overline{E_{j}}$ et $\{\dim E_{j}\,;\ 1\le j\le g-1\}$ est l'ensemble des points de non-d\'erivabilit\'e de la fonction $P_{\overline{E}}$.
\end{prop}
Cette filtration est appel{\'e}e \emph{filtration canonique}, ou \emph{filtration d'Harder-Narasimhan-Stuhler-Grayson}.
\begin{proof}
Le lemme~\ref{lemmefiltrationcanonique} montre l'existence et l'unicit{\'e} des $\overline{E_{j}}$. Il ne reste plus qu'{\`a} prouver l'inclusion $E_{j}\subseteq E_{j+1}$. Elle repose exactement sur le m{\^e}me argument que celui employ{\'e} dans la preuve du lemme~\ref{lemmefiltrationcanonique}, appliqu{\'e} aux sous-espaces $E_{j}$, $E_{j+1}$, $E_{j}\cap E_{j+1}$ et $E_{j}+E_{j+1}$. 
\end{proof}
La d{\'e}monstration donn{\'e}e ici de l'existence et de l'unicit{\'e} de la filtration canonique dans le cas hermitien s'adapte mal au cas g{\'e}n{\'e}ral en partie {\`a} cause de la minoration du degr{\'e} d'une somme de sous-espaces de $E$, qui fait intervenir la distance ad{\'e}lique \`a l'espace hermitien $(k^{n},\vert.\vert_{2})$. Cela emp{\^e}che les arguments de concavit{\'e} de $P_{\overline{E}}$ de fonctionner.
\begin{proprietes}\label{proprietes511}\noignorespaces Soit $\overline{L}$ un fibr{\'e} en droites ad{\'e}lique et $\overline{E}$ un fibr{\'e} vectoriel ad{\'e}lique.
\begin{enumerate}
\item[1)] Pour toute norme tensorielle ad\'elique $\alpha$ d'ordre $2$ et pour tout $x\in[0,n]$, on a \begin{equation*}P_{\overline{E}\otimes_{\alpha}\overline{L}}(x)=P_{\overline{E}}(x)+x\widehat{\deg}_{\mathrm{n}}\overline{L}\ .\end{equation*}En particulier, pour tout $i\in\{1,\ldots,n\}$, on a \begin{equation*}\widehat{\mu}_{i}(\overline{E}\otimes_{\alpha}\overline{L})=\widehat{\mu}_{i}(\overline{E})+\widehat{\deg}_{\mathrm{n}}\overline{L}\ .\end{equation*}
\item[2)] Pour tout $i\in\{1,\ldots,n\}$, on a$$\left\vert\widehat{\mu}_{i}(\overline{E^{\mathsf{v}}})+\widehat{\mu}_{n-i+1}(\overline{E})\right\vert\le\widehat{\deg}_{\mathrm{n}}L(\overline{E})-\widehat{\deg}_{\mathrm{n}}J(\overline{E})+\frac{n}{D}\log\Delta(\overline{E})\ .$$
\item[3)] Pour tout $i\in\{1,\ldots,n\}$, on a \begin{equation}\label{ineqminmax}\left\vert\widehat{\mu}_{i}(\overline{E})-\underset{E_{2}}{\mathrm{min}}\,\underset{E_{1}}{\mathrm{max}}\left\{\widehat{\mu}_{\mathrm{n}}\left(\overline{E_{1}/E_{2}}\right)\right\}\right\vert\le\frac{i}{D}\log\Delta(\overline{E})\end{equation}o\`u $E_{1}$ et $E_{2}$ varient parmi les sous-espaces vectoriels de $E$ v{\'e}rifiant $E_{2}\subseteq E_{1}$, $\dim E_{1}\ge i$ et $\dim E_{2}\le i-1$. On a aussi \begin{equation}\label{ineqmaxmin}\left\vert\widehat{\mu}_{i}(\overline{E})-\underset{E_{1}}{\mathrm{max}}\,\underset{E_{2}}{\mathrm{min}}\left\{\widehat{\mu}_{\mathrm{n}}\left(\overline{E_{1}/E_{2}}\right)\right\}\right\vert\le\frac{i}{D}\log\Delta(\overline{E})\ .\end{equation} 
\end{enumerate}
\end{proprietes}
\begin{rema}
La diff\'erence $\widehat{\deg}_{\mathrm{n}}L(\overline{E})-\widehat{\deg}_{\mathrm{n}}J(\overline{E})$ qui est dans l'estimation 2) est plus petite que $\frac{n}{D}\log\Delta(\overline{E})$ (voir~\eqref{eqref:estimationdifference}), terme lui-m{\^e}me inf{\'e}rieur {\`a} $n\log(2n)$ (voir~\eqref{eq:encadrementquotientadelique}).
\end{rema}
\begin{proof}
1) Cela r{\'e}sulte directement de la formule $\widehat{\deg}_{\mathrm{n}}(\overline{F}\otimes_{\alpha}\overline{L})=\widehat{\deg}_{\mathrm{n}}\overline{F}+(\dim F)\widehat{\deg}_{\mathrm{n}}\overline{L}$ qui d\'ecoule de la proposition~\ref{penteproduittensorielfibre52}.\par\noindent 2) Pour tout sous-fibr{\'e} ad{\'e}lique $\iota:\overline{F}\hookrightarrow\overline{E}$ de dimension $m$, on a une suite exacte $0\to\ker\iota^{\mathsf{v}}\to E^{\mathsf{v}}\to F^{\mathsf{v}}\to 0$ et donc (preuve de la proposition~\ref{propositiondegrequotient46})
$$\left\vert\widehat{\deg}_{\mathrm{n}}\overline{E^{\mathsf{v}}}-\widehat{\deg}_{\mathrm{n}}\overline{F^{\mathsf{v}}}-\widehat{\deg}_{\mathrm{n}}\overline{\ker\iota^{\mathsf{v}}}\right\vert\le\widehat{\deg}_{\mathrm{n}}L(\overline{E^{\mathsf{v}}})-\widehat{\deg}_{\mathrm{n}}J(\overline{E^{\mathsf{v}}}),$$ce qui entra{\^\i}ne $$\widehat{\deg}_{\mathrm{n}}\overline{F}\le\widehat{\deg}_{\mathrm{n}}\overline{\ker\iota^{\mathsf{v}}}+\widehat{\deg}_{\mathrm{n}}L(\overline{E^{\mathsf{v}}})-\widehat{\deg}_{\mathrm{n}}J(\overline{E^{\mathsf{v}}})-\widehat{\deg}_{\mathrm{n}}\overline{E^{\mathsf{v}}}$$gr{\^a}ce {\`a} la proposition~\ref{propositiondualitedegres}. Le degr{\'e} ad{\'e}lique normalis{\'e} de $\overline{\ker\iota^{\mathsf{v}}}$ est major{\'e} par $P_{\overline{E^{\mathsf{v}}}}(n-m)$. De plus, la dualit{\'e} $J(\overline{E^{\mathsf{v}}})=L(\overline{E})^{\mathsf{v}}$ entre les fibr{\'e}s hermitiens de John et L{\"o}wner donne \begin{equation*}\widehat{\deg}_{\mathrm{n}}L(\overline{E^{\mathsf{v}}})-\widehat{\deg}_{\mathrm{n}}J(\overline{E^{\mathsf{v}}})=\widehat{\deg}_{\mathrm{n}}L(\overline{E})-\widehat{\deg}_{\mathrm{n}}J(\overline{E})\ .\end{equation*}On en d{\'e}duit \begin{equation*}P_{\overline{E}}(m)\le P_{\overline{E^{\mathsf{v}}}}(n-m)+\widehat{\deg}_{\mathrm{n}}L(\overline{E})-\widehat{\deg}_{\mathrm{n}}J(\overline{E})-\widehat{\deg}_{\mathrm{n}}\overline{E^{\mathsf{v}}}\ .\end{equation*}En {\'e}changeant les r{\^o}les de  $\overline{E}$ et de son dual, on a aussi la majoration \begin{equation*}P_{\overline{E^{\mathsf{v}}}}(m)\le P_{\overline{E}}(n-m)+\widehat{\deg}_{\mathrm{n}}L(\overline{E})-\widehat{\deg}_{\mathrm{n}}J(\overline{E})-\widehat{\deg}_{\mathrm{n}}\overline{E}\ .\end{equation*}Par d{\'e}finition des pentes $\widehat{\mu}_{i}$, on a alors \begin{equation*}\left\vert\widehat{\mu}_{i}(\overline{E^{\mathsf{v}}})+\widehat{\mu}_{n-i+1}(\overline{E})\right\vert\le 2\left(\widehat{\deg}_{\mathrm{n}}L(\overline{E})-\widehat{\deg}_{\mathrm{n}}J(\overline{E})\right)-\left(\widehat{\deg}_{\mathrm{n}}\overline{E}+\widehat{\deg}_{\mathrm{n}}\overline{E^{\mathsf{v}}}\right)\ \cdotp\end{equation*}Au moyen de la proposition~\ref{prop:comparaisonjohnadelique} et de l'in\'egalit\'e $\widetilde{\vr}(\overline{E})\widetilde{\vr}(\overline{E^{\mathsf{v}}})\le\Delta(\overline{E})$, analogue \`a~\eqref{eq:encadrementquotientadelique}, l'on d\'eduit
 \begin{equation*}\begin{split}&\widehat{\deg}_{\mathrm{n}}L(\overline{E})-\widehat{\deg}_{\mathrm{n}}J(\overline{E})-\left(\widehat{\deg}_{\mathrm{n}}\overline{E}+\widehat{\deg}_{\mathrm{n}}\overline{E^{\mathsf{v}}}\right)\\ &\quad = \widehat{\deg}_{\mathrm{n}}L(\overline{E})+\widehat{\deg}_{\mathrm{n}}L(\overline{E^{\mathsf{v}}})-\left(\widehat{\deg}_{\mathrm{n}}\overline{E}+\widehat{\deg}_{\mathrm{n}}\overline{E^{\mathsf{v}}}\right)\le\frac{n}{D}\log\Delta(\overline{E})\ .\end{split}\end{equation*}On notera au passage l'{\'e}galit{\'e} $P_{\overline{E}}(m)=P_{\overline{E^{\mathsf{v}}}}(n-m)+\widehat{\deg}_{\mathrm{n}}\overline{E}$ obtenue lorsque $\overline{E}$ est hermitien.\par\noindent 3) On commence par {\'e}tablir l'{\'e}galit{\'e} \begin{equation}\label{egaliteeqcashermitien}\widehat{\mu}_{i}(\overline{E})=\underset{E_{2}}{\mathrm{min}}\,\underset{E_{1}}{\mathrm{max}}\left\{\widehat{\mu}_{\mathrm{n}}\left(\overline{E_{1}/E_{2}}\right)\right\}\end{equation}lorsque $\overline{E}$ est un fibr{\'e} ad{\'e}lique hermitien. On note $\alpha$ le minimax du membre de droite. Soit $\overline{E_{2}}\subseteq\overline{E_{1}}$ de dimensions respectives $n_{2}$ et $n_{1}$. On suppose $n_{2}\le i-1<n_{1}$. D'apr\`es la proposition~\ref{propositiondegrequotient46}, on a $$\widehat{\mu}_{\mathrm{n}}\left(\overline{E_{1}/E_{2}}\right)\le\frac{P_{\overline{E}}(n_{1})-\widehat{\deg}_{\mathrm{n}}\overline{E_{2}}}{n_{1}-n_{2}}$$et donc $\widehat{\deg}_{\mathrm{n}}\overline{E_{2}}\le P_{\overline{E}}(n_{1})-(n_{1}-n_{2})\alpha$ puis $P_{\overline{E}}(n_{2})\le P_{\overline{E}}(n_{1})-(n_{1}-n_{2})\alpha$. Inversement, on a $$\frac{\widehat{\deg}_{\mathrm{n}}\overline{E_{1}}-P_{\overline{E}}(n_{2})}{n_{1}-n_{2}}\le\widehat{\mu}_{\mathrm{n}}\left(\overline{E_{1}/E_{2}}\right)\le\underset{E_{1}}{\mathrm{max}}{\left\{\widehat{\mu}_{\mathrm{n}}\left(\overline{E_{1}/E_{2}}\right)\right\}}$$et donc $P_{\overline{E}}(n_{1})\le P_{\overline{E}}(n_{2})+(n_{1}-n_{2})\alpha$.  Ainsi, pour tout $n_{2}<i\le n_{1}$, on a \begin{equation*}P_{\overline{E}}(n_{1})=P_{\overline{E}}(n_{2})+(n_{1}-n_{2})\alpha\end{equation*}et donc, par d{\'e}finition, $\alpha=\widehat{\mu}_{i}(\overline{E})$. Le passage au cas g{\'e}n{\'e}ral s'effectue au moyen du fibr{\'e} $\overline{E}_{\varepsilon}$ associ{\'e} {\`a} $\overline{E}$ et pour un nombre r{\'e}el $\varepsilon>0$ donn\'e (voir~\eqref{eqref:important}). On observe que \begin{equation*}\begin{split}&-\frac{1}{D}\left(\log\Delta(\overline{E})+D\log(1+\varepsilon)\right)+\widehat{\mu}_{\mathrm{n}}\left(E_{1}/E_{2},\vert.\vert_{\varepsilon,\overline{E}}\right)\\ &\qquad\le\widehat{\mu}_{\mathrm{n}}\left(\overline{E_{1}/E_{2}}\right)\le\widehat{\mu}_{\mathrm{n}}\left(E_{1}/E_{2},\vert.\vert_{\varepsilon,\overline{E}}\right)\end{split}\end{equation*}puis, avec la formule~\eqref{egaliteeqcashermitien}, \begin{equation*}-\frac{1}{D}\log\Delta(\overline{E})-\log(1+\varepsilon)+\widehat{\mu}_{i}(\overline{E}_{\varepsilon})\le\underset{E_{2}}{\mathrm{min}}\,\underset{E_{1}}{\mathrm{max}}\left\{\widehat{\mu}_{\mathrm{n}}\left(\overline{E_{1}/E_{2}}\right)\right\}\le\widehat{\mu}_{i}(\overline{E}_{\varepsilon})\ \cdotp\end{equation*}Ce qui entra{\^\i}ne l'encadrement annonc{\'e} de $\widehat{\mu}_{i}(\overline{E})$ \textit{via} l'estimation~\eqref{comparaisonpentesijohn} et en faisant tendre $\varepsilon$ vers $0$. L'in\'egalit\'e~\eqref{ineqmaxmin} s'obtient de la m{\^e}me mani\`ere.
\end{proof}
\subsection{Fibr{\'e}s ad{\'e}liques semi-stables}
\begin{defi}
Un fibr{\'e} vectoriel ad{\'e}lique $\overline{E}$ est dit \emph{semi-stable} si toutes les pentes $\widehat{\mu}_{i}(\overline{E})$, $1\le i\le n$, sont {\'e}gales.
\end{defi}
De mani{\`e}re {\'e}quivalente, le fibr{\'e} $\overline{E}$ est semi-stable si $\widehat{\mu}_{1}(\overline{E})$ est la pente normalis{\'e}e $\widehat{\mu}_{\mathrm{n}}(\overline{E})$ de $\overline{E}$, ou si, pour tout $m\in [0,n]$, $P_{\overline{E}}(m)=m\widehat{\mu}_{\mathrm{n}}(\overline{E})$ ou bien encore si $\widehat{\mu}_{\mathrm{max}}(\overline{E})=\widehat{\mu}_{\mathrm{n}}(\overline{E})$ (car la somme des $\widehat{\mu}_{i}(\overline{E})$ {\'e}gale $\widehat{\deg}_{\mathrm{n}}\overline{E}=n\widehat{\mu}_{\mathrm{n}}(\overline{E})$). En vertu de la premi\`ere des propri\'et\'es~\ref{proprietes511}, le fibr\'e $\overline{E}$ est semi-stable si et seulement si, pour tout fibr\'e en droites ad\'elique $\overline{L}$ sur $\spec k$, le produit tensoriel $\overline{E}\otimes\overline{L}$ est semi-stable. Par ailleurs, comme l'a mentionn{\'e} J.-B.~Bost dans~\cite{bost2}, l'on dispose de conditions suffisantes pour v{\'e}rifier la semi-stabilit{\'e} de $\overline{E}$.
\begin{prop}
Supposons qu'il existe un sous-groupe $G$ de $\GL(E)$ ayant les propri{\'e}t{\'e}s suivantes~:
\begin{enumerate}\item[(i)] tout {\'e}l{\'e}ment $\varphi\in G$ induit une isom{\'e}trie de $E_{v}$ en toutes les places $v$ de $k$, 
\item[(ii)] pour tout sous-espace vectoriel $E'$ de $E$, si l'orbite $\{\varphi(E')\,;\ \varphi\in G\}$ est finie alors $E'=0$ ou $E'=E$.
\end{enumerate}
Alors $\overline{E}$ est semi-stable. Lorsque $\overline{E}$ est hermitien, cette m{\^e}me conclusion reste valide en rempla{\c c}ant l'hypoth\`ese (ii) par la condition plus faible \begin{enumerate}\item[(ii)'] si, pour tout $\varphi\in G$, on a $\varphi(E')=E'$ alors $E'=0$ ou $E'=E$ (l'action de $G$ sur $E$ est \emph{irr\'eductible}).
\end{enumerate}
\end{prop} 
\begin{proof}
La condition (i) implique que, pour tout $\overline{E'}\subseteq\overline{E}$, on a $\widehat{\deg}\overline{E'}=\widehat{\deg}\overline{\varphi(E')}$. En particulier, si $\overline{E'}$ est un sommet du polygone canonique construit avec $P_{\overline{E}}$, il en est de m{\^e}me pour $\overline{\varphi(E')}$ et l'orbite de $E'$ sous l'action de $G$ est finie. L'hypoth{\`e}se (ii) entra{\^\i}ne alors $E'=\{0\}$ ou $\overline{E'}=\overline{E}$, ce qui signifie que $\overline{E}$ est semi-stable. Dans le cas hermitien, la proposition~\ref{proposition:filtrationcanonique} donne l'unicit{\'e} des sous-fibr{\'e}s vectoriels ad{\'e}liques de $\overline{E}$ qui forment les sommets du graphe de $P_{\overline{E}}$. On a donc $\varphi(E')=E'$ et l'hypoth{\`e}se (ii)' suffit alors pour conclure \`a la semi-stabilit{\'e} de $\overline{E}$.    
\end{proof}
La g\'eom\'etrie d'Arakelov fournit une famille d'exemples de fibr\'es vectoriels hermitiens semi-stables, lorsque $k$ est un corps de nombres. En effet, soit $X$ une vari\'et\'e ab\'elienne d\'efinie sur $k$, munie d'un fibr\'e en droites $L$ ample et sym\'etrique. Soit $(\mathcal{X},\overline{\mathcal{L}})$ un mod\`ele de Moret-Bailly de $(X,L)$, au sens du \S~$4.3.1$ de~\cite{bost3}. Sans entrer dans les d\'etails, rappelons que $\mathcal{X}$ est un sch\'ema en groupes semi-stable sur l'anneau des entiers d'une extension finie $K$ de $k$, de fibre g\'en\'erique $A\otimes K$ et $\overline{\mathcal{L}}$ est un fibr\'e hermitien \emph{cubiste}, de fibre g\'en\'erique $L_{K}$. L'espace des section globales $\mathrm{H}^{0}(\mathcal{X},\mathcal{L})$ est un $\mathcal{O}_{K}$-module projectif de type fini. Les m{\'e}triques sur $\overline{\mathcal{L}}$ le munissent de m\'etriques hermitiennes aux places archim{\'e}diennes de $K$, ce qui lui conf{\`e}re alors une structure de fibr\'e ad\'elique hermitien $\overline{\mathrm{H}^{0}(X,L)}$ sur $K$. Au moyen du groupe de Mumford qui agit irr\'eductiblement sur $\mathrm{H}^{0}(X,L)$, J.-B.~Bost a montr\'e que $\overline{\mathrm{H}^{0}(X,L)}$ est \emph{semi-stable} (et on l'on conna{\^\i}t m{\^e}me une expression de sa pente normalis\'ee, voir le th\'eor\`eme~$4.2$ de~\cite{bost2}, ainsi que le th\'eor\`eme~$2.10$, ii), de~\cite{graftieaux1} et aussi le dernier paragraphe).

\subsection{Minima successifs ad{\'e}liques}
Comme le montre le travail de T.~Borek~\cite{ThomasBorek} dans le cas des fibr{\'e}s vectoriels hermitiens sur $\spec\mathcal{O}_{k}$, les pentes de $\overline{E}$ telles que nous venons de les d{\'e}finir sont {\'e}troitement li{\'e}es aux minima successifs ad{\'e}liques de $\overline{E}$. Si $a=(a_{v})_{v}$ est un {\'e}l{\'e}ment de $k_{\mathbf{A}}$, on note $\mathbb{B}(\overline{E},a)$ l'ensemble $$\mathbb{B}(\overline{E},a)=\{(x_{v})_{v}\in E_{\mathbf{A}}\,;\ \Vert x_{v}\Vert_{E,v}\le\vert a_{v}\vert_{v}\ \text{pour toute place $v$ de $k$}\}$$(ainsi $\mathbb{B}(\overline{E},1)=\mathbb{B}(\overline{E})$).

\begin{defi}
Soit $\overline{E}$ un fibr{\'e} vectoriel ad{\'e}lique sur $\spec k$, de dimension $n$. Pour $i\in\{1,\ldots,n\}$, le $i^{\text{{\`e}me}}$ minimum relatif {\`a} $\overline{E}$, que l'on notera $\lambda_{i}(\overline{E})$, est la borne inf{\'e}rieure de l'ensemble des nombres r{\'e}els positifs de la forme $\vert a\vert_{\mathbf{A}}$ o\`u $a\in k_{\mathbf{A}}$ est tel que $E\cap\mathbb{B}(\overline{E},a)$ contienne $i$ vecteurs $k$-lin{\'e}airement ind{\'e}pendants.
\end{defi}
Ces nombres $\lambda_{i}(\overline{E})$ d{\'e}pendent \textit{a priori} du corps $k$. Si $k$ est un corps de nombres, il arrive souvent que l'on se restreigne {\`a} des {\'e}l{\'e}ments $a$ de la forme $(\lambda,\ldots,\lambda,1,\ldots)$ o\`u $\lambda\in\mathbf{R}^{+}$ aux places archim{\'e}diennes de $k$, comme le font par exemple E.~Bombieri \& J.~Vaaler~\cite{BombieriVaaler}. Les puissances $D^{\text{\`emes}}$ de leurs minima sont donc plus grands que ceux consid{\'e}r{\'e}s ici. La d{\'e}finition donn{\'e}e ici est emprunt\'ee \`a l'article de J.L.~Thunder~\cite{jlthunder}. On notera que $$\lambda_{1}(\overline{E})\le\cdots\le\lambda_{n}(\overline{E})$$et que, si $\overline{E'}$ est un sous-fibr{\'e} ad{\'e}lique de $\overline{E}$ alors $\lambda_{i}(\overline{E})\le\lambda_{i}(\overline{E'})$ pour tout $i\in\{1,\ldots,\dim E'\}$. De plus, on a \begin{equation}\label{encadrementiememinima}\lambda_{i}(\overline{E}_{\varepsilon})\le\lambda_{i}(\overline{E})\le\Delta(\overline{E})(1+\varepsilon)^{D}\lambda_{i}(\overline{E}_{\varepsilon})\ .\end{equation}Existe alors une variante du \emph{second th{\'e}or{\`e}me de Minkowski}.
\begin{theo}\label{secondtheorminkowski}
Le produit des minima successifs de $\overline{E}$ v{\'e}rifie\begin{equation*}\lambda_{1}(\overline{E})\cdots\lambda_{n}(\overline{E})\le\kappa^{n}\frac{\covol(E)}{\vol(\mathbb{B}(\overline{E}))}\vr(\overline{E})^{n}\end{equation*}o\`u $\kappa=2^{[k:\mathbf{Q}]}$ si $k$ est un corps de nombres et $\kappa=q$ si $k$ est un corps de fonctions. Ici, $\vol$ d{\'e}signe une mesure de Haar quelconque sur $E_{\mathbf{A}}$ et $\covol(E)$ est la mesure de l'espace quotient $E_{\mathbf{A}}/E$, relative \`a $\vol$.
\end{theo} 
Dans le cas d'un corps de nombres, le th{\'e}or{\`e}me est {\'e}tabli dans~\cite{BombieriVaaler} (modulo la remarque ci-dessus sur les d{\'e}finitions l{\'e}g{\`e}rement diff{\'e}rentes des minima, diff{\'e}rence qui joue en notre faveur). Pour un corps de fonctions, on consultera~\cite{jlthunder}, corollaire~1. Ces deux r\'ef\'erences traitent seulement le cas d'un fibr\'e vectoriel hermitien (c.-\`a-d. $\vr(\overline{E})=1$). Le cas g\'en\'eral s'en d\'eduit au moyen du fibr\'e de John $J(\overline{E})$ en observant que $\lambda_{i}(\overline{E})\le\lambda_{i}(J(\overline{E}))$ et en utilisant la d\'efinition de $\vr(\overline{E})$. \par Il est possible d'obtenir {\'e}galement une minoration du produit de ces minima en observant que, pour tout $\varepsilon>0$, il existe des {\'e}l{\'e}ments $a_{1}(\varepsilon),\ldots,a_{n}(\varepsilon)\in k_{\mathbf{A}}$ et une base $(e_{1}(\varepsilon),\ldots,e_{n}(\varepsilon))$ de $E$ tels que, pour tout $i\in\{1,\ldots,n\}$, on a $\vert a_{i}(\varepsilon)\vert_{\mathbf{A}}\in[\lambda_{i}(\overline{E}),\lambda_{i}(\overline{E})+\varepsilon]$ et, pour toute place $v$ de $k$, $\Vert e_{i}(\varepsilon)\Vert_{v}\le\vert a_{i}(\varepsilon)\vert_{v}$. De ce fait, le volume de la boule unit{\'e} de $\overline{E}$ --- volume relatif \`a la mesure de Haar $\mu_{E_{\mathbf{A}}}$ d\'efinie au paragraphe~\ref{paragrapheadelescorps} --- est minor{\'e} \begin{equation}\label{minorationintermediaireun}\vol(\mathbb{B}(\overline{E}))\ge\vol_{n}(b^{1}_{n,\mathbf{R}})^{r_{1}}(2^{n}\vol_{2n}(b^{1}_{n,\mathbf{C}}))^{r_{2}}\left(\prod_{v,i}{\Vert e_{i}(\varepsilon)\Vert_{v}^{n_{v}}}\right)^{-1}\end{equation}car, pour tout $x=\sum_{i=1}^{n}{x_{i}e_{i}(\varepsilon)}\in E\otimes k_{v}$, on a \begin{equation*}\Vert x\Vert_{E,v}\le\begin{cases}\underset{1\le i\le n}{\mathrm{max}}\left\{\vert x_{i}\vert_{v}\Vert e_{i}(\varepsilon)\Vert_{v}\right\} & \text{si $v$ est ultram\'etrique}\\ \sum_{i=1}^{n}{\vert x_{i}\vert_{v}\Vert e_{i}(\varepsilon)\Vert_{v}} & \text{si $v$ est archim\'edienne}\end{cases}\end{equation*}et la formule~\eqref{transfoechelle} permet de calculer le degr{\'e} de $E$ muni des m{\'e}triques {\`a} droite. L'estimation~\eqref{minorationintermediaireun} entra{\^\i}ne (avec $\varepsilon\to 0$) \begin{equation}\lambda_{1}(\overline{E})\cdots\lambda_{n}(\overline{E})\frac{\vol(\mathbb{B}(\overline{E}))}{\covol(E)}\ge\begin{cases}\frac{2^{nD}\pi^{nr_{2}}}{(n!)^{r_{1}}(2n!)^{r_{2}}\vert D_{k}\vert^{n/2}} & \text{{\footnotesize si $k$ est un corps de nombres}},\\ q^{-n(g(k)-1)} & \text{{\footnotesize si $k$ est un corps de fonctions}}.\end{cases}\end{equation}Les quantit{\'e}s qui sont {\`a} droite proviennent des formules~\eqref{formuledevolumep} et~\eqref{formulevolumecomplexep} ainsi que du choix de la mesure $\mu$ sur $k_{\mathbf{A}}$ d{\'e}finie au paragraphe~\ref{paragrapheadelescorps}. Dans le cas d'un corps de nombres, cette minoration est la m{\^e}me que celle de E.~Bombieri \& J.~Vaaler~\cite{BombieriVaaler}, th\'eor\`eme~$6$.\par
Comme nous l'avons dit, le lien entre ces minima et les pentes de $\overline{E}$ a {\'e}t{\'e} pr{\'e}cis{\'e} par T.~Borek~\cite{ThomasBorek}. L'{\'e}nonc{\'e} suivant est une g{\'e}n{\'e}ralisation aux fibr{\'e}s vectoriels ad{\'e}liques. 
\begin{theo}Soit $\overline{E}$ un fibr{\'e} vectoriel ad{\'e}lique sur $\spec k$, de dimension $n$. Posons $$C(n,k):=\frac{1}{D}\log\frac{\kappa^{n}\covol(k^{n})}{\vol\mathbb{B}(k^{n},\vert.\vert_{2})}$$o\`u $\kappa$ est la constante qui intervient dans le th{\'e}or{\`e}me~\ref{secondtheorminkowski}. Alors, pour tout $i\in\{1,\ldots,n\}$, on a \begin{equation*}-\frac{i}{D}\log\Delta(\overline{E})\le\widehat{\mu}_{i}(\overline{E})+\frac{1}{D}\log\lambda_{i}(\overline{E})\le\frac{i}{n}C(n,k)+\frac{i}{D}\log\Delta(\overline{E})\ \cdotp\end{equation*} 
\end{theo}
Nous n'allons pas refaire la preuve de cet {\'e}nonc{\'e} qui --- dans le cas hermitien --- est l'objet du travail~\cite{ThomasBorek} de T.~Borek. M{\^e}me si ce dernier ne consid{\`e}re que des fibr{\'e}s vectoriels hermitiens sur $\spec\mathcal{O}_{k}$, les m{\'e}thodes de cet article se g{\'e}n{\'e}ralisent au cadre ad{\'e}lique car, comme nous l'avons vu, les degr{\'e}s et les pentes d{\'e}finis ici donnent les m{\^e}mes invariants que ceux d{\'e}finis habituellement pour $\overline{E}\to\spec\mathcal{O}_{k}$. De plus, l'on dispose {\'e}galement de la filtration canonique associ{\'e}e {\`a} un fibr{\'e} ad{\'e}lique hermitien, qui est l'outil essentiel dans~\cite{ThomasBorek}. Le cas g{\'e}n{\'e}ral d'un fibr{\'e} vectoriel ad{\'e}lique quelconque se traite {\`a} partir du cas hermitien au moyen du fibr{\'e} hermitien $\overline{E}_{\varepsilon}$ et des encadrements~\eqref{comparaisonpentesijohn} et~\eqref{encadrementiememinima}. On peut mentionner la formule asymptotique $C(n,k)=\frac{n}{2}\log n+\mathrm{O}(n)$ lorsque $n\to+\infty$ (le corps $k$ \'etant fix\'e), formule qui s'obtient au moyen de la mesure $\mu$ d\'efinie au \S~\ref{paragrapheadelescorps} et des formules~\eqref{formuledevolumep} et~\eqref{formulevolumecomplexep}.  

\section{In{\'e}galit{\'e}s de pentes}
L'objectif de ce paragraphe est de comparer les degr{\'e}s et les pentes de fibr{\'e}s vectoriels ad{\'e}liques reli{\'e}s entre eux par une application lin{\'e}aire. En r{\'e}alit{\'e} nous avons d{\'e}j{\`a} compar{\'e} de telles donn{\'e}es lorsque l'application sous-jacente est une inclusion. 
\begin{defi}\label{definitiondehauteur}
Soit $\overline{E},\overline{F}$ des fibr{\'e}s vectoriels ad{\'e}liques sur $\spec k$ et $\varphi:E\to F$ une application $k$-lin{\'e}aire. La \emph{hauteur} de $\varphi$ relative {\`a} $\overline{E}$ et $\overline{F}$, not\'ee $h(\overline{E},\overline{F};\varphi)$ ou, plus simplement, $h(\varphi)$ s'il n'y a pas d'ambigu{\"\i}t{\'e}, est la somme 
$$h(\overline{E},\overline{F};\varphi):=\frac{1}{D}\sum_{v}{n_{v}\log\Vert\varphi\Vert_{v}}\ \cdotp$$La norme $\Vert\varphi\Vert_{v}$ est la norme d'op{\'e}rateur de l'application induite $\varphi_{v}:E\otimes_{k}k_{v}\to F\otimes_{k}k_{v}$ et la somme ci-dessus ne comporte qu'un nombre fini de termes non nuls.\end{defi}La formule du produit assure l'invariance par multiplication par un scalaire non nul~:$$\forall\lambda\in k\setminus\{0\},\quad h(\overline{E},\overline{F};\lambda\varphi)=h(\overline{E},\overline{F};\varphi)\ \cdotp$$La hauteur est invariante par extension finie $K$ du corps de base $k$~:
\begin{equation*} h(\overline{E_{K}},\overline{F_{K}};\varphi_{K})=h(\overline{E},\overline{F};\varphi)\end{equation*}car la norme d'op\'erateur de $\varphi$ en une place $w$ de $K$ est \'egale \`a celle en la place $v$ de $k$ correspondante et car $\sum_{w\mid v}{[K_{w}:k_{v}]}=[K:k]$. De plus, il est utile de noter les deux encadrements suivants~: \begin{equation}\label{evaluationunhauteur}-\log(1+\varepsilon)-\frac{1}{D}\log\Delta(\overline{E})+h(\overline{E},\overline{F};\varphi)\le h(\overline{E},\overline{F}_{\varepsilon};\varphi)\le h(\overline{E},\overline{F};\varphi)\end{equation}et\begin{equation}\label{evaluationdeuxhauteur}h(\overline{E},\overline{F};\varphi)\le h(\overline{E}_{\varepsilon},\overline{F};\varphi)\le h(\overline{E},\overline{F};\varphi)+\frac{1}{D}\log\Delta(\overline{E})+\log(1+\varepsilon)\ \cdotp\end{equation}
\begin{lemma}\label{lemme:egalitedesdegres}
Si $\varphi:E\to F$ est un isomorphisme entre deux \emph{droites} vectorielles alors $$\widehat{\deg}_{\mathrm{n}}\overline{E}=\widehat{\deg}_{\mathrm{n}}\overline{F}+h(\overline{E},\overline{F};\varphi)\ .$$
\end{lemma}
La raison en est que, pour tout $x\in E_{v}\setminus\{0\}$, on a $\Vert\varphi\Vert_{v}=\Vert\varphi(x)\Vert_{F,v}/\Vert x\Vert_{E,v}$. On conclut au moyen du lemme~\ref{lemme4225} par exemple. En particulier, si $\overline{E}$ et $\overline{F}$ sont des fibr{\'e}s ad{\'e}liques \emph{hermitiens} (de dimension quelconque), et si $\varphi:E\to F$ est un isomorphisme, le lemme~\ref{lemme4226} entra{\^\i}ne \begin{equation}\label{casdegaliteimportant}\widehat{\deg}_{\mathrm{n}}\overline{E}=\widehat{\deg}_{\mathrm{n}}\overline{F}+h(\det\overline{E},\det\overline{F};\det\varphi)\ \cdotp\end{equation}
\begin{lemma}
Si $\varphi:E\to F$ est un isomorphisme alors les pentes normalis{\'e}es de $\overline{E}$ et $\overline{F}$ v{\'e}rifient\begin{equation}\widehat{\mu}_{\mathrm{n}}(\overline{E})\le\widehat{\mu}_{\mathrm{n}}(\overline{F})+h(\overline{E},\overline{F};\varphi)\ \cdotp\end{equation}
\end{lemma}
\begin{proof}
On peut supposer $E=F=k^{n}$ et $\varphi$ s'identifie {\`a} une matrice de $\GL_{n}(k)$. En toute place $v$ de $k$, on a l'inclusion $$\left\{x\in k_{v}^{n}\,;\ \Vert x\Vert_{E,v}\le\frac{1}{\Vert\varphi\Vert_{v}}\right\}\subseteq\left\{x\in k_{v}^{n}\,;\ \Vert\varphi(x)\Vert_{F,v}\le 1\right\}\ \cdotp$$Ce qui, pour une mesure de Haar $\vol_{v}$ quelconque sur $k_{v}^{n}$, entra{\^\i}ne $$\left(\frac{1}{\Vert\varphi\Vert_{v}}\right)^{n_{v}n}\frac{\vol_{v}(\mathbf{B}(k_{v}^{n},\Vert.\Vert_{E,v}))}{\vol_{v}(\mathbf{B}(k_{v}^{n},\vert.\vert_{2,v}))}\le \frac{\vol_{v}(\mathbf{B}(k_{v}^{n},\Vert.\Vert_{F,v}))}{\vol_{v}(\mathbf{B}(k_{v}^{n},\vert.\vert_{2,v}))}\vert\det\varphi\vert_{v}^{-n_{v}}\ \cdotp$$L'in{\'e}galit{\'e} souhait{\'e}e s'en d{\'e}duit au moyen de la formule du produit appliqu{\'e}e {\`a} $\det\varphi\in k\setminus\{0\}$. \end{proof}
\begin{lemma}\label{lemmefondamentalpentes}
Soit $\overline{E},\overline{F}$ des fibr\'es vectoriels ad\'eliques et $\varphi:E\to F$ une application lin\'eaire. Si $\varphi$ est injective alors \begin{equation}\widehat{\mu}_{\mathrm{max}}(\overline{E})\le\widehat{\mu}_{\mathrm{max}}(\overline{F})+h(\varphi)\ .\end{equation}
\end{lemma}
\begin{rema}La convention $\widehat{\mu}(0)=-\infty$ trouve une justification ici car l'on peut choisir $E=\{0\}$ dans cet \'enonc\'e.\end{rema} 
\begin{proof}
Soit $\overline{E'}\subseteq\overline{E}$. On pose $F'=\varphi(E')$, espace vectoriel sur $k$ que l'on munit des m{\'e}triques induites par $\overline{F}$. Ainsi, gr{\^a}ce {\`a} l'hypoth{\`e}se d'injectivit{\'e}, $\varphi$ induit un isomorphisme $\varphi_{\vert E'}^{\vert F'}$ entre $E'$ et $F'$. La majoration du lemme pr{\'e}c{\'e}dent donne alors \begin{equation*}\begin{split}\widehat{\mu}_{\mathrm{n}}(\overline{E'})&\le\widehat{\mu}_{\mathrm{n}}(\overline{F'})+h\left(\overline{E'},\overline{F'};\varphi_{\vert E'}^{\vert F'}\right)\\ &\le \widehat{\mu}_{\mathrm{max}}(\overline{F})+h(\overline{E},\overline{F};\varphi)\ .\end{split}\end{equation*}On choisit alors convenablement $E'$ pour obtenir la pente maximale de $\overline{E}$ {\`a} gauche. 
\end{proof}
Il y a plusieurs fa{\c c}ons d'{\'e}tendre cette in{\'e}galit{\'e}. La premi{\`e}re --- qui est c{\oe}ur de ce que l'on appelle \emph{m{\'e}thode des pentes} --- consiste {\`a} consid{\'e}rer une filtration $$0=F_{N}\subseteq F_{N-1}\subseteq\cdots\subseteq F_{0}:=F$$ d'un espace $F$ par des $k$-espaces vectoriels $F_{i}$. On ne suppose pas que $F$ est muni d'une structure de fibr{\'e} vectoriel ad{\'e}lique. En revanche, chacun des espaces quotients $G_{i}:=F_{i-1}/F_{i}$, $1\le i\le N$, est muni de normes $\Vert.\Vert_{i,v}$ aux places de $k$ telles que $(G_{i},(\Vert.\Vert_{i,v})_{v})$ soit un fibr{\'e} vectoriel ad{\'e}lique. Soit $\overline{E}$ un fibr{\'e} vectoriel ad{\'e}lique sur $\spec k$ de dimension $n$ et $\varphi: E\to F$ une application $k$-lin{\'e}aire. On pose $$\forall\,i\in\{1,\ldots,N+1\},\quad E_{i}:=\varphi^{-1}(F_{i-1})\ .$$Chacun des espaces $E_{i}$ est pourvu de la structure ad{\'e}lique $\overline{E_{i}}$ induite par $\overline{E}$. On pose $E_{N+1}:=\{0\}$. Soit $\varphi_{i}:E_{i}\to G_{i}$ la compos{\'e}e de la restriction de $\varphi$ {\`a} $E_{i}$ et de la projection canonique $F_{i-1}\to F_{i-1}/F_{i}$. 
\begin{prop}\label{propositionmethodedespentes}
Avec les notations ci-dessus, si $\varphi$ est injective alors \begin{equation}\label{inegalitedepentesgenerales}\widehat{\mu}_{\mathrm{n}}(\overline{E})\le\sum_{i=1}^{N}{\frac{\dim(E_{i}/E_{i+1})}{n}\left\{\widehat{\mu}_{\mathrm{max}}(\overline{G_{i}})+h(\overline{E_{i}},\overline{G_{i}};\varphi_{i})\right\}}+\frac{1}{D}\log\vr(\overline{E})\ \cdotp
\end{equation}
\end{prop}
Cette majoration n'est pas tout {\`a} fait une g{\'e}n{\'e}ralisation du lemme~\ref{lemmefondamentalpentes} {\`a} cause du terme $\frac{1}{D}\log\vr(\overline{E})$ qui subsiste m{\^e}me lorsque $N=1$.
\begin{proof}
Chacune des applications $\varphi_{i}$ se factorise en une injection $\widetilde{\varphi}_{i}:E_{i}/E_{i+1}\hookrightarrow G_{i}$ {\`a} laquelle l'on peut appliquer le lemme~\ref{lemmefondamentalpentes}. On a donc en particulier \begin{equation*}\begin{split}\widehat{\deg}_{\mathrm{n}}\left(\overline{E_{i}/E_{i+1}}\right)&\le\dim\left(E_{i}/E_{i+1}\right)\left\{\widehat{\mu}_{\mathrm{max}}(\overline{G_{i}})+h(\overline{E_{i}/E_{i+1}},\overline{G_{i}};\widetilde{\varphi}_{i})\right\}\\ &\le\dim\left(E_{i}/E_{i+1}\right)\left\{\widehat{\mu}_{\mathrm{max}}(\overline{G_{i}})+h(\overline{E_{i}},\overline{G_{i}};\varphi_{i})\right\}\ .\end{split}\end{equation*}La difficult{\'e} est que le degr{\'e} du quotient {\`a} gauche n'est pas en g{\'e}n{\'e}ral la diff{\'e}rence des degr{\'e}s \emph{sauf} si $\overline{E_{i}}$ est hermitien (voir proposition~\ref{propositiondegrequotient46}). Autrement dit, si $\overline{E}$ est hermitien, la proposition est {\'e}tablie en sommant les in{\'e}galit{\'e}s ci-dessus de $i=1$ {\`a} $N$. Le cas g{\'e}n{\'e}ral s'en d{\'e}duit alors en appliquant l'in{\'e}galit{\'e}~\eqref{inegalitedepentesgenerales} {\`a} $J(\overline{E})$ et en utilisant la proposition~\ref{prop:comparaisonjohnadelique} ainsi que la majoration $h((E_{i},\vert.\vert_{J(\overline{E})}),\overline{G_{i}};\varphi_{i})\le h(\overline{E_{i}},\overline{G_{i}};\varphi_{i})$.
\end{proof}
Une autre fa{\c c}on de g{\'e}n{\'e}raliser le lemme~\ref{lemmefondamentalpentes} est de consid{\'e}rer des pentes autres que la pente maximale, et m{\^e}me de supprimer l'hypoth{\`e}se d'injectivit{\'e} de $\varphi$, au prix de quelques modifications donn{\'e}es dans l'{\'e}nonc{\'e} suivant.
\begin{prop}\label{proposition66au}
Soit $\overline{E},\overline{F}$ des fibr{\'e}s vectoriels ad{\'e}liques et $\varphi:E\to F$ une application $k$-lin{\'e}aire. Alors, pour tout entier $i\ge 1$ inf\'erieur au rang de $\varphi$, la pente $\widehat{\mu}_{i+\dim\ker\varphi}(\overline{E})$ est inf\'erieure \`a \begin{equation*}\widehat{\mu}_{i}(\overline{F})+\frac{i}{D}\log\Delta(\overline{F})+\frac{(i+\dim\ker\varphi)}{D}\log\Delta(\overline{E})+h(\overline{E},\overline{F};\varphi)\ \cdotp\end{equation*}
\end{prop}
\begin{proof}
On commence par montrer que $$\widehat{\mu}_{i+\dim\ker\varphi}(\overline{E})\le\widehat{\mu}_{i}(\overline{F})+h(\overline{E},\overline{F};\varphi)$$lorsque $\overline{E}$ et $\overline{F}$ sont hermitiens. Tout d'abord, montrons que l'on peut supposer que $\varphi$ est injective. En effet, le quotient $E/\ker\varphi$ s'injecte dans $F$. Gr{\^a}ce {\`a} la formule~\eqref{ineqminmax} du minimax (qui est une {\'e}galit{\'e} dans le cas hermitien), pour tout $\epsilon>0$, il existe un sous-fibr{\'e} ad{\'e}lique hermitien $\overline{E}_{\epsilon}$ de $\overline{E}$, de dimension inf{\'e}rieure {\`a} $\dim\ker\varphi+i-1$ et contenant $\ker\varphi$, tel que \begin{equation*}\widehat{\mu}_{\mathrm{max}}\overline{\left(\frac{E/\ker\varphi}{E_{\epsilon}/\ker\varphi}\right)}\le\widehat{\mu}_{i}\left(\overline{E/\ker\varphi}\right)+\epsilon\ .\end{equation*}Or l'application $p:E/E_{\epsilon}\to(E/\ker\varphi)/(E_{\epsilon}/\ker\varphi)$ est un isomorphisme de hauteur n{\'e}gative car la norme d'op{\'e}rateur de $p$ est plus petite que $1$ en chaque place de $k$. La pente maximale de $$\overline{\left(\frac{E/\ker\varphi}{E_{\epsilon}/\ker\varphi}\right)}$$est donc sup{\'e}rieure {\`a} celle de $\overline{E/E_{\epsilon}}$ par le lemme~\ref{lemmefondamentalpentes}, elle-m{\^e}me plus grande que $\widehat{\mu}_{i+\dim\ker\varphi}(\overline{E})$ toujours en vertu de la formule~\eqref{ineqminmax}. On peut alors faire tendre $\epsilon$ vers $0$ et constater que $$\widehat{\mu}_{i+\dim\ker\varphi}(\overline{E})\le\widehat{\mu}_{i}\left(\overline{E/\ker\varphi}\right)\ \cdotp$$Comme $h(\overline{E/\ker\varphi},\overline{F};\varphi)\le h(\overline{E},\overline{F};\varphi)$, il suffit donc bien de traiter le cas $\ker\varphi=\{0\}$. Dans ce cas, on utilise {\`a} nouveau la formule~\eqref{ineqminmax} en observant que si $E_{2}\subseteq E_{1}\subseteq E$ avec $\dim E_{1}\ge i>\dim E_{2}$ alors $\varphi(E_{2})\subseteq\varphi(E_{1})\subseteq F$ avec $\dim\varphi(E_{1})\ge i>\dim\varphi(E_{2})$. La pente maximale de $\overline{E/E_{2}}$ est {\`a} la fois plus petite que $\widehat{\mu}_{\mathrm{max}}(\overline{F/\varphi(E_{2})})+h(\overline{E},\overline{F};\varphi)$ (lemme~\ref{lemmefondamentalpentes}) et plus grande que $\widehat{\mu}_{i}(\overline{E})$ (formule~\eqref{ineqminmax}). D'o\`u l'on d{\'e}duit $\widehat{\mu}_{i}(\overline{E})\le\widehat{\mu}_{i}(\overline{F})+h(\overline{E},\overline{F};\varphi)$ en faisant varier $E_{2}$, l'espace $\varphi(E_{2})$ parcourant les sous-espaces vectoriels de $\varphi(E)$ de dimension inf{\'e}rieure {\`a} $i-1$. \par 
Une fois la majoration {\'e}tablie dans le cas hermitien, on l'applique aux fibr{\'e}s ad{\'e}liques hermitiens $\overline{E}_{\varepsilon}$ et $\overline{F}_{\varepsilon}$ (introduits apr\`es les relations~\eqref{eqref:important}), puis l'on utilise l'encadrement~\eqref{comparaisonpentesijohn} ainsi que les {\'e}valuations~\eqref{evaluationunhauteur} et~\eqref{evaluationdeuxhauteur} de la hauteur de $\varphi$ pour en d\'eduire le cas g\'en\'eral.
\end{proof}
\begin{coro}
Soit $\overline{E},\overline{F}$ des fibr\'es vectoriels ad\'eliques et $\varphi:E\to F$ une application $k$-lin\'eaire. On note $m$ la dimension de $F$. Si $\varphi$ est surjective alors la pente maximale $\widehat{\mu}_{\mathrm{max}}(\overline{F})$ est inf\'erieure \`a\begin{equation*}\widehat{\deg}_{\mathrm{n}}\overline{F}-(m-1)\widehat{\mu}_{\mathrm{min}}(\overline{E})+(m-1)h(\overline{E},\overline{F};\varphi)+\frac{m}{D}\log\left(\Delta(\overline{E})\Delta(\overline{F})\right)\ \cdotp\end{equation*}  
\end{coro}
\begin{proof}
Des formules $\widehat{\deg}_{\mathrm{n}}\overline{F}=\sum_{i=1}^{m}{\widehat{\mu}_{i}(\overline{F})}$ et $\widehat{\mu}_{1}(\overline{F})=\widehat{\mu}_{\mathrm{max}}(\overline{F})$ (lemme~\ref{lemmamuunmumax}) l'on d\'eduit la majoration\begin{equation}\label{majorationpentedegre}\widehat{\mu}_{\mathrm{max}}(\overline{F})\le\widehat{\deg}_{\mathrm{n}}\overline{F}-(m-1)\widehat{\mu}_{m}(\overline{F})\ .\end{equation}La proposition~\ref{proposition66au} appliqu\'ee avec $i=m=\rg\varphi$ fournit une minoration de $\widehat{\mu}_{m}(\overline{F})$ qui est trop faible pour prouver le corollaire \`a cause des termes d'erreurs qui se sont ajout\'es. On va donc commencer par utiliser l'estimation~\eqref{majorationpentedegre} avec le fibr\'e hermitien $\overline{F}_{\varepsilon}$ construit \`a la suite de~\eqref{eqref:important} au lieu de $\overline{F}$. On a vu que $\widehat{\mu}_{\mathrm{max}}(\overline{F})\le\widehat{\mu}_{\mathrm{max}}(\overline{F}_{\varepsilon})$ (voir ce qui suit la d\'efinition~\ref{definitiondelapentemaximale}) et la proposition~\ref{proposition66au} appliqu\'ee au morphisme $\varphi:\overline{E}_{\varepsilon}\to\overline{F}_{\varepsilon}$ et $i=m=\rg\varphi$ fournit la minoration\begin{equation*}\widehat{\mu}_{m}(\overline{F}_{\varepsilon})\ge\widehat{\mu}_{n}(\overline{E}_{\varepsilon})-h(\overline{E}_{\varepsilon},\overline{F}_{\varepsilon};\varphi)\ .\end{equation*}On obtient ainsi \begin{equation}\begin{split}\label{tyeree}\widehat{\mu}_{\mathrm{max}}(\overline{F})\le\widehat{\mu}_{\mathrm{max}}(\overline{F}_{\varepsilon})&\le\widehat{\deg}_{\mathrm{n}}\overline{F}_{\varepsilon}-(m-1)\widehat{\mu}_{m}(\overline{F}_{\varepsilon})\\ & \le\widehat{\deg}_{\mathrm{n}}\overline{F}_{\varepsilon}-(m-1)\widehat{\mu}_{n}(\overline{E}_{\varepsilon})+(m-1)h(\overline{E}_{\varepsilon},\overline{F}_{\varepsilon};\varphi)\ .\end{split}\end{equation}D'apr\`es l'observation~\eqref{eqref:important}, on a $\widehat{\deg}_{\mathrm{n}}\overline{F}_{\varepsilon}\le\widehat{\deg}_{\mathrm{n}}\overline{F}+\frac{m}{D}\log\Delta(\overline{F})+m\log(1+\varepsilon)$. De plus, les estimations~\eqref{evaluationunhauteur} et~\eqref{evaluationdeuxhauteur} pour la hauteur donnent \begin{equation*}h(\overline{E}_{\varepsilon},\overline{F}_{\varepsilon};\varphi)\le h(\overline{E}_{\varepsilon},\overline{F};\varphi)\le h(\overline{E},\overline{F};\varphi)+\frac{1}{D}\log\Delta(\overline{E})+\log(1+\varepsilon)\ .\end{equation*}Enfin, gr{\^a}ce au lemme~\ref{lemmeliti} pour le fibr\'e hermitien $\overline{E}_{\varepsilon}$, on a \begin{equation*}\widehat{\mu}_{n}(\overline{E}_{\varepsilon})=\min_{E'\subsetneq E}{\widehat{\mu}_{\mathrm{n}}(E/E',\vert.\vert_{\varepsilon})}\ge \min_{E'\subsetneq E}{\widehat{\mu}_{\mathrm{n}}(E/E',\Vert.\Vert_{\overline{E}})}=\widehat{\mu}_{\mathrm{min}}(\overline{E})\ .\end{equation*}La majoration souhait\'ee pour la pente maximale $\widehat{\mu}_{\mathrm{max}}(\overline{F})$ d\'ecoule de ces estimations que l'on reporte dans la majoration~\eqref{tyeree}, en faisant tendre ensuite $\varepsilon$ vers $0$. 
\end{proof}

\section{Pentes maximales des puissances sym{\'e}triques}\label{subsectionpenteproduitsymetrique} L'objectif de ce paragraphe est de montrer l'{\'e}nonc{\'e} suivant.
\begin{theo}\label{pentefibresymetrique}Soit $\overline{E}$ un fibr{\'e} ad{\'e}lique \emph{hermitien} de dimension $n\ge 1$. Pour tout entier $\ell\ge 0$, la pente maximale de la $\ell^{\text{{\`e}me}}$ puissance sym{\'e}trique $S^{\ell}(\overline{E})$ v{\'e}rifie\noignorespaces 
\begin{enumerate}\item[1)] si $k$ est un corps de fonctions alors on a $\widehat{\mu}_{\mathrm{max}}(S^{\ell}(\overline{E}))=\ell\widehat{\mu}_{\mathrm{max}}(\overline{E})$. 
\item[2)] si $k$ est un corps de nombres alors 
\begin{equation*}0\le\widehat{\mu}_{\mathrm{max}}(S^{\ell}(\overline{E}))-\ell\widehat{\mu}_{\mathrm{max}}(\overline{E})\le 2\ell n\log n\ .\end{equation*}\end{enumerate}
\end{theo} Avant d'effectuer la d{\'e}monstration de ce th{\'e}or{\`e}me, rassemblons deux {\'e}nonc{\'e}s pr{\'e}paratoires, d'int{\'e}r{\^e}t ind{\'e}pendant. \par
Pour $\mathbf{i}=(i_{1},\ldots,i_{n})\in\mathbf{N}^{n}$, on note $\vert\mathbf{i}\vert$ la \emph{longueur} de $\mathbf{i}$ d{\'e}finie par $\vert\mathbf{i}\vert:=i_{1}+\cdots+i_{n}$ et $\mathbf{i} !:=i_{1}!\cdots i_{n}!$. {\'E}tant donn{\'e} $\ell\in\mathbf{N}$, si $A$ est la matrice d'un endomorphisme $u$ dans une base $(e_{1},\ldots,e_{n})$ de $k^{n}$, la matrice $S^{\ell}(A)$ est celle qui repr{\'e}sente $S^{\ell}(u)$ dans la base $e_{1}^{i_{1}}\cdots e_{n}^{i_{n}}$, o\`u les $n$-uplets $(i_{1},\ldots,i_{n})\in\mathbf{N}^{n}$ sont de longueur $\ell$ et ordonn{\'e}s lexicographiquement.
\begin{lemma}\label{lemme63determinant} Soit $v$ une place de $k$ et $u_{v}$ un endomorphisme de $k_{v}^{n}$. On suppose $k_{v}^{n}$ muni de la norme $\vert.\vert_{2,v}$ et on note $\Vert.\Vert_{v}$ la norme d'op\'erateur des endomorphismes de $(k_{v}^{n},\vert.\vert_{2,v})$ si $v$ est ultram\'etrique ou la norme de Hilbert-Schmidt si $v$ est archim\'edienne. 
\begin{enumerate}
\item[(i)] Si $u_{v}$ est un isomorphisme alors \begin{equation}\label{formuleinversematrice}\Vert u_{v}^{-1}\Vert_{v}\le\frac{\Vert u_{v}\Vert_{v}^{n-1}}{\vert\det u_{v}\vert_{v}}\ \cdotp\end{equation}
\item[(ii)] Pour tout entier naturel $\ell$, on a \begin{equation*}\det S^{\ell}(u_{v})=(\det u_{v})^{\binom{\ell+n-1}{n}}\ \cdotp\end{equation*}
\end{enumerate}
\end{lemma}
\begin{proof}[\'El\'ements de d\'emonstration] (i) En une place finie $v$, si l'on note $A_{v}$ la matrice qui repr\'esente $u_{v}$ dans la base canonique de $k_{v}^{n}$, l'estimation~\eqref{formuleinversematrice} est la cons\'equence de la formule matricielle $A_{v}^{-1}=\frac{{}^{\mathrm{t}}\mathrm{com}A_{v}}{\det A_{v}}$ o\`u ${}^{\mathrm{t}}\mathrm{com}A_{v}$ d{\'e}signe la transpos{\'e}e de la comatrice de $A_{v}$. Si $v$ est archim\'edienne, on observe que la norme de Hilbert-Schmidt $\Vert u_{v}\Vert_{v}$ est la racine carr\'ee de la somme des valeurs propres de l'op\'erateur hermitien ${}^{\mathrm{t}}\overline{u_{v}}u_{v}$ qui est d\'efini positif. En remarquant que ${}^{\mathrm{t}}\overline{u_{v}}u_{v}$ et $u_{v}{}^{\mathrm{t}}\overline{u_{v}}$ ont le m{\^e}me spectre $\lambda_{1}\le\cdots\le\lambda_{n}$, l'in\'egalit\'e~\eqref{formuleinversematrice} n'est rien d'autre que \begin{equation*}\left(\sum_{i=1}^{n}{\lambda_{i}^{-1}}\right)(\lambda_{1}\cdots\lambda_{n})\le\left(\sum_{i=1}^{n}{\lambda_{i}}\right)^{n-1}\ .\end{equation*}\noindent (ii) Pour d\'emontrer cette formule, on peut par exemple choisir une base trigonalisante de $u_{v}$ dans une extension alg{\'e}brique de $k_{v}$ et observer que le mon{\^o}me sym{\'e}trique $\prod_{\vert\mathbf{i}\vert=\ell}{X_{1}^{i_{1}}\cdots X_{n}^{i_{n}}}$ {\'e}gale $(X_{1}\cdots X_{n})^{\binom{\ell+n-1}{n}}$.
\end{proof}
 Le second lemme requiert deux notations suppl{\'e}mentaires. On pose \begin{equation}\label{defidedelta}\delta:=\begin{cases} 0 & \text{si $k$ est un corps de fonctions}\\ 1 & \text{si $k$ est un corps de nombres}\end{cases}\end{equation}et, pour $n,\ell\in\mathbf{N}$, on d{\'e}finit \begin{equation*}\gamma_{n,\ell}:=\left(\prod_{\vert\mathbf{i}\vert=\ell}{\frac{\ell !}{\mathbf{i} !}}\right)^{\binom{\ell+n-1}{n-1}^{-1}}\qquad\text{{\scriptsize(dans ce produit, $\mathbf{i}\in\mathbf{N}^{n}$).}} \end{equation*}On conna{\^\i}t une estimation asymptotique du logarithme de $\gamma_{n,\ell}$ lorsque $\ell\to+\infty$ et $n$ fix{\'e}, que l'on peut {\'e}crire en fonction des nombres harmoniques $H_{n}:=\sum_{i=1}^{n}{\frac{1}{i}}$. La quantit\'e $\log\gamma_{n,\ell}$ est alors \'egale {\`a} $(H_{n}-1)(\ell+\mathrm{o}(\ell))$ lorsque $\ell\to+\infty$  (voir annexe).
\begin{lemma}\label{corollaire64}Soit $\ell\in\mathbf{N}\setminus\{0\}$ et $\overline{E}$ un fibr{\'e} vectoriel ad{\'e}lique de dimension $n\ge 1$. Alors, pour toute norme tensorielle ad\'elique hermitienne $\alpha$ d'ordre $\ell$, on a\begin{equation*}\left\vert\widehat{\mu}_{\mathrm{n}}\left(S_{\alpha}^{\ell}(\overline{E})\right)-\ell\widehat{\mu}_{\mathrm{n}}(\overline{E})-\frac{\delta}{2}\log\gamma_{n,\ell}\right\vert\le\frac{\ell}{D}\log\Delta(\overline{E})\ .\end{equation*}  
\end{lemma}
\begin{proof}
Supposons dans un premier temps que $\overline{E}$ est hermitien. Dans ce cas, on peut omettre la r\'ef\'erence \`a $\alpha$ dans la notation $S_{\alpha}^{\ell}(\overline{E})$ (qui n'en d\'epend pas, car $\alpha$ est hermitienne). Il s'agit de montrer que \begin{equation*}\widehat{\mu}_{\mathrm{n}}\left(S^{\ell}(\overline{E})\right)=\ell\widehat{\mu}_{\mathrm{n}}(\overline{E})+\frac{\delta}{2}\log\gamma_{n,\ell}\ .\end{equation*}Quitte {\`a} fixer une $k$-base de $E$, on peut identifier $E$ {\`a} $k^{n}$. Il existe une matrice ad{\'e}lique $A=(A_{v})_{v}\in\GL_{n}(k_{\mathbf{A}})$ telle que, pour toute place $v$ de $k$ et tout $x\in k_{v}^{n}$, on ait $\Vert x\Vert_{E,v}=\vert A_{v}.x\vert_{2,v}$. Les m{\'e}triques sur $S^{\ell}(\overline{E})$ sont d{\'e}termin{\'e}es par $S^{\ell}(A)$, au sens o\`u $S^{\ell}(\overline{E})$ est l'image par $S^{\ell}(A)$ du fibr{\'e} ad{\'e}lique hermitien $S_{\alpha}^{\ell}(k^{n},\vert.\vert_{2})$. L'{\'e}galit{\'e}~\eqref{casdegaliteimportant} et le lemme~\ref{lemme4225} entra{\^\i}nent alors \begin{equation*}\widehat{\deg}_{\mathrm{n}}S^{\ell}(k^{n},\vert.\vert_{2})=\widehat{\deg}_{\mathrm{n}}S^{\ell}(\overline{E})+\frac{1}{D}\sum_{v}{n_{v}\log\frac{1}{\left\vert\det S^{\ell}(A_{v})\right\vert_{v}}}\ \cdotp\end{equation*}Le lemme~\ref{lemme63determinant} (ii) permet de calculer le d{\'e}terminant de $S^{\ell}(A_{v})$. En observant que $\dim S^{\ell}(E)=\binom{\ell+n-1}{n-1}$ et en appliquant {\`a} nouveau l'{\'e}galit{\'e}~\eqref{casdegaliteimportant} {\`a} $\overline{E}$ et $(k^{n},\vert.\vert_{2})$, on obtient \begin{equation*}\widehat{\mu}_{\mathrm{n}}\left(S^{\ell}(\overline{E})\right)=\ell\widehat{\mu}_{\mathrm{n}}(\overline{E})+\widehat{\mu}_{\mathrm{n}}\left(S^{\ell}(k^{n},\vert.\vert_{2})\right)\ .\end{equation*}Ce dernier terme donne la contribution $\frac{1}{2}\log\gamma_{n,\ell}$ aux ({\'e}ventuelles) places archim{\'e}diennes de $k$. Ceci est une cons{\'e}quence du calcul local suivant. Soit $v$ une place de $k$. Si $(e_{1},\ldots,e_{n})$ d{\'e}signe la base canonique de $k^{n}$ alors, pour tout $\mathbf{i}=(i_{1},\ldots,i_{n})\in\mathbf{N}^{n}$, la norme de l'{\'e}l{\'e}ment $e_{1}^{i_{1}}\cdots e_{n}^{i_{n}}\in S^{\ell}\left(k^{n},\vert.\vert_{2}\right)$ vaut $(\mathbf{i}!/\ell !)^{1/2}$ si $v$ est archim\'edienne et $1$ sinon (voir \S~\ref{paragraphetroistrois}). La formule souhait{\'e}e {\'e}tant {\'e}tablie dans le cas hermitien, le passage au cas g{\'e}n{\'e}ral s'effectue gr{\^a}ce aux fibr{\'e}s de John et L{\"o}wner associ{\'e}s {\`a} $\overline{E}$, en observant que $S^{\ell}(L(\overline{E}))\preceq S_{\alpha}^{\ell}(\overline{E})\preceq S^{\ell}(J(\overline{E}))$. On obtient alors \begin{equation*}-\frac{\ell}{D}\log\vr(\overline{E})\le\widehat{\mu}_{\mathrm{n}}\left(S_{\alpha}^{\ell}(\overline{E})\right)-\ell\widehat{\mu}_{\mathrm{n}}(\overline{E})-\frac{\delta}{2}\log\gamma_{n,\ell}\le\frac{\ell}{D}\log\widetilde{\vr}(\overline{E})\end{equation*}(voir~\eqref{comparaisonsousfibres} et la proposition~\ref{prop:comparaisonjohnadelique}). On conclut au moyen de la majoration $\max{\{\vr(\overline{E}),\widetilde{\vr}(\overline{E})\}}\le\Delta(\overline{E})$ (voir~\eqref{eq:encadrementquotientadelique}).
\end{proof}
\begin{proof}[D{\'e}monstration du th{\'e}or{\`e}me~\ref{pentefibresymetrique}]
La premi{\`e}re in{\'e}galit{\'e} $\ell\widehat{\mu}_{\mathrm{max}}(\overline{E})\le\widehat{\mu}_{\mathrm{max}}(S^{\ell}(\overline{E}))$ est une cons{\'e}quence directe du lemme~\ref{corollaire64} (cas hermitien)~: \begin{equation*}\forall\,\overline{F}\subseteq\overline{E},\quad\ell\widehat{\mu}_{\mathrm{n}}(\overline{F})\le\widehat{\mu}_{\mathrm{n}}(S^{\ell}(\overline{F}))\le\widehat{\mu}_{\mathrm{max}}(S^{\ell}(\overline{E}))\ .\end{equation*}Le r{\'e}sultat se d{\'e}duit de la d{\'e}finition de la pente maximale de $\overline{E}$. Montrons maintenant que \begin{equation*}\widehat{\mu}_{\mathrm{max}}(S^{\ell}(\overline{E}))-\ell\widehat{\mu}_{\mathrm{max}}(\overline{E})\le 2\ell\delta n\log n\end{equation*}($\delta$ est d\'efini par~\eqref{defidedelta}). D'apr{\`e}s le lemme de Siegel absolu~\ref{lemmedesiegelabsolu}, pour tout $\varepsilon>0$, il existe une $\overline{k}$-base $e_{1},\ldots,e_{n}$ de $E\otimes\overline{k}$ telle que $$h_{\overline{E}}(e_{1})+\cdots+h_{\overline{E}}(e_{n})+\widehat{\deg}_{\mathrm{n}}\det\overline{E}\le\frac{\delta n}{2}\log n+\varepsilon\ .$$ On peut supposer que les vecteurs $e_{1},\ldots,e_{n}$ sont d{\'e}finis sur une extension finie $K$ de $k$. Soit $\overline{E_{0}}$ le fibr{\'e} ad{\'e}lique hermitien dont l'espace sous-jacent $E_{0}$ est $E_{K}$ et dont la norme d'un vecteur $x=\sum_{i}{x_{i}e_{i}}\in E\otimes_{k}\mathbf{C}_{v}$ en une place $w$ de $K$ au-dessus de $v$ est $(\sum_{i}{\vert x_{i}\vert_{v}^{2}\Vert e_{i}\Vert_{E,v}^{2}})^{1/2}$ si $w$ est archim{\'e}dienne et $\max_{i}{\{\vert x_{i}\vert_{v}\Vert e_{i}\Vert_{E,v}\}}$ sinon. Comme $\overline{E_{0}}$ et $\overline{E_{K}}$ sont des fibr{\'e}s hermitiens avec le m{\^e}me espace sous-jacent, il existe une matrice $\Sigma=(\Sigma_{w})_{w}\in\GL(E\otimes_{k}K_{\mathbf{A}})$ telle que $\overline{E_{K}}=\Sigma.\overline{E_{0}}$ (au sens donn\'e peu apr{\`e}s la d{\'e}finition du degr{\'e} ad{\'e}lique, voir~\eqref{transfoechelle}). On note encore $\Sigma:E_{0}\to E_{K}$ l'application identique. Cet abus de notation se justifie par le fait que les normes d'op{\'e}rateur ou de Hilbert-Schmidt de l'application $x\in E_{0}\otimes K_{w}\mapsto x\in E\otimes K_{w}$ sont celles de $\Sigma_{w}$ pour toute place $w$ de $K$. En appliquant le lemme~\ref{lemmefondamentalpentes} {\`a} l'application inverse de $S^{\ell}(\Sigma):S^{\ell}(E_{0})\to S^{\ell}(E_{K})$, on obtient~:\begin{equation*}\begin{split}\widehat{\mu}_{\mathrm{max}}\left(S^{\ell}(\overline{E})\right)=\widehat{\mu}_{\mathrm{max}}\left(S^{\ell}(\overline{E_{K}})\right)&\le\widehat{\mu}_{\mathrm{max}}\left(S^{\ell}(\overline{E_{0}})\right)+h\left(S^{\ell}(\Sigma^{-1})\right)\\ &\le\widehat{\mu}_{\mathrm{max}}\left(S^{\ell}(\overline{E_{0}})\right)+\ell h\left(\Sigma^{-1}\right)\ \cdotp\end{split}\end{equation*}\'Evaluons chacun des termes du membre de droite de cette in{\'e}galit{\'e}. La pente maximale de $S^{\ell}(\overline{E_{0}})$ se calcule en observant que $S^{\ell}(E_{0})$ est la somme directe orthogonale des espaces $K.e_{1}^{i_{1}}\cdots e_{n}^{i_{n}}$ pour $\mathbf{i}:=(i_{1},\ldots,i_{n})$ de longueur $\ell$. D'apr{\`e}s la propri{\'e}t{\'e}~\ref{proprietes55}, 2), on a \begin{equation*}\begin{split}\widehat{\mu}_{\mathrm{max}}\left(S^{\ell}(\overline{E_{0}})\right)&=\max_{\vert\mathbf{i}\vert=\ell}{\left\{\widehat{\deg}_{\mathrm{n}}(K.e_{1}^{i_{1}}\cdots e_{n}^{i_{n}})\right\}}\\ & =\max_{\vert\mathbf{i}\vert=\ell}{\left\{\sum_{j=1}^{n}{i_{j}\widehat{\mu}_{\mathrm{max}}(\overline{K.e_{j}})}+\frac{\delta}{2}\log\frac{\ell !}{i_{1}!\cdots i_{n}!}\right\}}\ \cdotp
\end{split}
\end{equation*}On en d{\'e}duit \begin{equation*}\begin{split}\widehat{\mu}_{\mathrm{max}}\left(S^{\ell}(\overline{E_{0}})\right)&\le\ell\max_{1\le i\le n}{\left\{\widehat{\mu}_{\mathrm{max}}(\overline{K.e_{i}})\right\}}+\frac{\delta\ell}{2}\log n\\ & =\ell\widehat{\mu}_{\mathrm{max}}(\overline{E_{0}})+\frac{\delta\ell}{2}\log n\ .\end{split}\end{equation*}Le lemme~\ref{lemmefondamentalpentes} appliqu{\'e} {\`a} $\Sigma:E_{0}\to E_{K}$ donne $\widehat{\mu}_{\mathrm{max}}(\overline{E_{0}})\le\widehat{\mu}_{\mathrm{max}}(\overline{E})+h(\Sigma)$. La hauteur $h(\Sigma)$ est major{\'e}e par $\frac{\delta}{2}\log n$ car $\Vert\Sigma\Vert_{w}\le\sqrt{n}$ si $w$ est archim\'edienne et $\Vert\Sigma\Vert_{w}\le 1$ sinon. On obtient ainsi $$\widehat{\mu}_{\mathrm{max}}\left(S^{\ell}(\overline{E_{0}})\right)\le\ell\left(\widehat{\mu}_{\mathrm{max}}(\overline{E})+\delta\log n\right)\ .$$Par ailleurs, en ce qui concerne la hauteur de $\Sigma^{-1}$, on applique la majoration~\eqref{formuleinversematrice} du lemme~\ref{lemme63determinant}. Compte tenu de l'estimation des normes de $\Sigma$, on en d\'eduit la majoration \begin{equation*}h\left(\Sigma^{-1}\right)\le\delta(n-1)\log\sqrt{n}-\frac{1}{D}\sum_{\text{$w$ place de $K$}}{n_{w}\log\vert\det\Sigma_{w}\vert_{w}}\ \cdotp\end{equation*}L'\'egalit\'e~\eqref{casdegaliteimportant} obtenue {\`a} la suite du lemme~\ref{lemme:egalitedesdegres} montre que cette derni{\`e}re somme est aussi la diff{\'e}rence des degr{\'e}s ad{\'e}liques $\widehat{\deg}_{\mathrm{n}}\overline{E}-\widehat{\deg}_{\mathrm{n}}\overline{E_{0}}$, qui vaut $\widehat{\deg}_{\mathrm{n}}\det\overline{E}+h_{\overline{E}}(e_{1})+\cdots+h_{\overline{E}}(e_{n})$ car $\overline{E}$ est hermitien et le degr\'e de $\overline{E_{0}}$ se calcule au moyen de la formule~\eqref{transfoechelle}. Cette quantit{\'e} est inf{\'e}rieure {\`a} $\frac{\delta n}{2}\log n+\varepsilon$ en vertu du choix de la base $(e_{1},\ldots,e_{n})$. La synth{\`e}se de ces informations conduit {\`a} l'estimation~:\begin{equation*}\begin{split}\widehat{\mu}_{\mathrm{max}}\left(S^{\ell}(\overline{E})\right)&\le\ell\left(\widehat{\mu}_{\mathrm{max}}(\overline{E})+\delta\log n\right)+\ell\left(\delta(n-1)\log\sqrt{n}+\frac{\delta n}{2}\log n+\varepsilon\right)\\ & \le\ell\left(\widehat{\mu}_{\mathrm{max}}(\overline{E})+\left(n+\frac{1}{2}\right)\delta\log n+\varepsilon\right)\ \cdotp\end{split}\end{equation*}Dans cette in{\'e}galit{\'e}, on peut faire tendre $\varepsilon$ vers $0$ et l'on obtient le r{\'e}sultat voulu.
\end{proof}
\begin{coro}\label{corollaire742007}
Pour tout entier $\ell\ge 1$, pour toute norme tensorielle ad\'elique \emph{hermitienne} $\alpha$ d'ordre $\ell$ et tout fibr{\'e} vectoriel ad{\'e}lique $\overline{E}$ sur $\spec k$, on a $$-\frac{\ell}{D}\log\Delta(\overline{E})\le\widehat{\mu}_{\mathrm{max}}\left(S_{\alpha}^{\ell}(\overline{E})\right)-\ell\widehat{\mu}_{\mathrm{max}}(\overline{E})\le \ell\left(2n\delta\log n+\frac{1}{D}\log\Delta(\overline{E})\right)\ .$$ 
\end{coro}
\begin{proof} En vertu du lemme~\ref{corollaire64}, pour tout sous-fibr{\'e} ad{\'e}lique $\overline{F}\subseteq\overline{E}$, on a \begin{equation*}\ell\widehat{\mu}_{\mathrm{n}}(\overline{F})\le\widehat{\mu}_{\mathrm{n}}\left(S_{\alpha}^{\ell}(\overline{F})\right)+\frac{\ell}{D}\log\Delta(\overline{F})\quad\text{et}\quad\Delta(\overline{F})\le\Delta(\overline{E})\ .\end{equation*}\textit{A priori} les m\'etriques sur $S_{\alpha}^{\ell}(\overline{F})$ ne sont pas les restrictions des m\'etriques de $S_{\alpha}^{\ell}(\overline{E})$ \`a $S^{\ell}(F)$. Toutefois, d'apr\`es la proposition~\ref{propositioncompnormes}, on a $(F^{\otimes\ell},\Vert.\Vert_{\overline{E}^{\otimes_{\alpha}\ell}})\preceq \overline{F}^{\otimes_{\alpha}\ell}$, d'o\`u l'on d\'eduit \begin{equation*}(S^{\ell}(F),\Vert.\Vert_{S_{\alpha}^{\ell}(\overline{E})})\preceq S_{\alpha}^{\ell}(\overline{F})\quad\text{puis}\quad\widehat{\mu}_{\mathrm{n}}\left(S_{\alpha}^{\ell}(\overline{F})\right)\le\widehat{\mu}_{\mathrm{max}}\left(S_{\alpha}^{\ell}(\overline{E})\right)\ .
\end{equation*}La premi\`ere in\'egalit\'e du corollaire~\ref{corollaire742007} en d\'ecoule. Pour la seconde estimation, on applique le th{\'e}or{\`e}me~\ref{pentefibresymetrique} au fibr{\'e} ad{\'e}lique hermitien $\overline{E}_{\varepsilon}$ (voir~\eqref{eqref:important}). Comme $\overline{E}_{\varepsilon}\preceq\overline{E}$, on en d{\'e}duit $S^{\ell}(\overline{E}_{\varepsilon})\preceq S_{\alpha}^{\ell}(\overline{E})$ (lemme~\ref{lemmecomptrois}) et donc \begin{equation*}\widehat{\mu}_{\mathrm{max}}\left(S_{\alpha}^{\ell}(\overline{E})\right)\le\widehat{\mu}_{\mathrm{max}}\left(S^{\ell}(\overline{E}_{\varepsilon})\right)\ .\end{equation*}D'autre part, comme nous l'avons mentionn{\'e} {\`a} la suite de la d{\'e}finition~\ref{definitiondelapentemaximale}, la pente maximale de $\overline{E}_{\varepsilon}$ est major{\'e}e par celle de $\overline{E}$ plus $\log(1+\varepsilon)+\frac{1}{D}\log\Delta(\overline{E})$. Le r{\'e}sultat s'ensuit en faisant tendre $\varepsilon$ vers $0$.
\end{proof}
\begin{rema}Si $k$ est un corps de nombres, l'{\'e}galit{\'e} $\widehat{\mu}_{\mathrm{max}}\left(S_{\alpha}^{\ell}(\overline{E})\right)=\ell\widehat{\mu}_{\mathrm{max}}(\overline{E})$ est en g{\'e}n{\'e}ral fausse, m{\^e}me asymptotiquement lorsque $\ell\to+\infty$ ou m{\^e}me si $\overline{E}$ et $\alpha$ sont hermitiens. En effet l'estimation asymptotique $\log\gamma_{n,\ell}\underset{\ell\to+\infty}{\sim}(H_{n}-1)\ell$ et le lemme~\ref{corollaire64} entra{\^\i}nent la majoration\begin{equation*}\widehat{\mu}_{\mathrm{n}}(\overline{F})+\frac{1}{2}(H_{m}-1)\le\frac{1}{D}\log\Delta(\overline{E})+\liminf_{\ell\to+\infty}{\frac{\widehat{\mu}_{\mathrm{max}}\left(S^{\ell}(\overline{E})\right)}{\ell}}\ ,\end{equation*}vraie pour tout sous-fibr{\'e} ad{\'e}lique $\overline{F}\subseteq\overline{E}$ de dimension $m\ge 1$. Aussi, lorsque la pente maximale de $\overline{E}$ est atteinte pour un fibr{\'e} $\overline{F}$ de dimension $\ge 2$ (par exemple lorsque $\overline{E}$ est semi-stable) et si $\overline{E}$ est hermitien ($\Delta(\overline{E})=1$), on a \begin{equation*}\widehat{\mu}_{\mathrm{max}}(\overline{E})<\liminf_{\ell\to+\infty}{\frac{\widehat{\mu}_{\mathrm{max}}\left(S^{\ell}(\overline{E})\right)}{\ell}}\ \cdotp\end{equation*}
\end{rema}
Ces propri{\'e}t{\'e}s des pentes et pentes maximales des puissances sym{\'e}triques de $\overline{E}$ se transmettent aux fibr{\'e}s des sections globales des puissances tensorielles du faisceau canonique $\mathcal{O}(1)$ de l'espace projectif $\mathbf{P}(E)$.\par Plus pr{\'e}cis{\'e}ment, soit $\overline{E}$ un fibr{\'e} vectoriel ad{\'e}lique sur $\spec k$ de dimension $n\ge 1$. On note $\mathbf{P}(E)$ le sch\'ema $\Proj\mathbf{S}(E)$ de morphisme structural $\pi:\mathbf{P}(E)\to\spec k$ et $\mathcal{O}(1)$ le fibr{\'e} en droites universel sur $\mathbf{P}(E)$. Pour toute place $v$ de $k$ et tout point $x\in\mathbf{P}(E)(\mathbf{C}_{v})$, la surjection canonique $\pi^{*}E\twoheadrightarrow\mathcal{O}(1)$ conf\`ere une m\'etrique $\Vert.\Vert_{v,x}$ \`a la fibre $\mathcal{O}(1)_{x}$ par quotient. Et plus g{\'e}n{\'e}ralement, pour $\ell\in\mathbf{N}\setminus\{0\}$, on note $\Vert.\Vert_{v,x,\ell}$ la norme induite sur $\mathcal{O}(\ell)_{x}=\mathcal{O}(1)_{x}^{\otimes\ell}$. En une place archim{\'e}dienne $v$ de $k$, une mesure de Haar sur le groupe localement compact $(E\otimes_{k}\mathbf{C}_{v},+)$ induit une mesure sur l'ouvert $(E\otimes_{k}\mathbf{C}_{v})\setminus\{0\}$, qui est un espace homog\`ene sous l'action du groupe multiplicatif $\mathbf{C}_{v}^{*}$. On obtient de la sorte une unique mesure de probabilit\'e $\mu_{E,v}$ sur $\mathbf{P}(E\otimes_{k}\mathbf{C}_{v})\simeq\left((E\otimes_{k}\mathbf{C}_{v})\setminus\{0\}\right)/\mathbf{C}_{v}^{*}$. On peut munir alors le $k$-espace vectoriel $E_{\ell}$ des sections globales $\mathrm{H}^{0}(\mathbf{P}(E),\mathcal{O}(\ell))$ d'une structure de fibr\'e ad\'elique hermitien sur $\spec k$ en consid{\'e}rant les normes suivantes. \'Etant donn\'e une place $v$ de $k$ et un {\'e}l{\'e}ment $s$ de $E_{\ell}\otimes_{k}\mathbf{C}_{v}$, on pose \begin{equation*}\Vert s\Vert_{E_{\ell},v}:=\begin{cases}\sup_{x\in\mathbf{P}(E)(\mathbf{C}_{v})}{\Vert s(x)\Vert_{v,x,\ell}} & \text{si $v$ est ultram\'etrique},\\ \left(\int_{\mathbf{P}(E)(\mathbf{C}_{v})}{\Vert s(x)\Vert_{v,x,\ell}^{2}\,\mathrm{d}\mu_{E,v}(x)}\right)^{1/2} & \text{si $v$ est archim\'edienne.}\end{cases}\end{equation*}Aux places archim{\'e}diennes le choix peut sembler arbitraire mais il assure la compatibilit{\'e} avec la litt{\'e}rature arakelovienne dans le cas o\`u $\overline{E}$ est hermitien. 
\begin{lemma}[Lemme 4.3.6 de~\cite{BGS}] Soit $\overline{E}$ un fibr{\'e} ad{\'e}lique \emph{hermitien} sur $\spec k$ et $\ell$ un entier $\ge 1$. L'isomorphisme canonique $\mathrm{H}^{0}(\mathbf{P}(E),\mathcal{O}(\ell))\overset{\sim}{\longrightarrow}S^{\ell}(E)$ induit une isom\'etrie sur les $\mathbf{C}_{v}$-espaces vectoriels correspondants si $v$ est ultram{\'e}trique et une similitude de rapport $\binom{n-1+\ell}{\ell}^{1/2}$ si $v$ est archim{\'e}dienne.
\end{lemma} 
Autrement dit, si $a_{n,\ell}$ d\'esigne l'id{\`e}le $(\binom{n-1+\ell}{\ell}^{1/2},\ldots,\binom{n-1+\ell}{\ell}^{1/2},1,\ldots,1,\ldots)$ alors $S^{\ell}(\overline{E})$ s'identifie \`a $a_{n,\ell}.\overline{E_{\ell}}$ (au sens de la page~\pageref{pagefibre}, avec $a_{n,\ell}$ vu comme une homoth\'etie de $E\otimes k_{\mathbf{A}}$). Ce lemme peut se d\'emontrer en explicitant la norme $\Vert s(x)\Vert_{v,x,\ell}$ avec des coordonn\'ees. Le choix d'une base $(X_{1},\ldots,X_{n})$ de $E_{1}\simeq E$ permet d'\'ecrire $s$ comme un polyn{\^o}me $P$ en les variables $X_{i}$ et si l'on d\'esigne par $(x_{1},\ldots,x_{n})\in\mathbf{C}_{v}^{n}$ les coordonn\'ees de $x$ dans cette base, on a la relation \begin{equation*}\Vert s(x)\Vert_{v,x,\ell}=\frac{\vert P(x_{1},\ldots,x_{n})\vert_{v}}{\Vert x\Vert_{E,v}^{\ell}}\ \cdotp\end{equation*}On peut alors calculer les int\'egrales ci-dessus lorsque $v$ est archim\'edienne et si la norme $\Vert.\Vert_{E,v}$ est hermitienne.
\begin{prop}
Pour tout fibr{\'e} vectoriel ad{\'e}lique $\overline{E}$ sur $\spec k$ et tout entier $\ell\ge 1$, la pente normalis{\'e}e du fibr{\'e} ad{\'e}lique hermitien $\overline{E_{\ell}}$ v{\'e}rifie \begin{equation*}\left\vert\widehat{\mu}_{\mathrm{n}}(\overline{E_{\ell}})-\ell\widehat{\mu}_{\mathrm{n}}(\overline{E})-\frac{\delta}{2}\log\left(\binom{n-1+\ell}{\ell}\gamma_{n,\ell}\right)\right\vert\le\frac{\ell}{D}\log\Delta(\overline{E})\ \cdotp\end{equation*}Quant \`a la pente maximale de $\overline{E_{\ell}}$, la quantit\'e \begin{equation*}\widehat{\mu}_{\mathrm{max}}(\overline{E_{\ell}})-\ell\widehat{\mu}_{\mathrm{max}}(\overline{E})-\frac{\delta}{2}\log\binom{n-1+\ell}{\ell}\end{equation*}est comprise entre $-\frac{\ell}{D}\log\Delta(\overline{E})$ et $\ell\left(2\delta n\log n+\frac{1}{D}\log\Delta(\overline{E})\right)$.
\end{prop}
\begin{proof}
Dans le cas hermitien, cette proposition est une simple cons{\'e}quence du th{\'e}or{\`e}me~\ref{pentefibresymetrique} et du lemme~\ref{corollaire64} car $\overline{E_{\ell}}=a_{n,\ell}^{-1}.S^{\ell}(\overline{E})$. Une fois ce cas {\'e}tabli, le cas g{\'e}n{\'e}ral s'obtient en consid{\'e}rant le fibr{\'e} $\overline{E}_{\varepsilon}$ (voir~\eqref{eqref:important}). Il faut observer que les normes sur $E_{\ell}$ donn\'ees par celles de $\overline{E_{\varepsilon}}$ et $\overline{E}$ v\'erifient des in\'egalit{\'e}s du m{\^e}me type que~\eqref{eqref:important}, \`a savoir que, pour tout sous espace vectoriel $E'$ de $E_{\ell}$, on a \begin{equation*}(E',(\vert.\vert_{(E_{\varepsilon})_{\ell},v})_{v})\preceq(E',(\Vert.\Vert_{E_{\ell},v})_{v})\preceq\left(E',\left(\left(\mathrm{d}(E\otimes_{k}k_{v},\ell^{2}_{n,k_{v}})(1+\varepsilon)\right)^{\ell}\vert.\vert_{(E_{\varepsilon})_{\ell},v}\right)_{v}\right)\cdotp\end{equation*}\end{proof}Ce dernier \'enonc\'e offre une transition vers la g\'eom\'etrie alg\'ebrique et la g\'eom\'etrie diophantienne, domaines dans lesquels le formalisme des pentes trouve probablement sa \emph{raison d'\^etre}. 
\section{Perspectives g\'eom\'etriques}
\label{par:conclusion}
Comme nous l'avons d{\'e}j{\`a} mentionn{\'e} dans l'introduction, le formalisme des pentes ad{\'e}liques a des applications en g{\'e}om{\'e}trie diophantienne. Dans ce dernier paragraphe, nous souhaitons donner quelques pistes pour mieux comprendre l'usage que l'on peut faire de ces pentes dans un probl{\`e}me diophantien de nature g{\'e}om{\'e}trique. Dans cette optique, l'in{\'e}galit{\'e} \begin{equation}\label{inegalitedepentesgeneralesdeux}\widehat{\mu}_{\mathrm{n}}(\overline{E})\le\sum_{i=1}^{N}{\frac{\dim(E_{i}/E_{i+1})}{n}\left\{\widehat{\mu}_{\mathrm{max}}(\overline{G_{i}})+h(\overline{E_{i}},\overline{G_{i}};\varphi_{i})\right\}}+\frac{1}{D}\log\vr(\overline{E})\end{equation}de la proposition~\ref{propositionmethodedespentes} joue un r{\^o}le important m{\^e}me si elle n'est finalement qu'une version savante de la formule du produit, ou de ce que l'on appelle en Transcendance l'\emph{in\'egalit\'e de Liouville}. Quoi qu'il en soit, elle est l'axe autour duquel s'articule la \emph{m\'ethode des pentes}, m\'ethode con{\c c}ue par J.-B.~Bost~\cite{bost2} et utilis\'ee par exemple dans les articles~\cite{artepredeux,graftieaux1,graftieaux2,viada}. Dans un contexte g{\'e}om{\'e}trique, les objets $\overline{E},\overline{G_{i}},\varphi$ sont souvent choisis de la mani{\`e}re suivante. On consid{\`e}re une vari{\'e}t{\'e} projective $X$ sur un corps global $k$ et $L\to X$ un fibr{\'e} en droites ample sur $X$. On pose $E:=\mathrm{H}^{0}(X,L^{\otimes m})$ l'espace vectoriel des sections globales d'une puissance enti\`ere $m$ de $L$. On se donne {\'e}galement un sous-sch\'ema ferm\'e $\Sigma$ de $X$, fini, compos\'e de points \'epaissis de $X(k)$ dans certaines directions \`a divers ordres. Un point de $\Sigma$ proc\`ede d'un triplet constitu\'e d'un point $k$-rationnel $x$ de $X$, d'une direction d'\'epaississement (p.~ex. cela peut {\^e}tre un sous-espace vectoriel de l'espace tangent $t_{X,x}$ \`a $X$ au point $x$) et d'un ordre de d\'erivation. L'objectif de la m\'ethode des pentes est de d\'egager quelques propri\'et\'es diophantiennes de $\Sigma$. Typiquement, on dispose d'un certain nombres de points complexes de $X$ dont on cherche \`a prouver la transcendance. On les suppose alg\'ebriques et on forme un sch\'ema $\Sigma$ avec cette hypoth\`ese. L'application lin{\'e}aire $\varphi$ est le morphisme d'{\'e}valuation qui \`a une section $s\in E$ associe la restriction $s_{\vert\Sigma}$ de $s$ \`a $\Sigma$. La filtration $(F_{i})_{i\in\{0,\ldots,N\}}$ de $F:=\mathrm{H}^{0}(\Sigma,L_{\mid\Sigma}^{\otimes m})$ est choisie de mani{\`e}re \`a annuler successivement les points de $\Sigma$ et les ordres de d{\'e}rivations correspondants. L'espace quotient $G_{i}$ s'identifie \`a un sous-espace de $S^{\ell}(t_{X,x}^{\mathsf{v}})\otimes x^{*}L^{\otimes m}$ o\`u $\ell\in\mathbf{N}$ et $x\in\Sigma$. Afin de mettre en {\oe}uvre l'in{\'e}galit{\'e} de pentes~\eqref{inegalitedepentesgeneralesdeux}, il faut choisir des structures ad{\'e}liques sur $E$ et les $G_{i}$, suffisamment appropri{\'e}es pour que la pente normalis{\'e}e $\widehat{\mu}_{\mathrm{n}}(\overline{E})$, la pente maximale $\widehat{\mu}_{\mathrm{max}}(\overline{t_{X,x}^{\mathsf{v}}})$ et le degr\'e normalis\'e $\widehat{\deg}_{\mathrm{n}}\overline{x^{*}L}$ puissent {\^e}tre {\'e}valu{\'e}s en fonction d'invariants attach{\'e}s \`a $X,L,x$. Nous allons d\'etailler cela dans un instant. Soulignons auparavant que, pour esp\'erer obtenir quelques informations sur $\Sigma$, il importe en g\'en\'eral d'avoir des donn\'ees rendues \emph{dynamiques} par l'introduction de param\`etres tels que l'entier $m$ dans la d\'efinition de $E$ ou bien ceux qui sont dissimul\'es dans la d\'efinition de $\Sigma$ (ordres d'annulations en chaque point $x$).
\par Revenons maintenant sur le choix des m\'etriques. Pour ce qui est de l'espace tangent $t_{X,x}$, on peut le munir d'une structure enti\`ere provenant d'un mod\`ele $(\EuScript{X}\to\spec\mathcal{O}_{k},\varepsilon_{x}:\spec\mathcal{O}_{k}\to\EuScript{X})$ de $(X,x)$, lisse en $\varepsilon_{x}$, o\`u $\mathcal{O}_{k}$ d{\'e}signe l'anneau des entiers de $k$. Aux places archim{\'e}diennes $v$ de $k$, la premi{\`e}re forme de Chern $c_{1}(L_{v})$ du fibr{\'e} ample $L_{v}\to X\times\spec\mathbf{C}_{v}$, induit par $L$, fournit des m{\'e}triques hermitiennes sur $t_{X,x}\otimes_{k}\mathbf{C}_{v}$. \par En ce qui concerne les structures ad{\'e}liques de $E$ et de $x^{*}L$, on peut distinguer deux {\'e}coles~: soit l'on choisit avec soin les m{\'e}triques pour en obtenir de \emph{canoniques} et l'on calcule alors explicitement les quantit{\'e}s $\widehat{\mu}_{\mathrm{n}}(\overline{E})$ et $\widehat{\deg}_{\mathrm{n}}\overline{x^{*}L}$ (cette derni{\`e}re {\'e}tant alors une \emph{hauteur canonique} de $x$); soit, \textit{a contrario}, on opte pour une tr{\`e}s grande souplesse dans le choix des m{\'e}triques, avec le minimum de contraintes, et l'on se contente de formules asymptotiques pour ces quantit{\'e}s avec $L$ remplac{\'e} par $L^{\otimes m}$ et $m\to+\infty$. D'une mani{\`e}re g{\'e}n{\'e}rale, la premi{\`e}re option est celle qui pr{\'e}domine en g{\'e}om{\'e}trie d'Arakelov. Mentionnons, \`a titre d'exemple, le cas d'une vari{\'e}t{\'e} ab{\'e}lienne $X$ sur un corps de nombres $k$, munie d'un fibr\'e en droites $L$ ample et sym\'etrique. En s'appuyant sur les travaux de L.~Moret-Bailly, J.-B.~Bost a montr{\'e} comment choisir un mod{\`e}le dit \emph{cubiste} de $(X,L,\Sigma)$ afin de calculer explicitement la pente normalis{\'e}e de $\overline{\mathrm{H}^{0}(X,L)}$ en termes de la hauteur de Faltings $h_{F}(X)$ de $X$ et du degr{\'e} g{\'e}om{\'e}trique $\deg_{L}X$; la formule exacte est $$\widehat{\mu}_{\mathrm{n}}(\overline{\mathrm{H}^{0}(X,L)})=-\frac{1}{2}h_{F}(X)+\frac{1}{4}\log\frac{\deg_{L}X}{(2\pi)^{d}d!}\qquad (d:=\dim_{k}X)\ .$$Dans ce cas le degr{\'e} normalis{\'e} $\widehat{\deg}_{\mathrm{n}}\overline{x^{*}L}$ est la hauteur de N{\'e}ron-Tate de $x$ relative \`a $L$ (voir le th{\'e}or{\`e}me~$5.10$ de~\cite{bost3}). Toutefois les choix des m{\'e}triques et les calculs qui en d{\'e}coulent sont assez d{\'e}licats, utilisant des r{\'e}sultats profonds comme le \emph{th{\'e}or{\`e}me de Riemann-Roch arithm{\'e}tique}. Aussi --- lorsque le probl{\`e}me s'y pr{\^e}te --- il est parfois pr{\'e}f{\'e}rable de se contenter de formules asymptotiques pour $\widehat{\mu}_{\mathrm{n}}(\overline{\mathrm{H}^{0}(X,L^{\otimes m})})$, lorsque $m\to+\infty$ (th{\'e}or{\`e}me de Hilbert-Samuel arithm{\'e}tico-g{\'e}om{\'e}trique). D'autant plus, peut-{\^e}tre, depuis la parution du m{\'e}moire de R.~Rumely, C.F.~Lau \& R.~Varley~\cite{Rumelyetal} qui fournit un {\'e}nonc{\'e} de ce type, sous des hypoth{\`e}ses tr{\`e}s faibles. Pour mieux comprendre de quoi il s'agit, nous allons {\'e}noncer le r{\'e}sultat principal de~\cite{Rumelyetal} (th\'eor\`eme~(A) p.~4) avec des hypoth{\`e}ses un peu plus fortes (p.~ex. nous prenons des normes au lieu de semi-normes). Ce r{\'e}sultat donne une illustration de la mani{\`e}re avec laquelle on proc{\`e}de pour donner une structure ad{\'e}lique \`a $\mathrm{H}^{0}(X,L^{\otimes m})$ et $x^{*}L$ et d'un calcul de pente normalis{\'e}e non trivial. \'Etant donn{\'e} une place $v$ de $k$ et $x\in X(\mathbf{C}_{v})$, on consid{\`e}re une norme $\Vert.\Vert_{L,v,x}$ sur la fibre $L_{x}\otimes_{k}\mathbf{C}_{v}$, invariante sous l'action de $\gal(\mathbf{C}_{v}/k_{v})$, norme que l'on transmet {\`a} $L_{x}^{\otimes m}\otimes_{k}\mathbf{C}_{v}$ par produit tensoriel. On suppose que, pour tout $x\in X(k)$, pour tout {\'e}l{\'e}ment $e\in L_{x}\setminus\{0\}$, on a $\Vert e\Vert_{L,v,x}=1$ pour toute place $v$ en dehors d'un nombre fini. Pour $s\in\mathrm{H}^{0}(X,L^{\otimes m})\otimes_{k}\mathbf{C}_{v}$, on pose alors $\Vert s\Vert_{v}:=\sup_{x\in X(\mathbf{C}_{v})}{\Vert s(x)\Vert_{L^{\otimes m},v,x}}$. Cela d{\'e}finit une structure ad{\'e}lique sur $\mathrm{H}^{0}(X,L^{\otimes m})$ d\`es lors que $\Vert s\Vert_{v}$ est toujours fini et m{\^e}me un peu plus, \`a savoir que si $s\ne 0$ et pour toute place $v$ en dehors d'un ensemble fini (qui peut d{\'e}pendre de $s$), on requiert $\Vert s\Vert_{v}=1$. L\`a encore, le choix des normes $\Vert.\Vert_{L,v,x}$ aux places ultram{\'e}triques $v$ de $k$ peut se faire au moyen d'un mod{\`e}le entier $(\EuScript{X},\EuScript{L})$ sur $\spec\mathcal{O}_{k}$ de $(X,L)$ (voir \textit{op. cit.}).
\begin{theoRumely}
Soit $X$ une vari{\'e}t{\'e} projective sur un corps global $k$. On suppose que $X$ est {\'e}quidimensionnelle et g{\'e}om{\'e}triquement r{\'e}duite. Soit $L\to X$ un fibr{\'e} en droites ample, muni de m{\'e}triques comme ci-dessus, conf{\'e}rant \`a $\mathrm{H}^{0}(X,L^{\otimes m})$ une structure de fibr{\'e} vectoriel ad{\'e}lique, pour tout entier $m\ge 1$. Alors il existe un {\'e}l{\'e}ment $h_{\overline{L}}(X)\in\mathbf{R}\cup\{+\infty\}$ tel que \begin{equation}\label{deficapacitesectionnelle}\frac{1}{m}\widehat{\mu}_{\mathrm{n}}(\overline{\mathrm{H}^{0}(X,L^{\otimes m})})\underset{m\to+\infty}{\longrightarrow} h_{\overline{L}}(X),\end{equation}{\'e}l{\'e}ment qui s'exprime en fonction de la \emph{capacit{\'e} sectionnelle} $S_{\gamma}(\overline{L})$ de $\overline{L}$ par la formule \begin{equation*}h_{\overline{L}}(X)=-\frac{\log S_{\gamma}(\overline{L})}{D(d+1)\deg_{L}X}\end{equation*}($d=\dim_{k}X$, $D=[k:k_{0}]$).
\end{theoRumely} 
Nous avons not{\'e} \`a dessein $h_{\overline{L}}(X)$ la limite ci-dessus car, lorsqu'elle est finie (et l'on conna{\^\i}t un crit\`ere pour que ce soit le cas, \textit{ibid}.), les propri\'et\'es de la capacit\'e sectionnelle font que cette quantit{\'e} se comporte effectivement comme une hauteur canonique de $X$ relativement \`a $\overline{L}$ (voir~\cite{aclens,Zhang2} ainsi que le th\'eor\`eme~(B) de~\cite{Rumelyetal}).
\par Lorsque $X$ est une vari{\'e}t{\'e} ab{\'e}lienne sur un corps de nombres et $\overline{L}$ un fibr{\'e} en droites cubiste, on a $h_{\overline{L}}(X)=0$. On peut observer que l'utilisation conjointe d'un r\'esultat asymptotique tel que~\eqref{deficapacitesectionnelle} et de l'in\'egalit\'e de pentes~\eqref{inegalitedepentesgeneralesdeux} n'est pas entrav\'ee par la pr\'esence du terme d'erreur $\frac{1}{D}\log\vr(\overline{\mathrm{H}^{0}(X,L^{\otimes m})})$, logarithmique en $m$ et donc n\'egligeable devant $m$.

\section*{Annexe}
Nous \'etablissons l'estimation asymptotique suivante, mentionn\'ee au paragraphe~\ref{subsectionpenteproduitsymetrique}.
\begin{propAnnexe}Soit $n,\ell\in\mathbf{N}\setminus\{0\}$ et $\gamma_{n,\ell}$ le nombre r\'eel positif d\'efini par la formule \begin{equation*}\log\gamma_{n,\ell}=\frac{1}{\binom{\ell+n-1}{n-1}}\sum_{\genfrac{}{}{0pt}{}{\mathbf{i}\in\mathbf{N}^{n}}{\vert\mathbf{i}\vert=\ell}}{\log\frac{\ell !}{\mathbf{i}!}}\ \cdotp\end{equation*}Posons $H_{n}:=\sum_{i=1}^{n}{1/i}$. Alors, lorsque $n$ est fix\'e et $\ell$ tend vers $+\infty$, on a \begin{equation*}\log\gamma_{n,\ell}=(H_{n}-1)(\ell+\mathrm{o}(\ell))\ .\end{equation*}
\end{propAnnexe}
\begin{proof} Le cas $n=1$ est imm\'ediat et l'on suppose $n\ge 2$ pour cette d\'emonstration. Notons $i_{1},\ldots,i_{n}$ les composantes de $\mathbf{i}\in\mathbf{N}^{n}$. Soit $f:[0,+\infty[\to\mathbf{R}$ la fonction continue d\'efinie par $f(0)=0$ et, si $x>0$, $f(x)=-x\log x$. En vertu de la formule de Stirling, rappel\'ee au bas de la page~\pageref{Stirling}, le terme $\log\frac{\ell !}{\mathbf{i}!}$ s'\'ecrit $\sum_{h=1}^{n}{f(i_{h}/\ell)}+\mathrm{O}(\log\ell)$. On a alors \begin{equation*}\frac{\log\gamma_{n,\ell}}{\ell(n-1)!}\underset{\ell\to\infty}{\sim}\frac{1}{\ell^{n-1}}\sum_{\vert\mathbf{i}\vert=\ell}{\sum_{h=1}^{n}{f\left(\frac{i_{h}}{\ell}\right)}}\ \cdotp\end{equation*}On reconna{\^\i}t pour le membre de droite la somme\begin{equation*}\frac{1}{\ell^{n-1}}\sum_{\genfrac{}{}{0pt}{}{\mathbf{j}\in\mathbf{N}^{n-1}}{\vert\mathbf{j}\vert\le\ell}}{\left\{f\left(\frac{j_{1}}{\ell}\right)+\cdots+f\left(\frac{j_{n-1}}{\ell}\right)+f\left(1-\sum_{h=1}^{n-1}{\frac{j_{h}}{\ell}}\right)\right\}},
\end{equation*}qui est une somme de Riemann $(n-1)$-dimensionnelle convergeant vers l'int\'egrale \begin{equation}\label{cinquantecinq}\int_{\Omega}{\left(f(x_{1})+\cdots+f(x_{n-1})+f(1-x_{1}-\cdots-x_{n-1})\right)\,\mathrm{d}x_{1}\ldots\mathrm{d}x_{n-1}}\end{equation} sur l'ouvert $\Omega:=\{(x_{1},\ldots,x_{n-1})\in\,]0,+\infty[^{n-1};\ x_{1}+\cdots+x_{n-1}<1\}$. Par sym\'etrie, les int\'egrales $\int_{\Omega}{f(x_{i})}$, $1\le i\le n-1$, sont toutes \'egales. Par changement de variable $y:=1-x_{1}-\cdots-x_{n-1}$, qui ne change pas le domaine $\Omega$, la derni\`ere int\'egrande de~\eqref{cinquantecinq} peut {\^e}tre remplac\'ee par $f(x_{n-1})$, et l'int\'egrale~\eqref{cinquantecinq} vaut alors \begin{equation*}\begin{split}n\int_{\Omega}{f(x_{n-1})\,\mathrm{d}x_{1}\ldots\mathrm{d}x_{n-1}}&=n\int_{0}^{1}{f(u)\vol(x_{1}+\cdots+x_{n-2}\le 1-u)\,\mathrm{d}u}\\ &=\frac{n\vol(b_{n-2}^{1})}{2^{n-2}}\int_{0}^{1}{(1-u)^{n-2}f(u)\,\mathrm{d}u}\end{split}\end{equation*}(ici $\vol$ est la mesure de Lebesgue sur $\mathbf{R}^{n-2}$). La formule~\eqref{formuledevolumep} donne $n\vol(b_{n-2}^{1})/2^{n-2}=n/(n-2)!$. Et la derni\`ere int\'egrale se calcule au moyen de deux int\'egrations par partie (en observant que $\int_{0}^{1}{\frac{1-u^{n}}{1-u}\,\mathrm{d}u}=H_{n}$). On trouve $\int_{0}^{1}{(1-u)^{n-2}f(u)\,\mathrm{d}u}=\frac{1}{n(n-1)}(H_{n}-1)$, ce qui permet de conclure.
\end{proof}
 Cette d\'emonstration rev{\^e}t un caract\`ere \emph{ad hoc}. Avec une analyse du terme reste de la somme de Riemann en fonction de la premi\`ere d\'eriv\'ee de $f$, on peut montrer que $\log\gamma_{n,\ell}=(H_{n}-1)\ell+\mathrm{O}(\log\ell)$. Une vision plus savante et inspir\'ee de ce genre de calculs se trouve dans la th\`ese de H.~Randriambololona, dans laquelle le lecteur pourra trouver une \'evaluation asymptotique \`a un ordre quelconque de (variantes de) $\log\gamma_{n,\ell}$ et d'autres quantit\'es d'origine g\'eom\'etrique, plus g\'en\'erales (voir~\cite{theseRandriam}, chapitre~$4$, proposition 4.2.3 et \textit{supra}).  

\bibliographystyle{plain}

\noindent Universit\'e Grenoble I, Institut Fourier.\\ UMR $5582$, BP $74$\\ $38402$ Saint-Martin-d'H{\`e}res Cedex, France.\\ 
Courriel~: \texttt{Eric.Gaudron@ujf-grenoble.fr}\\ 
Page internet~: \texttt{http://www-fourier.ujf-grenoble.fr/\~{}gaudron}

\end{document}